\documentclass{amsart}
  \usepackage{amscd,amssymb,epsfig}
   \usepackage{epic,eepic}

                     \numberwithin{equation}{subsection}

                     \newtheorem{propo}{Proposition}[subsection]
                     \newtheorem{corol}[propo]{Corollary}
                     \newtheorem{theor}[propo]{Theorem}
                     \newtheorem{lemma}[propo]{Lemma}
                     \theoremstyle{definition}

                     \theoremstyle{remark}
                     \newtheorem{remar}[propo]{Remark}
\newtheorem{remars}[propo]{Remarks}

		     \newcommand{\CC}{\mathbb{C}}
		     \newcommand{\QQ}{\mathbb{Q}}
                     \newcommand{\ZZ}{\mathbb{Z}}
                     \newcommand{\RR}{\mathbb{R}}

		     \newcommand{\A}{\mathcal{A}}
 \newcommand{\M}{\mathcal{M}}
		     \newcommand{\str}{\mathcal{S}}

                     \newcommand{\Hom}{\operatorname{Hom}}
		     \newcommand{\Ker}{\operatorname{Ker}}
		     \newcommand{\Exp}{\operatorname{Exp}}
\newcommand{\Lbl}{\operatorname{Lbl}}
\newcommand{\lbl}{\operatorname{lbl}}
\newcommand{\rk}{\operatorname{rank}}
\newcommand{\Ann}{\operatorname{Ann}}

		    \newcommand{\Ima}{\operatorname{Im}}
		      \newcommand{\inc}{\operatorname{in}}
		      \newcommand{\Fix}{\operatorname{Fix}}
                     \newcommand{\id}{\operatorname{id}}

 \newcommand{\rrr}{r}

		     \newcommand{\sign} {\operatorname {sign}}
		    \newcommand{\modu}{\operatorname{mod}}
\newcommand{\Int}{\operatorname{Int}}
		     \newcommand{\arr}{\operatorname{arr}}
		  \newcommand{\rank}{\operatorname{rank}}
 \newcommand{\edg}{\operatorname{edg}}
		  \newcommand{\Perm}{\operatorname{Perm}}

\begin{document}
      \title{Virtual strings and their cobordisms}
                     \author[Vladimir Turaev]{Vladimir Turaev}
                     \address{%
              IRMA, Universit\'e Louis  Pasteur - C.N.R.S., \newline
\indent  7 rue Ren\'e Descartes \newline
                     \indent F-67084 Strasbourg \newline
                     \indent France }
                     \begin{abstract} 	A virtual string    is   a scheme of  self-intersections of a closed
curve on  a  surface.   We  study algebraic
invariants of strings as well as  two   equivalence relations on the set of strings:  homotopy and cobordism. 
We show that the homotopy invariants of strings  form an infinite dimensional  Lie   group.
 We also discuss connections  between virtual strings  and virtual
knots. 
                     \end{abstract}
                     \maketitle

\centerline {\bf Contents}
 \vskip1truecm

\noindent  {\bf  1. Introduction}

\noindent {\bf 2.  Generalities on virtual strings}

\noindent  {\bf 3.  Polynomial $u$}

\noindent {\bf 4.  Geometric realization  of virtual strings}

\noindent {\bf 5.  Cobordism  of strings and the slice genus}

\noindent  {\bf 6.  Based   matrices  of  strings}

\noindent {\bf 7.  Genus and cobordism for  based   matrices}

\noindent {\bf 8.  Genus estimates and sliceness of strings}

\noindent {\bf 9.  Lie cobracket for strings}

\noindent {\bf 10.  Virtual strings versus  virtual  knots}

\noindent {\bf 11.  Proof of Theorem \ref{th:eee}}

\noindent {\bf 12.      Open strings}

\noindent {\bf 13.   Questions}

                  \section{Introduction} 
		  
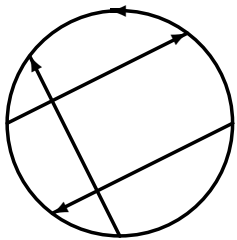
\begin{figure}[ht]
                     \setlength{\unitlength}{0.5cm}
                     \begin{picture}(18,6)(-10,-2.5)
                             \allinethickness{.5mm}
                             \put(0,0){\circle{6}}
                             \put(0,3){\vector(-1,0){.25}}
                             \put(-3,0){\vector(2,1){4.8}}
                             \put(0,-3){\vector(-1,2){2.4}}
                             \put(3,0){\vector(-2,-1){4.8}}
                     \end{picture}
                     \caption{A virtual string of rank 3}\label{fg:gdg}
                     \end{figure}

	A virtual string    is   a scheme of  self-intersections of a generic oriented closed
curve on  an 
oriented surface. More precisely, a    virtual string of rank $m\geq 0$   
is  an 
oriented circle  with $2m$ distinguished points  partitioned into $m$ ordered pairs. These 
$m$ ordered pairs   of points  are 
called   arrows of the virtual string. An example of a virtual string of rank 
$3$  is 
shown on Figure~\ref{fg:gdg} where the  arrows are represented by 
geometric 
vectors.
 
A  (generic  oriented)  closed curve  on an oriented 
surface   gives rise to an \lq\lq underlying" virtual string whose arrows correspond to the self-crossings of the curve.  The
usual homotopy of curves  on surfaces   suggests a  notion of homotopy for   strings. 
The
  homotopy of curves  in 3-manifolds   with boundary  suggests a  notion of  cobordism for   strings.  The main objective  of
the theory of virtual strings is  a study (and  eventually classification) of their homotopy classes and cobordism classes. To
this end, we   introduce     algebraic  invariants of virtual 
strings, specifically,  a one-variable polynomial  $u$ and a so-called
based matrix.   We
formulate  obstructions to  homotopy/cobordism of strings in terms of these invariants. This leads us to a purely algebraic
study of  analogues of homotopy and cobordism for skew-symmetric matrices. 

As an instance of cobordism, we call  a string  {\it  slice} if it can be realized by  a closed curve on the boundary of an
orientable 3-manifold $M$ that  is   contractible in $M$. We formulate  obstructions to the sliceness of  a string in
terms of the polynomial $u$ and the based matrix.  
 
We introduce a natural Lie cobracket in  the free abelian group generated by the homotopy classes of strings. Dually, the
abelian group  of   $\ZZ$-valued homotopy invariants of strings  becomes
 a Lie algebra.  This Lie algebra is integrated into  an infinite dimensional  Lie   group.  This  Lie group gives rise to  further 
algebraic objects including a Hopf algebra structure on the (commutative) polynomial algebra generated by the homotopy
classes of strings.

 Virtual strings are   closely related to   virtual knots   introduced by L.  Kauffman \cite{Ka}.  In particular, the  term  \lq\lq virtual knots" 
suggested to us the term virtual strings.  
		 Virtual knots  can be defined as  equivalence classes of  arrow diagrams which are just   virtual strings whose arrows  are 
provided with signs $+$ or $-$.  Forgetting these signs, we obtain  a map from the set of virtual knots into the
set of  homotopy classes of virtual strings. We give a more elaborate construction which associates with each virtual knot  a
polynomial  expression in virtual strings with coefficients in the ring
$\QQ[z]$. This leads to  an isomorphism between  a  \lq\lq skein algebra" of virtual knots   and  a polynomial
algebra generated by the homotopy classes of strings. 

We also study \lq\lq open strings" which are schemes of self-intersections of generic paths on surfaces.

As an application of this work, we obtain a new interesting relation of cobordism for skew-symmetric matrices and a homeomorphism  invariant 
of knots in cylinders  over   oriented surfaces with values in the polynomial ring $\QQ[z,t]$.

 A number of ideas and results     of  this paper have  predecessors in the
literature.  The homotopy  for  
  virtual strings can be defined in terms of the stable
equivalence    of curves  on surfaces    introduced by  J. S. Carter, S. Kamada,  and M. Saito
\cite{cks}. These authors also defined the notions of cobordism and sliceness for curves on surfaces which is essentially
equivalent to our cobordism and sliceness for strings.     Carter
\cite{ca2}   showed that there are closed curves on closed   oriented  surfaces that bound no  singular  disks in
3-manifolds bounded by these surfaces.  In  our language this means that there are non-slice strings
(see also \cite{hk}). 

   The main novelty of the present paper  lies       in the  
introduction of new homotopy invariants of strings and new operations on the  homotopy classes of strings.

 The organization of the paper should be clear from the Contents above.
	
		       \section{Generalities on virtual strings}
                    
		         \subsection{Definitions}\label{sn:g11} We give 
here a formal definition of a virtual string.   For an  integer $m\geq 
0$,   a 
{\it  virtual string $\alpha$ of rank $m$} (or briefly a {\it 
string}) 
is  an oriented circle, $S$, called the {\it core circle} of $\alpha$,  and a 
distinguished   set of $2m$ distinct points of $S$ partitioned into $m$  
ordered 
pairs.    We   call these $m$ ordered     pairs   of points  the {\it 
arrows} of 
$\alpha$. The set of arrows of $\alpha$ is denoted $\arr(\alpha)$. The 
endpoints 
$a,b\in S$ of  an arrow  $(a, b)\in \arr(\alpha)$  are called   its 
{\it  tail} 
and {\it   head}, respectively. The $2m$ distinguished 	points of $S$ 
are called 
the {\it endpoints} of $\alpha$.	
			 
 The  string formed by an oriented circle and an empty set of 
arrows is 
called a {\it trivial virtual string}.  
  An example of a virtual string of rank $3$  is shown on 
Figure~\ref{fg:gdg}.

			       By a {\it homeomorphism} of two virtual strings, 
we mean  an orientation-preserving  homeomorphism of the core circles 
transforming the  set of  arrows of the  first   string   onto the  set 
of arrows 
of the second string. Two virtual strings are {\it homeomorphic} if 
they are 
related  by   a homeomorphism. Clearly, homeomorphic strings have the same rank. 

By abuse of language, the homeomorphism classes of       virtual 
strings will be 
also called virtual strings.

		      \subsection{From curves   to    strings}\label{sn:g12} 
		     By a   surface,  we mean a smooth {\it oriented} 
2-dimensional manifold.   By a    {\it   closed curve} on a 
surface 
$\Sigma$, we mean a generic  smooth immersion  $\omega$ of an oriented circle $S$ into 
$\Sigma$. Recall that   a smooth map $ S\to \Sigma$ is  an {\it immersion} if   
its 
differential is  non-zero at all points of $S$.  An immersion 
$\omega:S\to \Sigma$ 
is {\it generic} if  $ \# (\omega^{-1}(x))\leq 2$  for all  $ x\in \Sigma$,  the 
set $\{x\in \Sigma\,\vert\, \# (\omega^{-1} (x))=2\}$ is finite, and   all 
its  points 
  are transverse intersections of   two branches.  Here and below the symbol 
$\# (A)$ 
denotes the cardinality of a set $A$. The points $x\in \Sigma$ such 
that $ \# 
(\omega^{-1} (x))=2 $ are called {\it double points} or {\it crossings} of 
$\omega$.  
		     
		     A closed curve $\omega:S\to \Sigma$   gives rise to an {\it underlying  virtual 
string} $\alpha_\omega$. The core circle of $\alpha_\omega$  is $S$ and the   arrows  of $\alpha_\omega$ 
are 
all ordered pairs $a,b\in S$ such that $\omega(a)=\omega(b)$ and the pair (a 
positive 
tangent vector of $\omega$ at $a$, a positive tangent vector of $\omega$ at 
$b$)  is a 
positive basis in   the tangent space of   $\omega(a)$.  For instance,  the underlying string of a simple closed curve on
$\Sigma$ is a trivial virtual string.

We say that a  virtual string  
 is {\it realized} by a closed curve $\omega:S\to \Sigma$ if  it is homeomorphic to $\alpha_\omega$. As we 
shall see 
below, every virtual string can be realized by a closed curve  on a 
surface.
		    
		     \subsection{Homotopy of   strings}\label{sn:g13}    The 
usual 
homotopy of closed curves on a surface suggests  to introduce a relation of homotopy for   virtual strings. Observe first   
that   two 
homotopic   curves on a surface can be related by a finite sequence of 
the 
following \lq\lq elementary" moves     (and the inverse 
moves):

		     (a) a local move adding a small curl to the curve;
		     
		     (b) a local move pushing a branch of the curve across 
another branch and creating  two new double points;
		     
		     (c) a local   move pushing a branch of the curve across a 
double point;
		     
		     (d) ambient isotopy in the surface.

	The move (a) has two forms (a${)}^{+}$   and (a${)}^{-}$   depending on	whether the curl lies on 
the left or   the right of  the curve where the left and the right are determined by the direction of the curve and the
orientation of the surface.  Considered up to ambient isotopy, the move (b) has three forms  depending on the direction 
of the two branches.   Similarly,  considered up to ambient isotopy, the move (c) has two
forms  (c${)}^{+}$   and (c${)}^{-}$  
depending on the direction of the   branches. Using the standard braid generators $\sigma_1, \sigma_2$ on 3 strands we
can encode  this move  as
$\sigma_1\sigma_2 
\sigma_1\mapsto \sigma_2\sigma_1 \sigma_2$ where the over/undercrossing 
information is forgotten.  The  moves  (c${)}^{+}$   and (c${)}^{-}$    are obtained by directing (before
and after the move)  the first and third strands up and the second strand up or down, respectively.   It is easy to see that  
(c${)}^{+}$,  (c${)}^{-}$  can be obtained from each other  using    ambient isotopy,    moves  (b), and  inverses to
 (b).  Similarly, the moves (a${)}^{+}$,  (a${)}^{-}$  can be obtained from each other  using    ambient isotopy,   
moves  (b),  (c${)}^{-}$, and  inverses to   (b).  Thus the moves (a${)}^{-}$, (b), (c${)}^{-}$ generate all the other moves.

It is clear that  ambient isotopy of a closed curve does not change 
the 
underlying virtual string. We now describe the analogues  for virtual 
strings of 
the moves (a${)}^{-}$, (b), (c${)}^{-}$. In this description and in the sequel, by an 
{\it arc} 
on  an oriented  circle $S$ we mean an {\it embedded arc} on $S$. The orientation 
of $S$ 
induces an orientation of all arcs on $S$. For two distinct points 
$a,b\in S$, we 
write   $ab$ for the unique oriented arc  in $S$ which      begins in 
$a$ and 
terminates in $b$.  Clearly, $S= ab \cup ba$ and $ab\cap ba= \{a,b\}$.

Let $\alpha$ be a virtual string with core circle $S$. Pick two 
distinct points 
$a,b \in S$  such that the arc   $ab\subset S$ is disjoint from the set of   
endpoints of 
$\alpha$. The move (a$)_s$, where $s$ stands for \lq\lq string",   
adds to 
$\alpha$  the pair $(a,b)$.    This amounts 
to 
attaching a small arrow to $S$  such that  the arc  in $S$ leading from its tail to its head  is disjoint from the   endpoints of
$\alpha$. The move (b$)_s$  acts on 
$\alpha$ as follows. Pick two   arcs   on $S$ disjoint from each 
other 
and   from the   endpoints of $\alpha$. Let $a,a'$ be the 
endpoints of 
the first arc (in an arbitrary order) and $b,b'$ be the endpoints of 
the second 
arc. The move   (b$)_s$ adds 
to $\alpha$ two arrows $(a,b)$ and $(b',a')$.  (This move has four forms depending on the
  two possible choices for  $a$  and two possible choices for  $b$.   However, two of these
 forms  of (b$)_s$  are  equivalent.) The move (c$)_s$     applies  to 
$\alpha$   when $\alpha$ has three arrows 
$(a^+,b), (b^+,c), (c^+,a)$ where  $a, a^+, b,b^+,c, c^+\in S$ such that
the   arcs   $aa^+$,   $bb^+$,     $cc^+$ are disjoint from each 
other and 
from the other   endpoints of $\alpha$. The   move (c$)_s$  replaces 
the  arrows  $(a^+,b), (b^+,c), (c^+,a)$ with the arrows  $ (a,b^+), (b,c^+), (c,a^+)$.

	We say that two virtual strings are {\it homotopic} if they can be 
related by a finite sequence of homeomorphisms, the  {\it homotopy moves} (a$)_s$, 
(b$)_s$, 
(c$)_s$, and the inverse moves. A virtual string homotopic to a trivial 
virtual 
string is said to be {\it homotopically trivial}.  For instance, as it follows directly from the definitions,  all virtual strings of rank 
$\leq 2$ 
are homotopically trivial. 

It is clear from what 
was said 
above that the underlying virtual strings of homotopic closed curves on 
a 
surface are themselves homotopic.  
	
\subsection{Transformations of  
strings}\label{sn:g1335}   For a   string $\alpha$, we define the {\it  opposite string}  $  
\alpha^-$ to be $\alpha$ with opposite orientation on 
the core 
circle. The {\it inverse  string} $ \overline \alpha$ is obtained from 
$\alpha$ 
by reversing all its arrows. On the level of closed curves on surfaces, 
these 
two transformations correspond to traversing the same curve in the 
opposite 
direction and to inverting the orientation of the ambient surface, respectively.
If two strings are homotopic, then their opposite (resp. inverse) 
strings are 
homotopic.  

One can raise a number of  questions concerning the transformations $\alpha\mapsto \alpha^-$,  $\alpha\mapsto \overline
\alpha$, $ \alpha\mapsto \overline
\alpha^-$. For instance, one can ask 
  whether  there is a  string $\alpha$ that is not homotopic to $\alpha^-$ (resp.  to $ \overline
\alpha$,  $\overline
\alpha^-$).    Below we will answer this
question   in the positive. 

A virtual string $\alpha$ with core circle $S$ is a {\it product} 
of 
virtual strings $\alpha_1$ and $\alpha_2$ if there are   disjoint 
arcs 
$a_1b_1, a_2b_2\subset S$ such that each arrow of $\alpha$ has both 
endpoints on either $a_1b_1$ or on $a_2b_2$ and the   string formed by 
$S$ and 
the arrows of $\alpha$ with endpoints on $ a_i b_i$ is homeomorphic to 
$\alpha_i$ for $i=1,2$.   
One can ask whether the product is  a  well-defined operation on strings (at least up to homotopy) and whether it is
commutative. Below we will answer these
questions  in the negative.

\subsection{Geometric invariants   of  
strings}\label{sn:g133} 
 We define three geometric characteristics of  strings: the genus, the 
homotopy 
genus,    and the homotopy rank.    The {\it genus} 
$g(\alpha)$ of a string
$\alpha$  is  the minimal integer $g\geq 0$   such that $\alpha$ can 
be 
realized by a closed curve on a surface of genus $g$. The {\it homotopy 
genus} 
$hg(\alpha)$  is     the minimal integer $g\geq 0$   such that $\alpha$ 
is 
homotopic to a string of genus $g$.      The {\it homotopy rank} 
$hr(\alpha)$  
is  the minimal integer $m\geq 0$   such that $\alpha$ is homotopic to 
a string 
of rank $m$. For example, if $\alpha$ is a trivial string, then 
$g(\alpha)=hg(\alpha)=hr(\alpha)=0$.  
 It is clear that  the homotopy genus and the homotopy rank are 
homotopy 
invariants of strings. Below we  compute the genus  explicitly  and show that it is not a 
homotopy 
invariant.

The numbers $g(\alpha), hg(\alpha), hr(\alpha)$  are  preserved under the transformations $\alpha \mapsto \alpha^-, 
\alpha
\mapsto 
\overline
\alpha$.

 \subsection{Encoding of strings}\label{sn:g1345} There are two 
simple 
methods allowing to encode virtual strings in a compact way. Although 
we do 
not use these methods in this paper, we briefly describe them   for 
completeness.
 
(1)  Consider a 
finite set 
$E$ consisting of $m$ elements and its disjoint copy $E^+=\{x^+\,\vert 
\, x\in 
E\}$. 
Let $y_1,y_2,\ldots, y_{2m}$  be  a sequence of elements of the set 
$E\cup E^+$ 
in which every element  appears  exactly once. (Such a sequence 
determines a 
total order in  $E\cup E^+$ and vice versa.)  The sequence  $ y_1,y_2,\ldots, y_{2m}$ defines  a 
string  
of rank $m$ whose  underlying circle   is    $S=\RR\cup \{\infty\}$ 
with    
right-handed orientation on $\RR$ and whose  $m$   arrows are the pairs 
$(a,b)$ 
such that $a,b\in \{1,2,\ldots, 2m\}\subset S, y_a\in E$, and 
$y_b=y_a^+\in E^+$.  Any string can be encoded in this way. For instance, the 
string drawn 
in Figure \ref{fg:gdg} is encoded by the sequence
$x^+,y,z^+,x,z,y^+$ where $E=\{x,y,z\}$.

(2) By Section \ref{sn:g12},  virtual strings can 
be 
encoded by  closed curves on surfaces. This   has an extension similar 
to 
Kauffman's  graphical  encoding of virtual knots in \cite{Ka}. Namely, 
consider 
a (generic) closed 
curve   on a surface and suppose that some of its crossings are marked 
as \lq\lq 
virtual". Take the   string of this curve as in Section 
\ref{sn:g12} 
and forget all its arrows corresponding to virtual crossings. 
 It is easy to see that every virtual string can be obtained in this 
way from a  
 closed curve  in $\RR^2$ with virtual crossings. This  yields a  
graphical 
encoding of strings by plane curves with virtual crossings. 
The relation of  homotopy for strings has a simple description in  this 
language: it is generated by  the moves shown in \cite{Ka}, 
Figure 2
(where the over/undercrossing  information should be forgotten).

 \subsection{Remarks}\label{sn:g14}   1.    We can point out certain classes of closed  curves on surfaces whose 
underlying  virtual strings are homotopically trivial.   Since  all closed curves on  
$ S^2$ are contractible,  their underlying   strings are  homotopically trivial.  Therefore the same is true for 
closed  curves on any subsurface of $S^2$, i.e.,   on any surface of genus 0.  In particular, all  closed  curves  on an
annulus have   homotopically trivial underlying  strings.  Since each  closed curve  on a torus can be deformed into an
annulus,  its  underlying string   is homotopically trivial.  The
same holds for    closed curves on a torus with holes. 
            
 2.  The move (a$)_s$     has a   version  (a${)}^{+}_s$ which 
is defined   as (a$)_s$ above but adds the arrow $(b,a)$ rather than $(a,b)$.  
 This  move 
underlies 
the move (a${)}^{+}$ on closed curves.  The
move  (a${)}^{+}_s$ preserves the homotopy class of  a string.  Indeed, it     can be 
expressed  as a composition of   (b$)_s$,  (c$)_s$, and an inverse  to   (a$)_s$.

3. The move (c$)_s$     has a  version  (c${)}^{+}_s$ which 
applies  to a string    
 when it has three arrows 
$(a,b), (a^+,c)$, $(b^+,c^+)$ 
such that
the   arcs   $aa^+$,   $bb^+$,    $cc^+$ are disjoint from each 
other and 
from the other   endpoints of the string. The   move   (c${)}^{+}_s$  replaces these 
three arrows 
with the arrows  $ (a^+,b^+), (a,c^+), (b,c)$.
 This  move 
underlies 
the move (c${)}^{+}$ on closed curves.  The   move   (c${)}^{+}_s$  can be 
expressed
as a composition of   (c${)}^{-}_s$    and  (b$)_s$.

            \section{Polynomial $u$}\label{abric} 
                    
		         \subsection{Invariants $\{u_k\}_k$}\label{sn:g21} Let 
$\alpha$ be a virtual string with core circle $S$. Each arrow 
$e=(a,b)\in 
\arr(\alpha)$ splits $S$ into two arcs   $ {ab}$ and $  {ba}$. We say 
that an 
arrow  $f=(c,d)$  of $\alpha$ (distinct from $e$) {\it links} $e$ if one of its endpoints lies on 
$ab$ and the 
other one lies on $ba$. More precisely,     $f=(c,d)$ 
  links $e$ {\it positively} (resp. {\it negatively}) if  $c\in 
ab, d\in 
ba$ (respectively, if $ c\in ba, d\in ab$).  If $f$ does not link $e$, then  $e$ and $f$ are {\it unlinked}.  Let
$n(e)\in
\ZZ$ be the algebraic  number of arrows of
$\alpha$  linking
$e$, i.e., 
  the 
number of arrows of $\alpha$  linking $e$ positively minus   the 
number of arrows of $\alpha$  linking $e$
negatively.
			 
			 It is easy to trace the behaviour of   $n(e)$ under the 
homotopy moves (a$)_s$, (b$)_s$,  (c$)_s$ on $\alpha$.
			 The move (a$)_s$ adds an arrow $e_0$ with  $ n( e_0)=0$ 
and keeps $n(e)$ for all other arrows. The move (b$)_s$ adds two arrows 
$e_1,e_2$ with $n(e_1)=-n( e_2)$ and keeps $n(e)$ for all other arrows. 
Consider 
  the move (c$)_s$ and use the notation of Section \ref{sn:g13}. It is 
obvious 
that for all arrows $e$ preserved under the move, the number $n(e)$ is 
also 
preserved.  
Each arrow $e=(a^+,b), (b^+,c), (c^+,a)$ occuring before the move gives 
rise to 
an arrow			 $ e'= (a,b^+), (b,c^+), (c,a^+)$, respectively, 
occuring after the move. We claim that $n(e)	=n(e')$.	Consider for 
concreteness $e=(a^+,b)$. Note that  the points $c,c^+$ lie either on 
$ab$ or on 
$ba$. Suppose that  $c,c^+\in ab$. Then  the arrows	$(b^+,c)$ and $ 
(c^+,a)$ 
contribute $1$ and $-1$ to $n(e)$, respectively, while the  corresponding 
arrows 
$(b,c^+)$ and $ ( c, a^+)$		contribute  $0$ to $n(e') $.
All other arrows contribute the same to  	$n(e)$ and  $n(e') $. 
Hence 
$n(e)=n( e')$. 
If $c,c^+\in ba$, then  the arrows	$(b^+,c)$ and $ (c^+,a)$ 
contribute $0$ to $n(e)$ while the corresponding arrows 
$(b,c^+)$ and $ ( c, a^+)$ 	contribute $-1$ and $1$ to $n( e')$, respectively.
All other arrows contribute the same to  	$n(e)$ and  $n( e')$. 
Hence 
$n(e)=n( e')$.  

For an integer $k\geq 1$,  set  $$u_k(\alpha)=\#\{e\in 
\arr(\alpha)\,\vert \,n(e) =k\}
- \#\{e\in \arr(\alpha)\,\vert \,n(e) =- k\} \in \ZZ.$$
  It is clear from what was said above that $u_k(\alpha)$ is preserved 
under the 
moves  (a$)_s$, (b$)_s$,  (c$)_s$. In other words, $u_k(\alpha)$ is a homotopy 
invariant 
of $\alpha$.  Clearly, $u_k(\alpha)=0$ for all $k$ greater than or equal to the rank of $\alpha$.  If $\alpha$ is homotopically
trivial,  then 
$u_k(\alpha)=0$ for all $k\geq 1$.

\subsection{Polynomial  $u(\alpha)$}\label{sn:g22} We can combine the     
invariants 
$u_k$ of a virtual string $\alpha$ into a polynomial 
$$u(\alpha)=\sum_{k\geq 1} u_k(\alpha)\, t^k$$
where $t$ is a  variable. The free term of this polynomial is always 
$0$ and its 
degree is bounded from above by $m-1$ where $m$ is the rank of 
$\alpha$. This 
polynomial is a homotopy invariant of $\alpha$. If $\alpha$ is 
homotopically 
trivial, then $u(\alpha)=0$. (The converse is not true, as we shall see 
below.) 
The polynomial $u(\alpha)$ yields an  estimate  for the 
homotopy rank 
$hr(\alpha)$ of $\alpha$ defined in Section \ref{sn:g133}: 
\begin {equation}\label{homr}
hr(\alpha)\geq \deg u(\alpha)+1.\end{equation}

We can rewrite  $u(\alpha)$ as follows:
\begin {equation}\label{homr11} u(\alpha)=\sum_{e\in \arr(\alpha), n(e)\neq 0} \sign (n(e))\, 
t^{\vert 
n(e)\vert}\end{equation}
where $\sign (n)=1$ for positive $n\in \ZZ$ and $\sign (n)=-1$ for 
negative 
$n\in \ZZ$.
 Therefore
$$\sum_{k\geq 1} k\, u_k(\alpha) \,t^{k-1}=u'(\alpha)=\sum_{e\in \arr(\alpha), n(e)\neq 0}  
n(e)\, 
t^{\vert n(e)\vert-1}=	\sum_{e\in \arr(\alpha)}  n(e)\, t^{\vert n(e)\vert-1}. $$
Substituting $t=1$, we obtain 
$$\sum_{k\geq 1} k\, u_k(\alpha)  =u'(1)=\sum_{e\in \arr(\alpha)}   n(e)=0.$$
The last equality follows from the fact that 	if an arrow $f$ links an 
arrow 
$e$ positively, then $e$ links $f$ negatively.	
	
		       \subsection{Examples}\label{sn:g23}  1.   For positive 
integers $p,q$, we define  $\alpha_{p,q}$ to be the lattice-looking 
virtual 
string formed by a   Euclidean circle in $\RR^2$ with counterclockwise 
orientation, $p$ 
disjoint vertical arrows $e_1,\ldots,e_p$ directed upward and numerated 
from 
left to right, and $q$ disjoint horizontal arrows   $e_{p+1},\ldots, 
e_{p+q}$  
crossing $e_1,\ldots,e_p$ from   right to left and numerated from 
bottom to top. 
(Here we identify arrows with geometric vectors in $\RR^2$  connecting 
two 
points of the core circle; the numeration of the arrows is compatible 
with the 
counterclockwise order of their tails.) Clearly, $n(e_i)= q$  for 
$i=1,\ldots , 
p$ and   
$n(e_{p+j})= -p$ for $j=1,\ldots , q$. Hence $u(\alpha_{p,q})=p t^{q}- 
q t^p$.  
We   conclude that the   strings $\{\alpha_{p,q}\}_{p\neq q}$ are 
pairwise 
non-homotopic and homotopically non-trivial.  The string    
$\alpha_{1,1}$ is 
homotopically trivial: it is obtained from a trivial  string by 
(b$)_s$. 
For $p\geq 2$, 
we have  $u(\alpha_{p,p})=0$. However 
$\alpha_{p,p}$ is homotopically non-trivial as will be shown  below. 

It follows from the definitions that $\overline {\alpha_{p,q}}=\alpha_{p,q}$ and  $ (\alpha_{p,q})^-=\alpha_{q,p}$.
Thus  the string $\alpha=\alpha_{p,q}$ with $p\neq q$ is not homotopic to  $\alpha^-, {\overline \alpha}^-$. 

Formula \ref{homr} implies that the strings  $\alpha_{p,1}$ and $\alpha_{1,p}$ with $p\geq 2$ have   minimal rank in
their homotopy classes.  We shall prove below that the same  holds for all  $\alpha_{p,q}$  except $\alpha_{1,1}$.

2. A  permutation $\sigma$ of the set $\{1,2,\ldots , m\}$ gives rise to 
a 
virtual string $\alpha_\sigma$ of rank $m$ as follows. Let $S^1=\{z\in 
\CC\,\vert \,  \vert z\vert=1\}$ be the unit circle with 
counterclockwise 
orientation. For $i=1,\ldots ,m$, let $a_i $ (resp. $b_i $) be the 
point of 
$S^1$ with real part $(i-1)/m$ and  negative (resp. positive) imaginary 
part. 
Then $\alpha_\sigma$ is formed by $S^1$ and the  $m$ arrows
$\{(a_i, b_{\sigma (i)})\}_{i=1}^m$. For the $i$-th arrow  $e_i=(a_i, 
b_{\sigma 
(i)})$,
\begin{equation}\label{pup}n(e_i)=    \sigma (i)-i.\end{equation}
This allows us to compute the polynomial $u(\alpha_\sigma)$ directly 
from 
$\sigma$.  This example generalizes the previous one since $\alpha_{p,q}=\alpha_\sigma$ for the permutation  $\sigma$
of the set $\{1,2,\ldots , p+q\}$  given by
$$\sigma(i)=\left\{\begin{array}{ll}
i+q ,~ {\rm {if}} 
\,\,\, 
1\leq i\leq p \\
\noalign{\smallskip}
i-p
,~ 
{\rm
{if}} \,\,\, p<i\leq p+q.
\end{array} \right.$$

\subsection{Properties of $u$}\label{sn:g24} We point out a few simple 
properties of the polynomial $u$. For a virtual string $\alpha$, we 
have  
$u(\alpha)=u (\overline \alpha)$. This  follows from the fact that 
if two arrows are linked positively (resp. negatively), then 
the reversed arrows are also linked positively (resp. negatively).
The transformation $\alpha\mapsto   \alpha^-$  transforms positively 
linked 
pairs of arrows into negatively linked pairs and vice versa. Therefore
  $u (  \alpha^-)= - u (\alpha)$. 
As an 
application, 
we   observe that if $u(\alpha)\neq 0$, then $\alpha$ is  not homotopic 
  to 
$\alpha^-, {\overline \alpha}^-$.

It is obvious that if a  string $\alpha$ is a  product 
of 
 strings $\alpha_1$ and $\alpha_2$, then 
  $u (\alpha)= u (\alpha_1)+ u  (\alpha_2)$.

\begin{theor}\label{th:t21}
                     An integral polynomial $u(t)$ can be realized as 
the 
$u$-polynomial   of a virtual string if and only if $u(0)=u'(1)=0$.
                     \end{theor}
                     \begin{proof}
                     We need only to prove the sufficiency of the 
condition  
$u(0)=u'(1)=0$. The proof goes by induction on the degree of $u$. If 
this degree 
  is  $\leq 1$, then $u=0$ is realized by a trivial virtual string. 
Assume that 
our claim is true for   polynomials of   degree $<m$ where $m\geq 2$. 
Let $u(t)$ be 
a polynomial of degree $m $ with highest term $a t^m$ where $a\in \ZZ$ 
and 
$a\neq 0$. Then      $v(t)=u(t)-a (t^m-mt)$ is a polynomial of degree 
$<m$ with 
$v(0)=v'(1)=0$.	     By the inductive assumption, $v(t)$ is realizable as the $u$-polynomial   of a  string. 
By 
Example \ref{sn:g23}, the polynomial $t^m-mt$ is also realizable. Taking a product  of strings we observe that  the sum of
realizable polynomials is realizable.  Hence for 
$a>0$, the polynomial $u(t)=v(t)+a (t^m-mt)  $ is realizable. If $a<0$, 
then this argument shows that $-u(t)$ is realizable by a string, $\alpha$. Then
$u(t)$ is realized by $\alpha^-$.      \end{proof}
		     
\subsection{Computation for curves}\label{sn:g25} We compute the polynomial $u$ 
 for the   string
$\alpha=\alpha_\omega $  underlying a  closed curve $\omega:S\to \Sigma$ 
on a surface $\Sigma$. The computation goes in terms of  the homological 
intersection form $B:H_1(\Sigma)\times H_1(\Sigma) \to \ZZ$ determined 
by the 
orientation of $\Sigma$. Here and below $H_1(\Sigma)=H_1(\Sigma;\ZZ)$.

 Let $e=(a,b)$ be an arrow of $\alpha$.  
Then $\omega(a)=\omega(b)$ so that    $\omega$ transforms the 
arcs $  
ab, ba\subset S$   into   loops $\omega(ab), \omega(ba)$ in $\Sigma$. Set 
$[e]=[\omega(ab)]\in 
H_1(\Sigma)$ and $[e]^*= [\omega(ba)]\in H_1(\Sigma)$ where the square 
brackets on 
the right-hand side stand for the homology class of a loop. 
We compute the intersection number $B([e], [e]^*)\in \ZZ$. The  loops 
$\omega(ab), 
\omega(ba)$      intersect transversely 
except at their common origin $\omega(a)=\omega(b)$.  
Drawing a picture of   $\omega(ab), \omega(ba)$ in a neighborhood of    
$\omega(a)=\omega(b)$, 
one  observes that   a small deformation makes these   loops   
disjoint in 
this neighborhood. 
The  transversal intersections of $\omega(ab), \omega(ba)$ bijectively correspond to the arrows of 
$\alpha$ 
linked with $e$, i.e., the arrows connecting an interior point of  $ab 
$ with  
an interior point of   $ba $. The intersection sign at such an intersection  is $+1$    if the  tail of the corresponding  arrow lies
on $ab$ and is
$-1$  otherwise. 
Adding these signs, we obtain that  $B([e], [e]^*)=n(e) $. This formula can 
be 
rewritten in a more convenient form. Set   $   s=[\omega]= [\omega(S)]\in 
H_1(\Sigma)$. Observe 
that $ s=[e]+  [e]^*$ and therefore 
$$ B([e], [e]^*)=B([e], s-[e])=B([e],s )-
B([e],  [e])=B([e],s).$$
Thus
\begin{equation}\label{pik}n(e)=B([e],s).\end{equation}
Therefore  for any $k\geq 1$, $$u_k(\alpha)=\#\{e\in \arr(\alpha)\,\vert \,   B( [e],s)=k\}
- \#\{e\in \arr(\alpha)\,\vert \,   B([e],s)=-k\}  $$
and
$$ u(\alpha)=\sum_{e\in \arr(\alpha), B( [e],s)  \neq 0} \sign ( B( 
[e],s) )\, 
t^{\vert  B( [e],s)  \vert}.$$
Using the bijective correspondence between the set $\arr(\alpha)$  and the set $\Join\! \!  (\omega)$ of   double
points of   
$\omega$, we obtain  \begin{equation}\label{orrr} u(\alpha)=\sum_{x \in \Join  (\omega),
B( [\omega_x],s)  \neq 0}
\sign ( B(   [\omega_x],s) )\, 
t^{\vert  B(   [\omega_x],s)  \vert}\end{equation}
where  for $x \in \,  \Join\! \! (\omega)$,  we let $\omega_x: [0,1]\to \Sigma$  be the loop beginning at  $x$
and following along $\omega$ until the first return to $x$
and such that the pair (a 
positive 
tangent vector of $\omega_x$ at $0$, a positive tangent vector of $\omega_x$ at 
$1$)  is a 
positive basis in   the tangent space of   $x$.

  \subsection{Coverings and higher polynomials.}\label{covvvr1} The polynomial $u$ gives rise to a   family of polynomial invariants of strings
numerated by sequences of positive integers. Their construction is based on the notion of a covering for strings. Let $\alpha$ be a string with
core circle $S$ and
$r\geq 1$ be an integer.  Let $\alpha^{(r)}$ be the string formed by $S$ and  the arrows $e\in \arr(\alpha)$ such that $n(e)\in r\ZZ$.  We call 
$\alpha^{(r)}$ the {\it $r$-th covering of $\alpha$}.  If  $\alpha$   underlies a  closed curve
$\omega:S\to
\Sigma$ on a surface $\Sigma$, then   $\alpha^{(r)}$   underlies a  lift of $\alpha$  to the $r$-fold covering of $\Sigma$
induced by the cohomology class in $H^1(\Sigma; \ZZ/r\ZZ)$ dual to $[\omega]\in H_1(\Sigma)$. Note that  $\alpha^{(1)}=\alpha$. 

	\begin{lemma}\label{covl} If   strings $\alpha$ and $\beta$ are homotopic, then $\alpha^{(r)}$ is  homotopic to  $\beta^{(r)}$ for all $
r\geq 1$.
			   \end{lemma} 
                     \begin{proof} If $\alpha$ is obtained from  $\beta$ by the  homotopy move  (a$)_s$, then the additional arrow $e$ verifies
$n(e)=0$ so that  $\alpha^{(r)}$ is  obtained from    $\beta^{(r)}$ by the move  (a$)_s$.
If $\alpha$ is obtained from  $\beta$ by the   move  (b$)_s$, then the additional arrows $e_1,e_2$ verify
$n(e_1)=-n(e_2)$. If $n(e_1)\in r\ZZ$, then   $\alpha^{(r)}$ is  obtained from    $\beta^{(r)}$ by    (b$)_s$; otherwise
$\alpha^{(r)}=\beta^{(r)}$.   Let  $\alpha$ be  obtained from  $\beta$ by the   move  (c$)_s$ replacing three arrows $e_1,e_2,e_3$ by
$e'_1,e'_2,e'_3$. It  was shown above that $n(e'_i)=n(e_i)$ for $i=1, 2, 3$.  It is easy to check that $n(e_1)+n(e_2)+n(e_3)=0$.
Three cases may occur:   the   numbers $n(e_1), n(e_2), n(e_3)$ are divisible by $r$;   one of these   numbers is divisible by
$r$ and the other two are not;  neither of these   numbers is divisible by $r$.  In the first case 
 $\alpha^{(r)}$ is  obtained from    $\beta^{(r)}$ by    (c$)_s$. In the second and third cases 
$\alpha^{(r)}=\beta^{(r)}$. This implies the claim of the lemma. \end{proof}

Iterating the coverings, we can define for a  string $\alpha$ and a finite sequence of positive integers $r_1,..., r_k $, a string
$$\alpha^{(r_1,...,r_k)}= (...( \alpha^{(r_1)})^{(r_2)}...)^{(r_k)}.$$ 
Set
$$u^{r_1,..., r_k } (\alpha)= u( \alpha^{(r_1,...,r_k)})\in \ZZ[t].$$
The results above imply that this polynomial  is a homotopy invariant of $\alpha$.

 \subsection{Exercises.}\label{sn:fgfd6}  1. Verify  that all 
virtual 
strings of rank $3$ are either homotopically trivial or  homeomorphic 
to     
$\alpha_{1,2},  \alpha_{2,1}$.

2.  For  an integer $r\geq 1$ and a virtual string $\alpha$,  define a 
virtual 
string $r\cdot \alpha$ as follows. Identifying the core circle of 
$\alpha$ with 
$S^1\subset \CC$ we can present arrows of $\alpha$ 
by  geometric 
vectors
  with endpoints on $S^1$. 
Replace each of these vectors, say $e$, by $r$   disjoint parallel 
vectors 
$e_1, \ldots , e_r$ running   closely to $e$ and having endpoints on  $S^1$. 
This 
gives  a virtual string $r\cdot \alpha$ of rank $rm$ where $m$ is the rank 
of 
$\alpha$.  Check that  $u(r\cdot \alpha) (t)=r\, u (\alpha) (t^r)$. In particular, if 
$u(\alpha)\neq 0$, then $r\cdot \alpha$ is  
homotopically non-trivial. 

3. Show that the rank  4 string   
$ \alpha_\sigma$ with 
  $\sigma=(1342)$ is homotopically trivial.    Hint:  apply to $\alpha_\sigma$ the   move ((c${)}^{+}_s)^{-1}$
as in Remark \ref{sn:g14}.3 where   
$a=a_2, a^+=a_3, b=a_4, b^+=b_4, c=b_2,c^+=b_1\in S^1$.  

4. Let $\sigma$ be a permutation of    $\{1,2,\ldots , m\}$. If $\sigma(m)=m$ and $\tau$ is the restriction of
$\sigma$ to $\{1,2,\ldots , m-1\}$, then $\alpha_\sigma$ is homotopic to $\alpha_\tau$. 
 If $\sigma(1)=m-1, \sigma(2)=1, \sigma (3)=m, \sigma (m)=2$, then $\alpha_\sigma$ is homotopic to $\alpha_\tau$
where  $\tau$ is the permutation of $\{1,2,\ldots , m-2\}$ defined by $\tau(1)=m-2$ and $\tau (i)=\sigma (i+1)-1$ for
$i>1$. 

5. Let  a string $\alpha$ be a product of the strings $\alpha_{1,3}$, $\alpha_{1,4}$, $\alpha_{2,1}$, $\alpha_{2,4}$, $\alpha_{ 3,5}$,
$\alpha_{4,3}$, $\alpha_{5,1}$, $\alpha_{5,2}$. Show that $u(\alpha)=0$ and  $\alpha^{(2)}$ is homotopic to $\alpha_{2,4}$. Thus  $u^{(2)}
(\alpha) =u(\alpha^{(2)})=2t^4-4t^2$.

6. More generally, for any $p,q,s,m\geq 1$, 
if  a string $\alpha$  is a product of the strings $\alpha_{s+m,p}$, $\alpha_{p,s}$, $\alpha_{p,m}$, $\alpha_{s+m,q}$, $\alpha_{q,s}$,
$\alpha_{q,m}$, $\alpha_{p+q, s+m}$, $\alpha_{s,p+q}$, $\alpha_{m,p+q}$, then 
$u(\alpha)=0$. If $p,s\in r\ZZ$ and $q,m$ are prime to $r$, then $\alpha^{(r)}=\alpha_{p,s}$.

 \section{Geometric realization  of virtual strings}
                    
		         \subsection{Realization   of   strings}\label{fi:g31}   
We explain here that every virtual string   admits a canonical 
realization by a closed curve on a surface and moreover describe   all  its 
realizations.   

Let  $\alpha$ be a virtual string of rank $m$ with core circle $S$.    Identifying the  tail with the  head for all arrows 
of 
$\alpha$, we transform $S$ into   a 1-dimensional CW-complex  
$\Gamma=\Gamma_\alpha$.  We   thicken $\Gamma$ to a surface $\Sigma_\alpha$ as follows.
 If $m=0$, then 
$\Gamma=S$ and we set $\Sigma_\alpha= S\times [-1,1]$. Assume 
that $m\geq 1$.  The 0-cells (vertices) of $\Gamma$ have 
valency 4 
and their number is equal to  $m$. A neighborhood of a 
vertex 
$v\in \Gamma$ embeds into a copy $D_v$ of the unit 2-disk   $ \{(x,y)\in 
\RR^2\,\vert \,  
x^2+y^2\leq 1\}$ as follows. Suppose that $v$ is obtained from an arrow 
$(a,b)$ 
where $a,b\in S$.  Note that any 
point 
$x\in S$ splits its small neighborhood in $S$ into two oriented arcs,  one 
of them being
incoming and the other one being outgoing with respect to $x$. 
A  neighborhood  of 
$v$ in 
$\Gamma$ consists of four  arcs which can be identified with small  
incoming 
and outgoing arcs of $a, b$ on $S$. We embed this neighborhood into $D_v$ so 
that $v$ 
goes to the origin and the incoming (resp. outgoing) arcs of $a,b$ 
go to 
the intervals $[-1,0]\times 0$,  $0  \times [-1,0]$ (resp. $[0,1]\times 0$, 
$0  \times 
[0,1]$), respectively. In this way the   vertices of 
$\Gamma$ can be thickened 
 to (disjoint) copies of the unit 2-disk endowed with counterclockwise 
orientation. 
Each  1-cell  of $\Gamma$ connects two (possibly coinciding) vertices 
and can be 
thickened  to a ribbon connecting the corresponding 2-disks. The 
thickening is 
uniquely determined by the condition that the orientation of these   
2-disks 
  extends to their union with the ribbon.  Thickening in this way 
all  the 
vertices and 1-cells of $\Gamma$ we embed $\Gamma$ into a 
surface  $ 
\Sigma_\alpha$.  By   
construction,  $\Sigma_\alpha$ is a   compact connected   oriented 
surface with non-void 
boundary and 
Euler characteristic   $\chi(\Sigma_\alpha)=\chi (\Gamma)=-m$. Composing the natural projection $S\to \Gamma$
with the  inclusion $\Gamma \hookrightarrow \Sigma_\alpha$ we obtain a 
closed 
curve   $\omega_\alpha:S\to \Sigma_\alpha$ realizing  $\alpha$.  
 The construction of $\Sigma_\alpha$ is well known, see 
\cite{f}, 
\cite{ca1}, \cite{cw}.

			 It is clear that for any surface $\Sigma$ and any 
(generic) closed curve $\omega:S\to \Sigma$   realizing  $\alpha$, a  
regular  
neighborhood   of $\omega(S)$
in $\Sigma$ is homeomorphic to $\Sigma_\alpha$.  Moreover, the 
homeomorphism   can  be chosen to transform $\omega$ into $\omega_\alpha$. In other 
words, $\omega$ 
can be obtained as a composition of $\omega_\alpha$ with an 
orientation-preserving 
embedding $\Sigma_\alpha\hookrightarrow \Sigma$. In particular, $\Sigma_\alpha$ is 
a   
surface of minimal genus containing a closed curve realizing  $\alpha$.  Therefore the genus   $g(\alpha)$ of $\alpha$
defined in Section 
\ref{sn:g133} is equal to the genus of $\Sigma_\alpha$. It  will be  explicitly 
computed in 
the next subsection. Note finally that a  
closed 
surface of minimal genus containing   a curve   realizing  $\alpha$  is 
obtained 
from $\Sigma_\alpha$ by gluing 2-disks to all components of $\partial 
\Sigma_\alpha$. 
			 
			  \subsection{Homological computations}\label{fi:g32} 
Consider  again a virtual string $\alpha$ of rank $m$ with core circle $S$. 
Let   
$\Gamma=\Gamma_\alpha$,  $\Sigma=\Sigma_\alpha$, and $\omega=\omega_\alpha:S\to \Sigma_\alpha$ be the
graph,  the  surface,  and the closed curve    constructed in  the 
previous subsection. 	    The 
orientation of $\Sigma$ determines a homological intersection pairing   
$B=B_\alpha: 
H_1(\Sigma)\times H_1(\Sigma) \to \ZZ$. This bilinear pairing   is  
skew-symmetric  and its rank is equal to twice the genus of 
$\Sigma$. Thus
		 $$g(\alpha) = (1/2)\, {\text{rank}}  B_\alpha.$$ In particular, 
  $\alpha$ can be realized by a closed curve on $S^2$ or $\RR^2$  if and only if $B_\alpha=0$.

Since $\Gamma$ is a deformation retract of $\Sigma 
$,  the inclusion homomorphism $H_1(\Gamma)\to 
H_1(\Sigma 
)$ is an isomorphism. Since   $\Gamma$ is  a connected graph with  $\chi 
(\Gamma)=-m$, the 
group $H_1(\Gamma)=H_1(\Sigma 
)$ is a free abelian group of rank ${m+1}$. We  
describe a  canonical 
basis in $H_1(\Sigma 
)$. 
 Set $s= [\omega ] \in H_1(\Sigma)$.  For an arrow $e=(a,b)\in 
\arr(\alpha)$, the map $\omega$ transforms the arc $  ab\subset S$, leading 
from $a$ 
to $b$ in the positive direction, into a loop $\omega(ab)$ in $\Sigma$. Set 
$[e]=[\omega(ab)]\in H_1(\Sigma)$. An easy induction on $m$ shows  that 
$s\cup 
\{[e]\}_{e\in \arr(\alpha)} $ is a basis of $H_1(\Sigma)$. 
 Our next aim is to compute    the matrix of $B$ in this 
basis.   Note 
for the 
record that $B(x,y)=-B(y,x)$ and $B(x,x)= 0$ for all elements $x,y$ of 
this 
basis.

	By Formula \ref{pik},  $B([e],s)=n(e)$ for any $e\in 
\arr(\alpha)$. To compute  the other values of $B$, 
 we need more notation. Let $a,b$ be distinct point of  
$S$.   The {\it interior} of the arc  $ab\subset S$ is the set 
$(ab)^\circ= 
ab-\{a,b\}$. For  any   arcs $ab, cd \subset S$,  we define   $ab\cdot 
cd\in 
\ZZ$ to be the number of arrows of $\alpha$ with tail in 	
$(ab)^\circ$ and head in $(cd)^\circ$	    minus the number of arrows 
of 
$\alpha$ with tail in 	$(cd)^\circ$ and head in  
$(ab)^\circ$. 
Note that the arrows with both endpoints in $(ab)^\circ\cap (cd)^\circ$ 
appear 
in this expression twice with opposite signs and therefore cancel out. 
Clearly, 
$ab\cdot cd=-cd\cdot ab$.
 In particular,  $ ab \cdot   ab =0$. If   $e=(a,b)$ is 
an arrow 
of $\alpha$, then it follows from the  definitions that $n(e)=ab\cdot  ba$.

		\begin{lemma}\label{l:t33} Let $e=(a,b)$ and $f=(c,d)$  be two 
  arrows of $\alpha$. Then $B([e],[f])=ab \cdot cd 
+\varepsilon$ where 
$\varepsilon=0$ if $e$ and $f$ are unlinked, $\varepsilon=1$ if   $f$ 
links 
$e$ positively, and  $\varepsilon=-1$ if   $f$ links $e$  negatively.
			   \end{lemma} 
                     \begin{proof}   If $e=f$, then $a=c, b=d$ and all terms of the stated equality are equal to 0.  (Note
that an arrow is unlinked  with itself.) Assume from now on that $e\neq f$ so that $a,b,c,d$ are  pairwise distinct points of
$S$.

Suppose first that
$e$ and
$f$ are unlinked. There are four cases to consider depending on whether the 
endpoints  of 
$e,f$ lie on $S$ in the cyclic order
(i)   $a,b, c,d$, or  (ii)   $a,b, d,c$, or  (iii)   $a,c,d, b$, or 
(iv)  
$a,d,c, b$.

		      In the case (i), the arcs $ ab , cd \subset S$ are 
disjoint so that   $[e],[f]\in H_1(\Sigma)$ are represented by transversal loops
$\omega(ab), 
\omega(cd)$, respectively. Then $B([e],[f])=ab \cdot cd $, cf. Section \ref{sn:g25}.
		     
		       In the case (ii), the arcs $ 
ab , dc \subset S$ are disjoint so that   $[e],[f]^{\ast}=s-[f]\in H_1(\Sigma)$ are 
represented by 
transversal loops $\omega(ab), \omega(d c)$, respectively. Hence $B([e],[f]^{\ast})= ab \cdot 
dc $  and $$B([e],[f])=B([e],s 
-[f]^{\ast})=B([e],s) - B([e],  [f]^{\ast}) =n(e) - ab \cdot 
dc= 
ab\cdot cd.$$
		     
		     In the case (iii), we have $B([e],[f])=-B([f],[e])=- c d 
\cdot ab=ab \cdot c d$ since the pair $(f,e)$ satisfies the conditions 
of (ii).
		     
		     In the case (iv),   the arcs $ 
ba , dc \subset S$ are disjoint so that   $[e]^{\ast}=s-[e], [f]^{\ast}=s-[f]\in H_1(\Sigma)$ are 
represented by 
transversal loops $\omega(ba), \omega(d c)$, respectively. 
Therefore 
$$B([e],[f])=B([e],s)+ B(s,[f])+B(s-[e],s-[f])=n(e)-n(f)+ B([e]^{\ast}, 
[f]^{\ast})
		     = n(e)- n(f) + ba \cdot dc.$$  
  It remains to observe that
		     $$ ba \cdot dc=-n(e) - ba \cdot  cd =-n(e)+cd\cdot ba= -n(e)+n(f) -cd\cdot ab=-n(e) +n(f) + ab \cdot 
cd.$$

		  Suppose that $f$ links $e$ positively. Then  their endpoints   
lie on $S$ in the cyclic order $a,c,b, d$.  The   loops 
$X=\omega(ab), Y=\omega(cd)$ representing $[e], [f]\in H_1(\Sigma)$ are not transversal since   both contain $\omega(c
b)$. Pushing 
$Y$  slightly  to its left in $\Sigma$, we obtain a loop, $Y^+$,  transversal to $X$. It is understood that the 
point 
$\omega(c)=\omega(d)\in Y$ is pushed  to a point lying between   $\omega(a c)$ and 
$\omega(d a)$	
in a small neighborhood of $\omega(c)=\omega(d)$.
Introducing coordinates $(x,y)$ in this neighborhood we can locally 
identify 
$X,Y, Y^+$ with the axis $y=0$, the union of two half-lines  $x=0, 
y\leq 0$ and 
$y=0, x\geq 0$, and the union
of two half-lines $x=-1, y\leq 1$ and $y=1, x\geq -1$, respectively. 
To compute the intersection number $B([e], [f])= X\cdot Y=X\cdot Y^+$, we  split the set $X\cap Y^+$ into 
four 
disjoint subsets. The first of them consists of a single point near 
$\omega(c)=\omega(d)$, given in the coordinates   above by $x=-1,y=0$. This point contributes $1$ to  $X\cdot Y^+$.
The second subset of $X\cap Y^+$ is   $\omega(ac)\cap Y^+$;  its points are numerated by  arrows of 
$\alpha$ with one endpoint  in the interior of $ac$ and the other 
endpoint in 
the interior of $cd$. The contribution of these crossings to 
$X\cdot Y^+$ is equal to $a c \cdot c  d$.
The third subset of $X\cap Y^+$ is numerated by   the   crossings of   
$\omega(cb)$ 
with the part of $Y^+$ obtained by pushing $\omega(bd)\subset Y$ to the 
left; they 
are numerated by arrows of $\alpha$ with one endpoint  in the interior 
of $c b$ 
and the other endpoint in the interior of $bd$. The contribution of 
these 
crossings to 
$X\cdot Y^+$ is   $c b \cdot bd$. The forth subset of $X\cap Y^+$ is 
numerated 
by   the   self-crossings of   $\omega(cb)$: each of them gives rise to two 
points of 
$X\cap Y^+$ with opposite intersection  signs. Therefore this forth subset 
contributes 0 to 
$X\cdot Y^+$.  Summing up these contributions we obtain
$$B([e],[f])=  1+  a c \cdot c  d + c b \cdot bd+0=a c \cdot c  d + 
cb\cdot cb + c 
b \cdot bd+1$$
$$=a c \cdot c  d + cb\cdot cd+1= a b\cdot cd+1.$$ 
	     
	If $f$ links  $e$ negatively, then $e$ links $f$ positively and by the 
results above,   $$B([e],[f])=-B([f],[e])=
		     -(c d \cdot a  b   +1)
		     =a  b \cdot c d -1.$$ \end{proof}			 
			    
     \subsection{Examples}\label{fi:g33}  (1) Consider the  string   
$\alpha=\alpha_{p,q}$ 
with $p,q\geq 1$ introduced in Section \ref{sn:g23}.1.  Recall the  arrows 
$e_1,\ldots, e_{p+q}$ of $\alpha$.  We  compute the  matrix of the
bilinear form $B=B_{\alpha}:H_1(\Sigma_\alpha)\times  
H_1(\Sigma_\alpha)\to \ZZ$ with respect to the basis $s\cup 
\{[e_i]\}_{i=1}^{p+q} $.  By Formula \ref{pik}, 
     $B([e_i],s) =q$ for $i=1,\ldots, p$ and $B([e_{p+j}],s)= -p$ for 
$j=1,\ldots , q$. Each pair of arrows $e_i, e_{i'}$  with $i, 
i'=1,\ldots, p$  is unlinked and by   Lemma \ref{l:t33}, $B([e_i],
[e_{i'}]) =0$. 
Similarly, each pair of arrows   $e_{p+j}, e_{p+j'}$ with   $j, 
j'=1,\ldots , q$ 
is unlinked and $B([e_{p+j}], [ e_{p+j'}])=0$. The arrow $e_{p+j}$ links 
$e_i$ 
positively and   by   Lemma \ref{l:t33},  $B([e_i],
[e_{p+j}])=(p-i)+(q-j)+1$. 
  It is easy to compute that the rank of $B$  is  equal to  2 
if 
$p=q=1$, to  6 if $\min(p,q)\geq 3$, and to $4$ in all the other cases.
 The genus $g(\alpha
) $, as we know, is half of this rank. 
 In particular, $g(\alpha_{1,1})=1$ which shows that the genus is not a 
homotopy 
invariant. 
 
 (2) Consider the  string   
$\alpha=\alpha_{\sigma}$ 
  defined  in Section \ref{sn:g23}.2  for a permutation $\sigma$ of the set $\{1,2,\ldots , m\}$. 
Recall the  arrows 
$e_1,\ldots, e_{m}$ of $\alpha$. 
We  compute the  
matrix of the bilinear form $B=B_{\alpha}:H_1(\Sigma_\alpha)\times  
H_1(\Sigma_\alpha)\to \ZZ$ with respect to the basis $s\cup 
\{[e_i]\}_{i=1}^m $. The number  $B([e_i],s) = n(e_i)$ is computed by Formula     \ref{pup}.
Pick two indices $i,j$ with  $1\leq i < j\leq m$.  Lemma \ref{l:t33} implies that if $\sigma(i)<\sigma(j)$, then
$$B([e_i], [e_j])=\#\{k\,\vert\, i<k<j, \sigma(j)<\sigma(k)\}\,-\, \#\{k\,\vert\, j<k\leq m,
\sigma(i)<\sigma(k)<\sigma(j)\}.$$
If $\sigma(j)<\sigma(i)$, then
$$B([e_i], [e_j])=\#\{k\,\vert\, i<k<j, \sigma(j)<\sigma(k)\}\,+\, \#\{k\,\vert\, j<k\leq m,
\sigma(j)<\sigma(k)<\sigma(i)\}+1.$$

			\section{Cobordism  of strings and the slice genus}\label{fi:g40}

\subsection{Cobordism of strings}\label{sn:nbee}  Two strings  $\alpha$ and
$\beta$ are {\it cobordant} if there are an oriented   3-manifold $M$  and two disjoint closed curves   on $\partial M$
realizing  
$\alpha$ and
$\overline \beta$, respectively,  and  homotopic to each other in $M$. 
Here two  curves  on a surface are {\it disjoint} if their images are disjoint subsets of the surface. The   involution
$\beta\mapsto  \overline \beta$  is needed in the definition of cobordism to ensure the reflexivity, cf. the proof of the 
next  lemma. 

In the definition of cobordism, we do not require   $M$ to be compact or connected. However, replacing 
$M$ by a regular neighborhood of a homotopy   relating the two curves  in $M$ we can always assume that $M$ is compact
and connected.

	\begin{lemma}\label{ldd:nbl} Cobordism is an equivalence relation on the  set of (homeomorphism classes of) strings.  If
strings   $\alpha$ and
$\beta$ are   cobordant, then $\overline \alpha$ is cobordant to $\overline \beta$ and $\alpha^-$ 
is cobordant to $ \beta^-$.
			   \end{lemma} 
                     \begin{proof} To see that a string $\alpha$ is cobordant to itself, we  realize  $\alpha$  by a
closed curve
$\omega$ on a  closed surface $\Sigma$ and   take $M=\Sigma \times [0,1]$ with   the orientation obtained as the
product of the orientation in
$\Sigma$ and the right-handed orientation in $[0,1]$. Then
$\omega\times 1$ and
$\omega\times  0$ are disjoint closed curves 
  on  
$\partial M=(\Sigma \times 1) \cup (-\Sigma \times 0)$ realizing  respectively $\alpha$ and
$\overline \alpha$. These curves are homotopic   in $M$. 

The relation of cobordism is symmetric:  if  two    disjoint closed curves   on   $\partial M$ realize strings
$\alpha, \overline \beta$, then    the same curves on  $-\partial M=\partial (-M)$
realize  $ \overline \alpha,\beta$. 

Suppose that a string $\alpha_1$ is cobordant to $\alpha_2$
and $\alpha_2$ is cobordant to $\alpha_3$. Let $M,  M'$ be  disjoint oriented 3-manifolds such that $\alpha_1, \overline
\alpha_2$ are realized by  disjoint closed curves 
$\omega_1,
\omega_2$  on  
$\partial M$  homotopic  in $M$ and  $\alpha_2, \overline \alpha_3$ are realized by  disjoint closed curves  $\omega'_2,
\omega_3$  on   $\partial M'$
 homotopic   in $M'$.  Observe that a regular neighborhood, $U$,  of $\omega_2$ in $\partial M$ is
homeomorphic to a regular neighborhood, $U'$,  of $\omega'_2$ in $\partial M'$ via an orientation reversing
homeomorphism transforming $\omega_2$ into $\omega'_2$.  Gluing $M$ to $M'$ along   $U\approx  U'$ 
we obtain an oriented 3-manifold
$N$. The curves $\omega_1, \omega_3$ lie  on $\partial N$ and realize  
 $\alpha_1 , \overline \alpha_3$, respectively. These curves are disjoint  and   homotopic  in $N$. Hence
$\alpha_1 $ is cobordant to $ \alpha_3$. 

The last claim of the lemma is obtained by inverting orientation on the ambient 3-manifold (resp.  on the circle). \end{proof}

 \begin{theor}\label{th:103654} Homotopic strings are cobordant. \end{theor}
  \begin{proof}    For strings $\alpha,
\beta$, we  write $\alpha\sim \beta$ if  these   strings can be realized by homotopic  (possibly intersecting) closed curves
on the same surface. The relation $\sim$ is reflexive and symmetric  but not transitive.  The    following lemma    shows that
the relation of homotopy is  precisely the equivalence relation generated by $\sim$. 

	\begin{lemma}\label{l:nbl} Two strings $\alpha,
\beta$ are homotopic if and only if there is a sequence of strings $\alpha_1=\alpha, \alpha_2,..., \alpha_n=\beta$ such that
$\alpha_i\sim \alpha_{i+1}$ for $i=1,..., n-1$.
			   \end{lemma} 
                     \begin{proof} As we know, the underlying virtual strings of homotopic closed curves on 
a surface are themselves homotopic.  Therefore if there is a sequence of strings $\alpha_1=\alpha, \alpha_2,...,
\alpha_n=\beta$ such that
$\alpha_i\sim \alpha_{i+1}$ for $i=1,..., n-1$, then $\alpha$ is homotopic to $\beta$.  To prove the converse, it suffices
to show that if $\beta$ is obtained from $\alpha$ by  a homotopy move  (a$)_s$, 
(b$)_s$,  or
(c$)_s$, then  $\alpha\sim \beta$.  

Let $S$ be the core circle  of $\alpha$ and $\omega:S\to 
\Sigma $ be 
a   curve    realizing  $\alpha$ on a surface $\Sigma$.  Pick   distinct points $a,b 
\in S$  
such that the arc   $ab\subset S$ does not contain  endpoints of 
$\alpha$.  
Let 
   $\beta$  be   obtained from $\alpha$ by the move (a$)_s$ adding    the arrow  $ (a,b)$.  
   Attaching to   $\omega$ a small curl 
on the 
right  of   $\omega(ab)$, we obtain a closed curve 
$\omega':S\to \Sigma $ realizing $\beta$.   Clearly, $\omega'$ is homotopic to $\omega$ in $\Sigma$.  Hence
 $\alpha\sim \beta$.

Pick two    arcs  $x,y$ on $S$ disjoint from each other and  
from the   
endpoints of $\alpha$. Let $a, a'$ be the endpoints of $x$ (in an 
arbitrary 
order) and $b,b'$ be the endpoints of $y$. Let $\beta$ be obtained from $\alpha$ by the move  (b$)_s$ adding 
to $\alpha$ the  arrows $ (a,b)$ and $ (b',a')$.    Let $D_x, D_y\subset
\Sigma - \omega(S)$ be two small closed disks  lying near 
the arcs $\omega(x), \omega(y)$, respectively. Removing the interiors of 
these disks 
from $\Sigma$ and gluing the circles $\partial D_x,\partial D_y$ along an 
orientation reversing homeomorphism, we obtain a new (oriented) surface, 
$\Sigma'$, 
containing $\omega(S)$. In $\Sigma'$ the arcs $\omega(x)$ and $\omega(y)$ are 
adjacent 
to the component of 
 $\Sigma'-\omega(S)$ containing $\partial D_x=\partial D_y$.  We can push 
$\omega(x)$ 
across this component towards $\omega(y)$  and eventually across $\omega(y)$. This gives a  curve $\omega':S\to
\Sigma'$   realizing $\beta$ and homotopic to $\omega:S\to \Sigma'$. 
Hence
 $\alpha\sim \beta$. Note that the four 
possible 
forms of the   move   (b$)_s$ (depending on whether $x$ leads from $a$ to $a'$ or  from $a'$ to $a$   and similarly  for
$y$) are realized by choosing  
$D_x, D_y$  on the 
left or on the right of $\omega(x), \omega(y)$.

Suppose that $\alpha$ has three arrows 
$(a^+,b), (b^+,c), (c^+,a)$ where  $a, a^+, b,b^+$, $c, c^+\in S$ such that
the  (positively oriented) arcs   $aa^+$,   $bb^+$,    $cc^+$ are disjoint from each 
other and 
from the other   endpoints of $\alpha$. Let $\beta$ be obtained from $\alpha$ by the move  (c$)_s$ replacing  the 
arrows 
$(a^+,b), (b^+,c), (c^+,a)$ with the arrows  $ (a,b^+), (b,c^+), (c,a^+)$.
  Consider the canonical realization $ \omega_\alpha:S
\to
\Sigma_\alpha$ of $\alpha$.  Observe that  the    arcs $ \omega_\alpha  (aa^+),  \omega_\alpha (bb^+), \omega_\alpha  (cc^+)$     form  a simple
closed curve in
$\Sigma_\alpha$  isotopic to a  component of   $\partial \Sigma_\alpha$. Gluing  a 2-disk  $D$
  to this  component of   $\partial \Sigma_\alpha$ we embed  $\Sigma_\alpha$ into a bigger 
surface, $\Sigma$. Pushing     $ \omega_\alpha(aa^+)$ 
across $D\subset \Sigma$ 
and then across the double point $\omega(b^+)=\omega (c)$, we obtain a  curve $\omega':S\to
\Sigma$   realizing $\beta$ and homotopic to $ \omega_\alpha:S\to \Sigma$. 
Hence
 $\alpha\sim \beta$. 
\end{proof} 

We   now accomplish the  proof of Theorem \ref{th:103654}.   The obvious cylinder construction shows   that 
the   underlying virtual strings of   homotopic   curves on 
a  surface are  cobordant.  Thus if   we have  strings $\alpha, \beta$ with  $\alpha\sim \beta$, then $\alpha$ is
cobordant to $\beta$. Lemmas  \ref{ldd:nbl}   and 
\ref{l:nbl} now imply the claim of the theorem. 
\end{proof}

 \begin{theor}\label{th:0388} The polynomial $u$ is a cobordism invariant of strings.   \end{theor}

\begin{proof}  We begin with a lemma.

	\begin{lemma}\label{l:nblddddd}  Let $ F $ be   a  compact (oriented) surface     whose boundary is a union
of $r\geq 1$ circles $S_1,...,S_r$  with  the induced orientation.    Let $M$ be an oriented  3-manifold and $\omega:F\to
M$ be a map
  such that  		$\omega(\partial F)\subset \partial M$,   $\omega(S_i)\cap
\omega(S_j)=\emptyset$ for all $i\neq j$,  and  the restriction of $\omega$ to each circle $S_i$ is a generic (closed) curve in 
$\partial M$ for $i=1,...,r$. 	  Let $\alpha_i$ be the virtual string underlying the latter curve. If the genus of $F$ is 0, then
\begin{equation}\label{osn} \sum_{i=1}^r u(\alpha_i)=0. \end{equation}
\end{lemma} 
                     \begin{proof}   We need a few general facts about maps       $F\to M$.   
A point  $ a\in  \Int F$  is a {\it simple branch point} of  a map $\omega:F\to M$ if   there is a closed 3-ball $D^3\subset  \Int M$
such that  $\omega(F)\cap D^3$ is the cone over a figure eight curve in $S^2=\partial D^3$ with cone point $\omega(a) \in \Int D^3$.
Here by a figure eight curve in $S^2$ we mean a closed curve with one transversal self-intersection. 
The set  of simple branch points of $\omega$  is  denoted $Br(\omega)$.

 Set $  \RR_+=\{r\in \RR, r\geq 0\}$.  We say that  a map   $\omega:F\to M$    is {\it
generic} if   $\omega^{-1} (\partial M)=\partial F$ and 

 (i)  the restriction of $\omega$ to   $\partial F$ is an immersion into  $\partial M$; any point of $ \omega(\partial F )$ has  a
neighborhood
$V\subset M$ such that the pair $(V, V\cap \omega(F))$ is homeomorphic to either  $(\RR^2, \RR  \times 0)  \times \RR_+
$  or
$(\RR^2,  \RR  \times 0  
  \cup 0\times
\RR )\times \RR_+$;

(ii) the  restriction of $\omega$ to  
$\Int F-Br(\omega) $ is   an immersion into   $\Int M- \omega(Br(\omega))$; 
    any point of $ \omega(\Int F-Br(\omega) )$ has  a neighborhood $V\subset \Int M$ such that the pair $(V, V\cap \omega(F))$ is
homeomorphic to either  $(\RR^3, \RR^2 \times 0  )$,  or $(\RR^3, \RR^2 \times 0 \cup 0\times
\RR^2)$ or  $(\RR^3, \RR^2 \times 0 \cup 0\times
\RR^2 \cup  \RRÊ\times 0\times \RR)$.

 Fix  a generic map $\omega:F\to M$. Set   $T=Br(\omega)\cup  \{a\in F\,\vert\,   \# \omega^{-1}
(\omega(a))
\geq 2\}\subset F$.   It is clear from (i), (ii) that  
$T$ consists of a finite number of   immersed  circles  in $\Int F$ and  immersed  proper intervals in $F$.   These circles
and intervals have only double transversal crossings.   The set of these crossings 
$\Join\! \! (T)$ is   finite and  consists precisely of the preimages of the triple points of
$\omega$.  The set $\partial T=  T \cap \partial F$ is   finite and  consists precisely of the preimages of the double points of
$\omega\vert_{\partial F}:\partial F \to  \partial M$.   The set $Br(\omega)\subset T$ is also finite and disjoint from 
$  \Join\! \! (T)$ and
$\partial T$.

Let  $\tilde T$ be  an abstract  1-dimensional manifold parametrizing $ T$. The   projection $p: \tilde T
\to T$ is 2-to-1 over  $\Join\! \! (T)$ and   1-to-1 over $T-  \Join\! \! (T)$.   
We shall identify the points of $T - \Join\! \! (T)$ with their preimages under $p$. 

For any  point $a\in T -(\Join\! \! (T)\cup Br(\omega))$
there is exactly one   other point
$b \in T-(\Join\! \! (T)\cup Br(\omega))$ such that $\omega(a)=\omega(b)$.  The correspondence $a    \leftrightarrow b$  extends by continuity
to an involution 
$\tau$  on
$\tilde T$. The set   of   fixed
points  of $\tau$  is    $Br(\omega)$.  It is clear that $\omega\, p\,\tau=\omega\, p: \tilde T\to M$.  For   $a\in \partial
T=\partial 
\tilde T$, the  point
$ \tau(a)  $ is the unique point   $b\in \partial T-\{a\}$ such that $\omega(a)=\omega(b)$.  

We   define an  involution  $\mu:\partial T\to   \partial T$.   For   $a\in \partial T=\partial 
\tilde T$,  let $I_a \subset  \tilde T$ be the  
  interval    adjacent to $a$ and  let  $\mu(a)\in
\partial T$     be its   endpoint   distinct from $a$.  Clearly, $\mu^2=\id$.   We claim that $\mu$ commutes with $\tau\vert_{\partial T}$.
Indeed, if  $
\tau(I_a)=I_a$, then   $\tau$ exchanges the endpoints of
$I_a$ so that   $ \tau=\mu $ on $\partial I_a$.  (In this case $\tau$ must have a unique fixed point on $I_a$,  so that    
$I_a$ contains   a unique branch point of $\omega$.)  If $\tau(I_a)\neq I_a$,  then  $\tau(I_a)$ has the endpoints $\tau(a), 
\tau(\mu(a) )$ so that $\mu(\tau(a))=\tau (\mu(a))$.    Since $\tau \vert_{\partial T}$ and $\mu$ commute,   $\mu$
induces an involution on $\partial T/\tau=\Join
\! \!(\omega\vert_{\partial F})$. The latter involution is denoted $\nu$.

Assume from now on that $\omega$ maps the components $S_1,...,S_r$ of $\partial F$ to disjoint subsets of $\partial M$. 
Each crossing point $x\in \,\,
 \Join
\! (\omega\vert_{\partial F})$  of $\omega\vert_{\partial F}$  is then a self-crossing of  $\omega(S)$ for a certain component $S=S_i$ of
$\partial F$.  Consider the  loop $\omega_x$   in
$\partial M$  beginning at 
$x$ and following along $\omega(S)$ until the first return to $x$
and such that the pair (a 
positive 
tangent vector of $\omega_x$ at $0$, a positive tangent vector of $\omega_x$ at 
$1$)  is a 
positive basis in   the tangent space of   $x$ in $\partial M$.  Let 
 $[\omega_x]\in H_1(\partial M)$ be the homology class of $\omega_x$  and $\inc:H_1(\partial M)
\to H_1(M)$ be the inclusion homomorphism.   We have either  $ \nu(x)=x$ or $ \nu
(x)\neq x$. We   prove    that in the first case   $\inc ([\omega_x]
)\in \omega_*(H_1(F))\subset H_1(M)$ and in the second case 
$$ \inc ([\omega_x]  + [\omega_{\nu (x)}])\in  \omega_*(H_1(F) ) \subset H_1(M) .
$$

Let $\omega^{-1}(x)=\{a,b\} $ where  $a,b\in S\cap \partial T $  and the corresponding  arrow
  is directed from $a$ to $b$ as in  Section \ref{sn:g12}.  Denote the positive arc
$ab\subset S$ by
$\gamma_x$ and observe that $\omega_x=\omega(\gamma_x)$.
Suppose   that  $ \nu(x)=x$.  Then
$\mu(a)\in
\{a,b\}$. Inspecting the   orientations of the  sheets of $\omega(F)$ meeting along $\omega(T)$, we observe
that
$\mu:\partial T\to   \partial T$ transforms arrowtails into arrowheads and vice versa (this was  first pointed out in \cite{ca2}).  Therefore
$\mu(a)=b$.   By the definition of $\tau$, we have $\tau(a)=b$ and $
\tau(b)=a$.  Since       $\tau$ preserves the set
$\partial I_a=\{a,  b\}$, we    have  $\tau(I_a)=I_a$. Observe that
the product of the path  $\gamma_x=ab\subset S$   with the immersed  interval 
$p(I_a)\subset F$ oriented  from
$ b$ to $a$  is a loop in $F$, say  $\rho$. The loop  $\omega(\rho) $ in $M$  is a product of  
$\omega(\gamma_x)=\omega_x$  with the loop $\omega p\vert_{I_a}$. The latter  loop    has the form $ \delta 
\delta^{-1}$ where
$\delta $ is the path in $  M$ obtained by restricting
$ \omega\,p$ to the arc in $I_a$ leading from $b$ to   the unique   branch point of $\omega$  on $ I_a$.  This loop   is
contractible in $M$.  Hence
$
\inc ([\omega_x])=[\omega(\rho)]\in \omega_*(H_1(F))$. 

Suppose    that  $ \nu(x)\neq x$.  Note that  the  path $\gamma_{\nu(x)} $ begins  at    $\mu(b)$ and  
terminates  at $ \mu(a)$.
Consider   the  loop $\rho= 
\gamma_x \,p (I_{b})\,
\gamma_{\nu(x)} \, (p(I_{ a}))^{-1}
$ in $F$ beginning and ending at $a$. Here the intervals   $ I_{b}
,   I_{a}$ are   oriented  from
$b$ to $\mu(b)$ and from $a$ to $\mu(a)$, respectively.   Then   $\omega(\rho)$ is the product
of the loop $\omega_x$ beginning and ending in $x$, the path $\omega p ( I_{b})$ beginning in $x$ and ending
in $ \nu(x)$, the   loop $\omega_{\nu(x)}$ beginning and ending in $ \nu(x)$, and  the path $(\omega p( I_{ a}))^{-1}$
beginning in
$ \nu(x)$ and ending in $x$.   The  
paths 
$\omega p( I_{b})$,  $(\omega p( I_{ a}))^{-1}$ are mutually inverse  since  $  I_{b} = \tau( I_{a}) $ and $\omega p \tau
=\omega p$. Hence 
$$ \inc ([\omega_x]  + [\omega_{\nu (x)}]) =[\omega(\rho)]\in
\omega_*(H_1(F)).$$

  Denote by
$B$ the intersection form $H_1(\partial M)\times H_1(\partial M) \to \ZZ$. Set $s_i=[\omega(S_i)]\in H_1(\partial M)$ and
$s=s_1+s_2+\cdots +s_r\in H_1(\partial M)$. By Section
\ref{sn:g25},  
\begin{equation}\label{sedrrttt}
\sum_{i=1}^r u(\alpha_i)= \sum_{i=1}^r   \,\,  \sum_{x \in
\Join  (\omega(S_i)), B( [\omega_x], s_i)  \neq 0}
\sign ( B(   [\omega_x],s_i) )\, 
t^{\vert  B(   [\omega_x],s_i)  \vert}.\end{equation}
   For  $x \in \,  \Join \!\!  (\omega(S_i))$,  the loop $\omega_x$ lies on
$\omega(S_i)$ and is      disjoint from $\cup_{j\neq i} \, \omega(S_j)$. Hence   
$B( [\omega_x], s)=B( [\omega_x], s_i)$. Note also that $\Join
\! (\omega\vert_{\partial F})=\cup_{i=1}^r  \Join
\! (\omega(S_i)) $. Therefore Formula  \ref{sedrrttt} can be rewritten as
\begin{equation}\label{se968t}
\sum_{i=1}^r u(\alpha_i)= \sum_{x \in
\Join  (\omega\vert_{\partial F}), B( [\omega_x], s)  \neq 0}
\sign ( B(   [\omega_x],s) )\, 
t^{\vert  B(   [\omega_x],s)  \vert}.\end{equation}

Assume   that the genus of $F$ is zero.  We  shall show  that each orbit of the involution $\nu$ on $\Join
\! \! (\omega\vert_{\partial F}) $ contributes $0$ to the right hand side of Formula \ref{se968t}.  This will imply the claim of
the lemma. It suffices to prove that   for any  $x\in\,  \Join
\! \! (\omega\vert_{\partial F})  $ we have $B( [\omega_x], s)=0$ in the case $\nu (x)=x$ and $B( [\omega_x],
s)=-B( [\omega_{\nu(x)}], s)$ in the case $\nu (x)\neq x$.  Set $h_x=[\omega_x]\in H_1(\partial M)$  if  $\nu (x)=x$ and $h_x=  [\omega_x]+
[\omega_{\nu(x)}]\in H_1(\partial M)$ if  $\nu (x)\neq x$. As we know, 
$\inc (h_x)
\in \omega_*(H_1(F))$.  Since the genus of $F$ is 0, the group $H_1(F)$ is generated by the homology classes of the boundary
components. Hence there is  an integral  linear combination
$h=h(x)$ of 
$s_1,...,s_r
\in H_1(\partial M)$ such that  $\inc (h_x)=\inc(h)$. Then 
$h_x- h \in K$ where $K$ is  the kernel of  $\inc: H_1(\partial M) \to H_1(M)$.  
The sum  $s=s_1+...+s_r$ being represented by the 1-cycle $\partial \,\omega(F)$   also lies in $K$.  It is well known that  
$B(K\times K)=0$.  Hence 
$B(h_x- h, s)=0$. Since the curves $\omega(S_1),...,\omega(S_r)$ are pairwise disjoint and $B$ is skew-symmetric,   $B (h,
s)=0$. Therefore
$B(h_x , s)=0$.    This implies  our claim.
 \end{proof}

We can now finish the proof of the theorem.  Let $\alpha,
\beta$  be   cobordant
strings.      By assumption, there is an oriented 3-manifold
$M$ and a   homotopy $\{\omega_t:S^1\to M\}_{t\in [0,1]}$   such that $\omega_0,\omega_1$  are disjoint (generic) 
closed curves on
$\partial M$ realizing $\alpha$ and $\overline \beta$, respectively.   The homotopy $\{\omega_t \}_{t}$ defines
a map   $\omega:S^1\times [0,1]\to M$.  We provide $S^1\times [0,1] $ with the orientation obtained as the product of the
counterclockwise orientation in
$S^1$ and the right-handed orientation in $[0,1]$.  Applying Lemma \ref{l:nblddddd} to $\omega$ we obtain that
$u(\alpha_0)+u(\alpha_1)=0$ where $\alpha_i$ is the string underlying the restriction of $\omega$ to $S^1 \times i$ where the
orientation of $S^1 \times i$ is induced by the one in $S^1\times [0,1] $. This  is the  counterclockwise  orientation on 
$S^1 \times 1$ and the opposite one on  $S^1 \times 0$.  Therefore
$\alpha_1=\overline
\beta$ and
$\alpha_0=\alpha^-$. Hence   
$u(\alpha)=-u(\alpha^-)=-u(\alpha_0)=u(\alpha_1)=u(\overline \beta)=u(\beta)$.  
\end{proof} 

\begin{corol}\label{aadddal:gg1} The strings $\alpha_{p,q}, \alpha_{p',q'}$ with $p\neq q, p'\neq q'$ are
cobordant if and only if $p=p'$ and $q=q'$.
		   \end{corol}

This follows   from the previous theorem and the formula $u(\alpha_{p,q})=p t^{q}- 
q t^p$. The strings   $\alpha_{p,q}$ with $p=q$ are all cobordant to each other as will be shown in the next subsection.

\begin{corol}\label{aadddal:gg222} For any integers $r_1,...,r_k\geq 1$, the polynomial
  $u^{(r_1,...,r_k)}$ is a cobordism invariant of strings.
		   \end{corol}
\begin{proof} It suffices to prove that if a string $\alpha$ is cobordant to a string $\beta$, then $\alpha^{(r)}$ is cobordant to
$\beta^{(r)}$ for  $r\geq 1$.  Let $M $ be  a compact    oriented 3-manifolds such that $\alpha , \overline
\beta$ are realized by  disjoint closed curves 
$\omega,
\omega'$  on  
$\partial M$  homotopic  in $M$. Let $f: S^1\times [0,1] \to M$ be a homotopy between $\omega=f \vert_{S^1\times 0}$ and 
$\omega'=f \vert_{S^1\times 1}$.  By the Poincar\'e duality,  there is a unique $y\in H^1(M;\ZZ)$ such that $y
\cap [M]=[f]\in H_2(M,\partial M; \ZZ)$.  Let $\tilde M\to M$ be the $r$-fold covering determined by $y(\modu r)\in H^1(M; \ZZ/r\ZZ)$. 
The mapping $f$ lifts to $\tilde M$ and yields a homotopy between $\alpha^{(r)}$ and $(\overline \beta)^{(r)}= \overline
{\beta^{(r)}}$.
\end{proof}

  \subsection{Slice strings}\label{fbm}  A virtual string cobordant to a trivial string is   {\it slice}. 
Clearly, a trivial string is
slice.  Lemma \ref{ldd:nbl} implies that   strings  cobordant  to a slice string are 
slice.  By Theorem \ref{th:103654}, a string homotopic   to a slice string is  
slice.  By the proof of Corollary \ref{aadddal:gg222}, all coverings of a slice string are slice. 

    Theorem \ref{th:0388} gives   obstructions to the
sliceness:  the  polynomial 
$u$ and the higher polynomials $u^{(r_1,...,r_k)}$ of a slice   string are equal to $0$. For example, the strings $\alpha_{p,q}$ with $p\neq q$ are
not slice.  

It is easy to see that   a string is slice  if  and only if it can be realized on a closed surface $\Sigma$ by a  closed  curve  
contractible in an orientable 3-manifold bounded by $\Sigma$.  Using the gluing of 3-manifolds along 2-disks in the 
boundary, we 
 obtain that a string that is a product of  slice strings is itself slice.  Similarly, using  the gluing of 3-manifolds along
(subsurfaces of) their boundary we obtain the following cancellation:  if a product of a string $\alpha$ with a slice
string is slice  then $\alpha$ is slice.

We outline a construction  of slice strings which mimics the well known fact    that a
sum of a knot with its  mirror image is slice.  Namely, for any  virtual string $\alpha$, its appropriate product with
$\overline
\alpha^-$ is slice.  Indeed, let  
$S$ be the core circle of $\alpha$  and let $ab\subset S$  be an arc containing all the   
endpoints of 
$\alpha$. Let $(\alpha', S', a'b'\subset S')$ be a disjoint copy of the triple $(\alpha, S, ab)$.  Consider the circle  
$S''= (ab\cup a'b')/a=a', b=b'$  and provide  it with the orientation extending the one on $ab$. The arrows
of
$\alpha$ and $\overline {\alpha'}$ are attached to $ab\cup a'b'$ and form in this way a virtual string, $\alpha''$, with core
circle 
 	$S''$.  It is clear that $\alpha''$ is a product of $\alpha$ with $\overline \alpha^-$. We claim that $\alpha''$ is slice.  To
see this, represent $\alpha$ by a closed curve $\omega:S\to \Sigma$ on a  surface $\Sigma$.  The map $\omega$
transforms   $S- ab$ onto an embedded arc in $\Sigma$ disjoint from the rest of the curve. Let $D\subset
\Sigma$ be a 2-disk such that  $D\cap \omega(S)=\omega(S-ab)$ and   $\partial D\cap
\omega(S)=\{\omega( a), \omega(b)\}$.  Consider the 3-manifold $M=(\Sigma-\Int D) \times [0,1]$. The
four paths
$\omega(ab)\times 0, \omega(ab)\times 1, \omega(a) \times [0,1], \omega(b) \times [0,1]$ form  a closed curve on
$\partial M$   realizing $\alpha''$ and contractible in $M$.

The analogy with knot theory  suggests the following definition.  A string $\alpha$   is  
 {\it ribbon} if  its core circle   has
 an orientation reversing involution  $j$   such that for any
arrow $(a,b)$ of $\alpha$ the pair $(j(b),j(a))$ is also an arrow of $\alpha$.  Note that such an involution  $j$   is
topologically equivalent to the complex conjugation on $S^1$ and, in particular,   has two fixed points.  The assumptions
on $\alpha$ imply that these two points are not endpoints of $\alpha$.    For example, it is obvious that the string
$\alpha''$ constructed in the previous paragraph is ribbon.  Another   example:  the string $\alpha_{p,p}$  with  $p\geq 1$ is
ribbon. The ribbonness of a string is not a homotopy   property: a string obtained from a ribbon string
by homotopy moves may be non-ribbon.

The next lemma shows that all ribbon strings are  slice.

	\begin{lemma}\label{ribb} Ribbon strings    
are  slice.
			   \end{lemma} 
                     \begin{proof} Let $S$ be the core circle of  a ribbon string $\alpha$ and  let  $j:S\to S$ be an orientation
reversing involution transforming     arrows of $\alpha$ into   arrows of $\alpha$  with opposite
orientation. Recall    the canonical realization 
$\omega_\alpha: S\to
\Sigma_\alpha$ of  $\alpha$.  The involution  $j$  induces an involution  $j'$   on  the graph 
$\Gamma_\alpha=\omega_\alpha(\Sigma_\alpha)$.  We  extend  $j'$ to the disks  $\{D_v\}_v$ used to
construct 
$\Sigma_\alpha$ by $D_v\to D_{j'(v)},\, (x,y)\mapsto (-y,-x)$  where $v$ runs over the  vertices  of  $ 
\Gamma_\alpha$ and $x,y$ are the canonical coordinates in these disks,
cf.  Section
\ref{fi:g31}.  The resulting involution  extends  to the ribbons in the obvious way and  yields  an orientation reversing
involution
$ J: \Sigma_\alpha\to \Sigma_\alpha$ such that $J \omega_\alpha=\omega_\alpha j$.     The set $\Fix (J) $ of fixed
points of 
$J$ consists of two disjoint embedded intervals  in $\Sigma_\alpha$  with endpoints on $\partial \Sigma_\alpha$.
Consider the cylinder $ \Sigma_\alpha\times [0,1]$ and identify $a\times 0= J(a) \times 1$ for all $a\in
\Sigma_\alpha-\Fix (J)$. For each  $a\in \Fix (J)$, contract   $a\times [0,1]\subset \Sigma_\alpha\times [0,1]$ into a
point. This transforms  $
\Sigma_\alpha\times [0,1]$ into an oriented  3-manifold 
$M$  such that $\partial M \supset \Sigma_\alpha$ and $\omega_\alpha$ is contractible in $M$.   \end{proof}
 
\begin{remar}\label{rered}  Not all slice strings are ribbon.  To give an example, consider the ribbon    string
$\alpha_{1,1}$.  Since   $\alpha_{1,1}$  is slice, any  string obtained  as a product of $\geq 2$  copies of 
$\alpha_{1,1}$  is slice. Some of such  products are not ribbon.   For example, consider   the permutation
$\sigma=(12)(34)$ on the set
$\{1,2,3,4\}$ permuting 1 with 2 and 3 with 4  and consider the rank  4  string
$ \alpha_\sigma$  defined in Section \ref{sn:g23}.2.  Drawing a picture,
one  observes that  $\alpha_\sigma$  is a product of two copies of 
$\alpha_{1,1}$
and is ribbon. Inverting  orientation of any arrow of  $ \alpha_\sigma$ we obtain a string which is also 
a product of two copies of  $\alpha_{1,1}$ but which is not ribbon by obvious geometric reasons.
\end{remar}

 \subsection{Slice genus}\label{sgenu}   The {\it slice  genus}
$sg(\alpha)$ of a string
$\alpha$ is the minimal  integer
$k\geq 0$ satisfying the following condition:   there are an oriented 3-manifold $M$, a compact (oriented) surface  $F$ of
genus
$  k$   bounded by a circle,  and a proper map
$\omega:F \to M$ such that
 $\omega\vert_{\partial F }:\partial F \to
\partial M$ is a (generic) closed curve on $\partial M$ realizing $\alpha$.  (The word \lq\lq proper" means that 
 $\omega (\partial F)\subset \partial M$.) Such 
$k$ exists because   any loop on a closed surface is homologically trivial in a certain  handlebody bounded
by this surface. It is clear that $sg(\alpha)$ is a cobordism invariant of  
$\alpha$. A string  $\alpha$ is slice if and only if $sg(\alpha)=0$.

We   similarly define a	   slice genus for tuples of strings.    The {\it   slice  genus}
$sg(\alpha_1,..., \alpha_r)$ of   $r\geq 1$   strings
$\alpha_1,..., \alpha_r$ is the minimal  integer
$k\geq 0$ satisfying the following condition:   there are an oriented 3-manifold $M$, a compact
(oriented) surface $F$   of genus
$ k$   bounded by $r$ circles $S_1,...,S_r$,   and a proper map
$\omega:F\to M$ such that
  the  maps $\omega\vert_{S_i}:S_i \to
\partial M $ with $i=1,...,r$ are disjoint (generic) closed curves on $\partial M$ realizing $\alpha_1,..., \alpha_r$,
respectively.  The existence of such 
$k$ can be obtained by realizing $\alpha_1,..., \alpha_r$ by curves on disjoint surfaces, taking the  connected sum of
these surfaces  and presenting the result as a boundary of   an appropriate   handlebody.  We do not require
$M$ or $F$
 to be connected although it is always possible to achieve their connectedness by taking connected sum. Note that  the
genus of a disconnected surface is by definition the sum of the genera of its components. 

Clearly $sg(\alpha_1,..., \alpha_r)\geq 0$.   If  $sg(\alpha_1,..., \alpha_r)=0$ then we call the sequence 
$\alpha_1,..., \alpha_r$  {\it slice}.   The same argument as in the proof of Corollary \ref{aadddal:gg222} shows that 
 if $\alpha_1,..., \alpha_r$ is slice, then for any integer $m\geq 1$, the sequence of the  $m$-th coverings 
$  \alpha_1^{(m)},..., \alpha_r^{(m)}$ is slice.
By Lemma \ref{l:nblddddd}, if $\alpha_1,..., \alpha_r$ is slice, then for any finite sequence of positive integers $m_1,...,m_k$ we have $
u^{(m_1,...,m_k)}(\alpha_1)+...+ u^{(m_1,...,m_k)}(\alpha_r)=0$. 
   
The number
$sg(\alpha_1,...,
\alpha_r)$  does not depend on the order in the tuple $\alpha_1,..., \alpha_r$. This number is preserved if we replace  
$\alpha_1,..., \alpha_r$ with  cobordant strings. If $\alpha_r$ is slice, then $sg(\alpha_1,..., \alpha_r)=sg(\alpha_1,...,
\alpha_{r-1})$.
Reversing orientations in 3-manifolds and/or surfaces $F$, we obtain  
$$sg(\overline \alpha_1,..., \overline \alpha_r)=sg(\alpha^-_1,..., \alpha^-_r)=sg(\alpha_1,..., \alpha_r).$$
 Using the gluing as  in the proof of
Lemma
\ref{ldd:nbl}, we  obtain that for any $0\leq s\leq r$ and any strings $\alpha, \alpha_1,..., \alpha_r$,
$$sg(\alpha_1,..., \alpha_r)\leq  sg(\alpha_1,..., \alpha_s, \alpha)
+ sg(\overline \alpha^-, \alpha_{s+1},..., \alpha_r).$$
When $\alpha$ is a trivial string, this gives the obvious inequality 
$sg(\alpha_1,..., \alpha_r)\leq  sg(\alpha_1,..., \alpha_s)
+ sg(\alpha_{s+1},..., \alpha_r)$.

For $r=2$, we can rewrite the slice genus    in the equivalent form
$sg'(\alpha, \beta)= sg(\alpha, \overline   \beta^-)= sg(\alpha^-, \overline   \beta)$ for   strings $\alpha, \beta$.    The
results of the previous paragraph imply  that    the number  $sg'(\alpha, \beta)$ depends only on the cobordism classes of
$\alpha, \beta$ and     defines a metric on the set of cobordism classes of strings: 
$sg'(\alpha, \beta)=0$   if and only if
$\alpha$ and $\beta$ are cobordant  (cf. the
end of the proof of Theorem
\ref{th:0388}),    $sg'(\alpha,
\beta)=sg'(  \beta, \alpha)$, and  $sg'(\alpha, \beta)\leq sg'(\alpha,
\gamma)+sg'(\gamma, \beta)$ for any strings
$\alpha, \beta, \gamma$.  Note also that $sg(\alpha)=sg'(\alpha,
O)$ where $O$ is a trivial string and $
sg'(\alpha, \beta)\leq sg(\alpha) +sg(\beta)$.

			  \subsection{Adams operations on strings}\label{sn:nbccc} We can define \lq\lq Adams operations" $\{\psi^n\}_{n\in
\ZZ}$ on  the  set    of homotopy classes of  strings.  Let $\alpha$ be a  virtual string. Replacing $\alpha$ by a
homeomorphic string, we can identify its core circle with   
$S^1=\{z\in 
\CC\,\vert \,  \vert z\vert=1\}$. Consider  a   curve
$\omega:S^1\to
\Sigma$ realizing $\alpha$ on a surface $\Sigma$.  The   mapping  $  S^1\to \Sigma$ sending 
$z\in S^1$  to $ \omega(z^n)$    is homotopic to a generic curve $   S^1\to \Sigma$. We 
define 
$\psi^n(\alpha) $ to be the homotopy class of its underlying string. Lemma \ref{l:nbl} implies that   
$\psi^n  (\alpha)
$  depends neither  on the choice of $\omega$ nor on the choice of $\alpha$ in its homotopy class. Clearly,
$\psi^{1} (\alpha)=\alpha$, $\psi^{-1} (\alpha)=\alpha^{-}$, and  $\psi^{mn}=\psi^m\circ 
\psi^n$ for any  $m,n\in \ZZ$.  It is  also clear that  $\psi^n 
$ transforms cobordant strings into cobordant strings and induces thus an \lq\lq Adams operation"  
on  the  set    of cobordism classes of  strings. As an exercise, the reader may check that
$u(\psi^n(\alpha))=\sign (n)\, n^2\, u(\alpha)$.

 \subsection{Remarks.}\label{sn:fvvvfd6}   1.  J. S. Carter
\cite{ca2}  first   observed that there are closed curves on    surfaces that bound no  singular  disks in
3-manifolds bounded by these surfaces.  One of his results can be rephrased by saying that the string
$\alpha_{2,1}$ is not slice.    Carter's technique consists in studying  certain partitions  (called filamentations) of the set of
double points of a curve into pairs and singletons.  In the notation of the proof of Lemma  \ref{l:nblddddd} (where $F$
should be a disk), a filamentation  is formed by the orbits of the involution
$\mu$.  Carter's  obstruction to the sliceness is formulated in terms of intersection numbers of  the  intervals of double points
on the disk.     Note also a relevant result of \cite{hk} (Theorem 4.10):  if a closed curve on a surface has a
filamentation then any homotopic  closed curve    also has a filamentation. 
 
2. The definition of cobordism for strings requires only the {\it existence} of  
disjoint  realizations homotopic in a 3-manifold.  Note  that {\it any} two disjoint  curves realizing
  cobordant strings   on a closed surface are homotopic in a certain  oriented 3-manifold bounded by this surface.  To see
this, one needs the following  two observations.  

(i) For any disjoint  realizations 
$\omega_1
$,
$\omega_2 $  of   strings $\alpha_1, \alpha_2$ on a closed surface  $\Sigma$,  there are 
realizations 
$\omega'_1
$,
$\omega'_2 $  of    $\alpha_1, \alpha_2$ on disjoint closed surfaces  $\Sigma_1, \Sigma_2$ and  an
oriented 3-manifold $M$ with $\partial M=\Sigma  \cup (-\Sigma_1)\cup (-\Sigma_2)$ such that $\omega_i$ is homotopic
to
$\omega'_i$ in
$M$ for $i=1,2$. This can be proven by adding 2-handles along the simple closed curves  in $\Sigma$  bounding a regular
neighborhood  of
$\omega_1$   in $\Sigma$.

 (ii) For any   realizations 
$\omega
$,
$\omega' $  of the same string on closed surfaces $\Sigma$, $\Sigma'$,  there is an
oriented 3-manifold $M$ with $\partial M=\Sigma  \cup (-\Sigma')$ such that $\omega$ is homotopic to $\omega'$ in $M$.
This can be deduced from the fact that  both  $\omega, \omega'$ can be obtained from  the canonical      realization  of  the
string on a closed surface  of minimal genus by adding 1-handles.

  \section{Based     matrices of strings}

\subsection{Based     matrices}\label{fi:g51}  Fix an abelian group $H$.  A  {\it based skew-symmetric matrix over 
$H$}  or shortly a {\it based matrix}  is a triple $(G, s, b:G^2=G\times G \to H)$ where   $G$ is a finite set, $s\in G$, 
and  the mapping
$b $ is   skew-symmetric   in the sense that 
$b(g,h)=-b(h,g)$ for all $g,h\in G$  and  $b(g,g)=0$ for all 
$g\in 
G$.

We   call an element $g\in G-\{s\}$ {\it annihilating} (with respect to 
$b$) if 
$b(g,h)=0$ for all $h\in G$. We call $g\in G-\{s\}$ a {\it core 
element} if 
$b(g,h)=b(s,h)$ for all $h\in G$.   We   call two elements $g_1,g_2\in 
G-\{s\}$ 
{\it complementary} if $  b (g_1, h)+ b
 ( g_2,h)=  b(s,h)$ for all $h\in  G$.  A based matrix 
$( G, s, b)$ is  {\it primitive}    if it has no annihilating elements, no core 
elements, and 
no complementary pairs of elements. An example of a primitive based matrix 
is 
provided by the {\it trivial based matrix} $(G,s, b)$ where $G$ consists of 
only one 
element $s$ and $b(s,s)=0$.

  We define 
three 
operations $M_1, M_2, M_3$ on based matrices, called {\it elementary 
extensions}. They add to a  based matrix $(  G, s, b)$ an 
annihilating 
element, a core element, and  a pair  of complementary elements, 
respectively. 
More precisely, $M_1$ transforms $(  G, s, b)$ into the (unique) based matrix $(  
\overline G=G\amalg \{g\},  s, \overline b)$
such that    $\overline b:\overline G\times \overline G \to H$  extends
$b$ and  $\overline b (g,h)= 0$ for all $h\in \overline G$.  
The move  $M_2$  transforms $(  G, s, b)$ into the (unique) based matrix $(  
\tilde 
G=G\amalg \{g\},   s, \tilde b)$
such that    $\tilde  b:\tilde  G\times \tilde  G \to H$  
extends $b$ and  $\tilde  b (g,h)=\tilde  b (s,h)$ for all $h\in 
\tilde G$. 
 The   move  $M_3$  transforms $(  G, s, b)$ into a   based matrix 
$(  \hat 
G=G\amalg \{g_1,g_2\},  s, \hat b)$
where   $\hat b:\hat G\times \hat G \to H$ is any skew-symmetric   map   
extending 
$b$ and such that $\hat b (g_1, h)+\hat b
 ( g_2,h)=  \hat  b(s,h)$ for all $h\in \hat G$.   It is clear that a based matrix 
  is primitive  if  and only if it cannot be obtained from another 
based matrix by  an elementary extension.  

 Two based   matrices $(  G, s, b)$ and $(  G', s', b')$ are {\it 
isomorphic} 
if there is a bijection $G\to G'$  sending $s$ into $s'$  and  
transforming $b$ 
into $b'$.  To specify the isomorphism class of a   based matrix  $(G, s, b)$, it suffices  to specify the
matrix
$(b(g,h))_{g,h\in G}$ where  it is understood that the first  column and row correspond to $s$. In this way every
skew-symmetric square matrix over $H$ (with zeroes on the diagonal) determines a  based matrix.

Two based matrices are {\it   homologous}
if one can be obtained from   the other by a finite sequence of 
elementary 
extensions $M_1,M_2, M_3$, the inverse transformations, and 
isomorphisms.   The homology is an equivalence relation on the set of based matrices.

\begin{lemma}\label{l:t51} Every based matrix is  obtained from   a 
primitive based matrix by elementary extensions. Two homologous primitive based matrices are isomorphic.
		   \end{lemma}
                     \begin{proof}  The first claim is obvious: 
eliminating    
annihilating elements,   core elements, and   complementary pairs of 
elements by 
the moves $M_i^{-1}$ with $i=1,2,3$ we can  transform any 
based matrix $T$ into a  primitive based matrix $T_\bullet$.  Then  $T$ is 
obtained from 
$T_\bullet$ by elementary extensions. 
  
  To prove the second claim, we need the following assertion:
  
  $(\ast)$ a  move $M_i$ followed by   $M_j^{-1}$ yields the same 
result as    
an isomorphism, or a   move $ M_k^{\pm 1}$, or     a  move $M_k^{-1}$ 
followed 
by $M_l$ with $k,l\in \{1,2, 3\}$.

  This assertion will imply the second claim of the lemma. Indeed,    
suppose that two primitive  based matrices $T,T'$ are related by a finite 	
sequence 
of transformations $M_1^{\pm 1}, M_2^{\pm 1}, M_3^{\pm 1}$ and 
isomorphisms. 	
 An isomorphism of based matrices followed by $ M_i^{\pm 1}$     
can be also 
obtained as  $ M_i^{\pm 1}$    followed by an isomorphism.  Therefore 
all 
isomorphisms in our sequence   can be   accumulated at the end. The 
claim  
$(\ast)$ implies that $T,T'$
can be  related by a finite sequence of moves consisting of  several 
  moves of type $  M_i^{-1}$ followed by  several 
  moves of type $ M_i$ and   isomorphisms. However, since   $T$ is 
primitive  
  we cannot apply to it a move of type $  M_i^{-1}$. Hence there are no 
such 
moves in our sequence. Similarly, since $
  T'$ (and any isomorphic based matrix) is primitive, it cannot be obtained by 
an 
application of $M_i$. Therefore our sequence  consists solely of 
isomorphisms so 
that $T$ is isomorphic to $T'$.

 Let us now prove $(\ast)$. We have to consider  nine  cases depending on 
$i,j\in \{1,2,3\}$. 
 
  For $i, j \in \{1,2\}$, the move   $M_i$ on a  based matrix 
$(G,s,b)$ adds 
one element $g$   and then  $M_j^{-1}$ removes   one  element  $g'\in G\amalg \{g\}$. If 
$g'=g$, 
then   $  M_j^{-1}\circ M_i$   is the identity. If   $g'\neq g$, then 
$g'\in G$ 
is annihilating (resp. core) for $j=1$ (resp. $j=2$). The 
transformation $  
M_j^{-1}\circ M_i$   can be achieved by first applying $M_j^{-1}$ that 
removes 
$g'$ and then applying $M_i$ that adds $g$.

   Let $i=1, j=3$. The move 
 $M_i$ on $(G,s,b)$ adds an annihilating element $g$  and  
  $M_j^{-1}$ removes two complementary  elements $g_1,g_2\in G\amalg \{g\}$. 
  If $g_1\neq g$ and $g_2\neq g$, then $g_1,g_2\in G$  and $  
M_j^{-1}\circ M_i$ 
 can be achieved by first removing $g_1,g_2$ and then adding $g$. If 
$g_1=g$, 
then   $g_2 $ is a core element of $G$ and $  M_j^{-1}\circ M_i$ is the 
move 
$M_2^{-1}$    removing $g_2$. The case $g_2=g$ is similar.
  
   Let $i=2, j=3$. The move 
 $M_i$ on $(G,s,b)$ adds a core element $g$  and  
  $M_j^{-1}$ removes two complementary  elements $g_1,g_2\in G\amalg \{g\}$. 
  If $g_1\neq g$ and $g_2\neq g$, then  $  M_j^{-1}\circ M_i$  can be 
achieved 
by first removing $g_1,g_2$ and then adding $g$. If $g_1=g$, then   
$g_2\in G$ 
is   an annihilating element of $G$ and $  M_j^{-1}\circ M_i$ is the 
move 
$M_1^{-1}$    removing $g_2$. The case $g_2=g$ is similar.

   Let $i=3, j=1$. The move 
 $M_i$  on $(G,s,b)$ adds two complementary  elements $g_1,g_2$   and  
$M_j^{-1}$ removes an annihilating element  $g\in G\amalg \{g_1,g_2\}$. 
  If $g\neq g_1$ and $g\neq g_2$, then $g\in G$ and $  M_j^{-1}\circ 
M_i$  can 
be achieved by first removing $g$ and then adding $g_1,g_2$.
   If $g=g_1$, then   $ g_2$ is a core element of $G\amalg \{g_2\}$   
and $  
M_j^{-1}\circ M_i=M_2$. The case $g=g_2$ is similar.
  
   Let $i=3, j=2$. The move 
 $M_i$  on $(G,s,b)$ adds two complementary  elements $g_1,g_2$   and  
$M_j^{-1}$ removes a core element  $g\in G\amalg \{g_1,g_2\}$. 
  If $g\neq g_1$ and $g\neq g_2$, then $g\in G$ and  $  M_j^{-1}\circ 
M_i$  can 
be achieved by first removing $g$ and then adding $g_1,g_2$.
   If $g=g_1$, then   $ g_2$ is  an annihilating  element of $G\amalg 
\{g_2\}$   
and $  M_j^{-1}\circ M_i=M_1$. The case $g=g_2$ is similar.

      Let $i=j=3$. The move 
 $M_i$  on $(G,s,b)$ adds two complementary  elements $g_1,g_2$   and  
$M_j^{-1}$ removes two  complementary elements  $g'_1,g'_2\in G\amalg \{g_1,g_2\}$. If these 
two pairs  
 are disjoint, then 
  $  M_j^{-1}\circ M_i$  can be achieved by first removing  $g'_1,g'_2\in G$ 
and then 
adding $g_1,g_2$. If these two pairs  
  coincide, then $  M_j^{-1}\circ M_i$  is the identity. It remains to 
consider 
the case where these pairs have one common element, say $g'_1=g_1$, 
while 
$g'_2\neq g_2$. Then $g'_2\in G$ and for all $h\in G$,
  $$\hat b(g_2,h)= \hat  b(s,h)-\hat b(g_1,h)= \hat  b(s,h)-\hat b(g'_1,h)=\hat 
b(g'_2,h)= 
b(g'_2,h).$$
 Therefore the move   $  M_j^{-1}\circ M_i$ gives a based matrix isomorphic 
to 
$(G,s,b)$. The isomorphism $G\to (G-\{g'_2\}) \cup \{g_2\}$ is the 
identity on 
$G- \{g'_2\}$ and  sends $g'_2$ into $g_2$.
\end{proof}

Lemma \ref{l:t51}  implies that each  based matrix $T=(G,s,b)$ is homologous to a   primitive    based matrix
 $T_\bullet=(G_\bullet, s_\bullet, b_\bullet)$
unique up to  isomorphism.  This  reduces   classification   of 
based matrices up to homology to a classification of primitive based matrices up to   isomorphism. Note that   we can
choose 
$T_\bullet $ in its isomorphism class so that $G_\bullet\subset 
G$ and 
$b_\bullet$ is the restriction of $b$ to $G_\bullet\times G_\bullet$.

 We define  two more operations on based   matrices.   For  a  based matrix $T=(G,s,b)$, set  $-T= (G,s,-b)$ and 
$T^-=(G,s, b^-)$  where     $b^- (s,h)=-b( s,h)$, $b^- (h,s)=-b( h,s) $  for all $h\in G$ and
$b^-(g,h)= b(g,h)- b(g,s)-b(s,h) $ for all
$g,h\in G-\{s\}$.   The transformations  $T\mapsto -T$, $ T\mapsto T^-$ are commuting involutions on the set
of based matrices.   It is easy to  check that they are compatible
with   homology and preserve the class of   primitive based matrices.
It follows from the definitions that   $(-T)_\bullet =- T_\bullet$ and $(T^-)_\bullet=(T_\bullet)^-$.

\begin{remars}\label{fi:g52} 1.  The moves $M_1, M_2, M_3$ on
based matrices are 
not independent. It is easy to present $M_2$ as a composition of $M_3$ 
with 
$M_1^{-1}$.			 

2. Each   isomorphism invariant $v$ of primitive  
based matrices  
 extends to a  homology        invariant   of based matrices by  $v(T)=v(T_\bullet)$.   
The most important numerical  invariant of a primitive based matrix 
$(G,s,b)$ 
is the number  $\#(G)$. It is   easy to define further  
invariants  of primitive  based matrices. For instance, for   $k\in H$,  we can 
set
$$v_k(G,s,b)= \#\{g\in G\,\vert\, b(g,s )=k\}.$$
Similarly, for $k\in H$ and   a finite set  $A$ of elements of $ H$ endowed 
with 
non-negative multiplicities,
set 
$$v_{k,A}(G,s,b)=\#\{g\in G\,\vert\,  b(g,s)=k\,\,{\text {and}}\,\,\, 
\{b(g,h)\}_{h\in G-\{s\}}=A\},$$
where the latter equality is understood as  an  equality of sets with 
multiplicities.  Clearly, $v_k=\sum_A v_{k,A}$.

3. If $H\subset \RR$ is a subgroup of the additive group of real numbers, then the 1-variable polynomial
  $$ u(T) \,(t)=\sum_{g\in G, b(g,s)  \neq 0} \sign ( b(g,s) )\, 
t^{\vert  
b(g,s)  \vert} $$
is a homology invariant of
a  based matrix
$T=(G,s,b)$. 
 \end{remars}

\subsection{The based matrix of a string}\label{fi:g53} With each 
virtual 
string $\alpha$ we associate a  based matrix 
$T(\alpha)=(G,s,b)$ over $\ZZ$. Set $G=G(\alpha)=\{s\}\amalg 
\arr(\alpha)$. To 
define $b=b(\alpha):G\times G\to \ZZ$, we identify $G$ with  the basis $s\cup 
\{[e]\}_{e\in \arr(\alpha)} $ of 
$ H_1(\Sigma_\alpha)$, see 
Section \ref{fi:g32}.  The map $b$ is obtained by restricting the 
homological 
intersection pairing
$ H_1(\Sigma_\alpha)\times H_1(\Sigma_\alpha)\to \ZZ$
to $G$. It is clear that $b$ is skew-symmetric. We can   compute 
$b$  combinatorially using   Formula \ref{pik} and Lemma  
\ref{l:t33}. In 
particular,    $b(e,s)=n(e)$ for all $e\in \arr(\alpha)$. 

The map  $b$ can be computed from any closed curve $\omega $ realizing $\alpha$ on a 
surface $\Sigma$. 
Indeed,   such a  curve is  obtained from  the 
canonical realization of $\alpha$ in $\Sigma_\alpha$  
via 
 an orientation-preserving embedding  $\Sigma_\alpha\hookrightarrow \Sigma$. 
It remains   to observe that such an embedding preserves  
intersection numbers    and  transforms the basis $s\cup 
\{[e]\}_{e\in \arr(\alpha)} $ of 
$ H_1(\Sigma_\alpha)$ into the subset $[\omega  ], \{[\omega_x] \}_{x\in \Join (\omega)}$ of $ H_1(\Sigma)$, cf.
Section
\ref{sn:g25}.

 \begin{lemma}\label{th:e53}  If two virtual strings   are homotopic, 
then their 
based matrices   are homologous. \end{lemma}
  \begin{proof} By Lemma \ref{l:nbl} it is enough to show that if  two closed curves $\omega, \omega' $ on a
surface
$\Sigma$ are homotopic, then the based matrices of their underlying strings are homologous.  
By the discussion in Section \ref{sn:g13},  it suffices to consider the  case where $\omega'$ is obtained from
$\omega$ by one of the   local moves listed there. 		

  If $\omega'$ is obtained from $\omega$ by   
adding a small  curl, then $\Join\! \! (\omega')=\,\Join\! \! (\omega)\cup \{y\}$  where $y$ is a new crossing. 
Clearly
$[\omega'_y]=0\in H_1(\Sigma)$  or  $[\omega'_y]=[ \omega']=[\omega] \in H_1(\Sigma)$ depending on
whether  the curl lies on the  right or  on the left  of
$\omega$. Also      $[\omega'_x]=[\omega_x]$  for all $x\in \,  \Join\!\! (\omega)$. 
 Hence   
$T(\beta)$ is obtained from  
$T(\alpha)$ by $M_1$ or $M_2$.

Suppose that  $\omega'$  is obtained from $\omega$ by  the move  pushing a branch of $\omega$ across  another branch
and creating  two new double points $y,z$.  Clearly,   $[\omega'_x]=[\omega_x]$  for all $x\in \, 
\Join\! \!(\omega)\subset\,
\Join\! \!(\omega')$.  It is easy to see that  
 $[\omega'_y]+[\omega'_z]= [ \omega']=[\omega]\in H_1(\Sigma)$.
Therefore 
 $T(\beta)$ is obtained from  $T(\alpha)$ by $M_3$.
 
If  $\omega'$  is obtained from $\omega$ by   pushing  a branch of $\omega$ across a 
double point, then the
  subsets $[\omega  ], \{[\omega_x]\}_{x\in \Join (\omega)}   $ and $[\omega'  ],
\{[\omega'_x]\}_{x\in \Join (\omega')}   $ of $ 
H_1(\Sigma)$ coincide so that   $T(\alpha)$ is isomorphic to $T(\beta)$.  \end{proof} 
  
    \subsection{Invariants of strings from 
based matrices}\label{fi:g54}   
Every virtual string $\alpha$ gives rise 
to a   
primitive based matrix $T_\bullet(\alpha) $  over $\ZZ$ by $T_\bullet(\alpha)=(T(\alpha))_\bullet$. This is the only
primitive based matrix (up to  isomorphism) homologous to
$T(\alpha)$.   By
 Lemma  \ref{th:e53}, the based matrix  
$T_\bullet(\alpha)=(G_\bullet, s_\bullet, b_\bullet)$   is a 
  homotopy invariant of    $\alpha$.  This based matrix  determines the polynomial $u(\alpha)$ introduced in Section
\ref{abric}:  it follows from    Formulas \ref{homr11} and  \ref{pik} that 
  $  u(\alpha)=u(T(\alpha))=u(T_\bullet(\alpha)) $.  The number
$\rho (\alpha) =\#  (G_\bullet ) -1$ 
is a useful  homotopy invariant of $\alpha$
  which may be non-zero even when $u(\alpha)=0$, cf.  the 
examples 
below.  Note that  if    $\alpha$ is homotopically trivial, 
then 
$T_\bullet(\alpha)$ is a trivial based matrix and $\rho(\alpha)=0$.

 It follows from the definitions that $T(\alpha^-)= (T(\alpha))^-$ and therefore  $T_\bullet(\alpha^-)= (T_\bullet(\alpha))^-$.
Similarly, $T(\overline \alpha)= -(T(\alpha))^-$ and    $T_\bullet(\overline \alpha)= -(T_\bullet(\alpha))^-$.

   The based matrix  $T_\bullet(\alpha)=(G_\bullet, s_\bullet, b_\bullet)$ can be used to estimate the 
homotopy 
rank and the  homotopy genus of $\alpha$. Namely,
   $hr(\alpha)\geq \rho(\alpha)$ since any string homotopic to $\alpha$ 
must 
have at least $\rho(\alpha)$ arrows. Similarly,  $hg(\alpha)\geq (1/2) 
\rank 
b_\bullet$ where   $\rank 
b_\bullet$  is     the 
rank of the integral matrix $(b_\bullet(g,h))_{g,h\in G_\bullet}$.   Indeed, if $\alpha'$ is a string homotopic to $\alpha$ and $T 
(\alpha')=(G',s',b')$, then  
   $g(\alpha')= (1/2) \rank b' \geq (1/2) \rank b_\bullet$ since the matrix 
of $b'$ 
contains the matrix of $b_\bullet$ as a submatrix. 

Combining the  inequalities   $hr(\alpha)\geq \rho(\alpha)$, 
$hg(\alpha)\geq (1/2) 
\rank 
b_\bullet$  with the obvious inequalities $\rank  \alpha \geq
hr(\alpha)$ and
$g(\alpha)\geq hg(\alpha)$, we obtain that  if  
$T(\alpha)  $ is  primitive, then  $hr(\alpha) =\rank  \alpha$ and 
$hg(\alpha)  =g(\alpha)$.

     \subsection{Applications}\label{fi:g55}  (1) The based matrix 
$T(\alpha_{p,q})$ 
of the string   $\alpha_{p,q}$ with $p,q\geq 1$ was computed in Section 
\ref{fi:g33}. It is easy to check   that except in the case $p=q=1$, this 
based matrix  
is primitive. Thus $ T_\bullet(\alpha_{p,q})=T(\alpha_{p,q})$,  $hr(\alpha_{p,q})=\rank \alpha_{p,q}=p+q$  
and  $hg(\alpha_{p,q})=
 g(\alpha_{p,q})$ provided $p\neq 1$ or
$q\neq 1$.  In particular,   $\alpha_{p,p}$
  is a homotopically non-trivial  string with zero $u$-polynomial 
for all 
$p>1$. 
 
 (2)  The product of strings defined in Section \ref{sn:g1335} does not
induce  a well-defined operation on the set of homotopy classes of strings.  To see this, we exhibit a homotopically
non-trivial  string which is   
a  product of two copies of  a  homotopically trivial string.  Namely, the string   $\alpha_\sigma$ considered in Section
\ref{rered}  has the required properties.   It is observed there  that $\alpha_\sigma$   is a product of two copies of 
the    homotopically trivial string $\alpha_{1,1}$.  The based matrix  $T({\alpha_\sigma})$
  can be explicitly computed, cf. Section    \ref{fi:g33}.2.   It is determined by the
following  skew-symmetric   matrix over $\ZZ$:  
$$
 \left [ \begin{array}{ccccc}  0& -1&1& -1& 1 \\
        1& 0& 1& -1& 1 \\
-1& -1& 0& -1&1 \\
1& 1&1& 0&1 \\
-1& -1& -1& -1& 0 
\end{array} \right ].$$
It is easy to check that this based matrix is primitive. Hence $\alpha_\sigma$ is not homotopically trivial.  Moreover,  it is not
homotopic to a string with   $<4$ arrows.    

(3) We   prove that the involution $\alpha\mapsto \overline \alpha$   acts  non-trivially on
the set of homotopy classes of strings.   Consider the permutation $\sigma=(134)(2)$ on the set $\{1,2,3,4\}$ sending 1 to
3, 3 to 4, 4 to 1, and 2 to 2.  Drawing the  string $ \alpha_\sigma$ we obtain that
$\overline {\alpha}_\sigma=\alpha_\tau$ where $\tau$ is the permutation $(124)(3)$.  The based matrices $T({\alpha_\sigma})$
and
$T(\alpha_\tau)$ can be explicitly computed. They are determined by the
following  skew-symmetric   matrices:  
$$
 \left [  \begin{array}{ccccc}  0& -2&0& -1& 3 \\
        2& 0& 1& 0& 3 \\
0& -1& 0& 0&2 \\
1& 0&0& 0&1 \\
-3& -3& -2& -1& 0 
\end{array} \right ],\,\,\,\,\,\,\,\,\,\,\,\,\,\,\,\,\,\,\,\, \left [  \begin{array}{ccccc}  0& -1&-2& 0& 3 \\
        1& 0& -1& 1& 3 \\
2& 1& 0& 1&2 \\
0&-1&-1& 0&1 \\
-3& -3& -2& -1& 0 
\end{array} \right ].$$
The based matrices   $T({\alpha_\sigma})$ and
$T(\alpha_\tau)$ are not isomorphic; this is clear for instance from the fact that the first matrix has a row with
three zeros while the second matrix does not have such a row. It is clear also that these based matrices are primitive. By Lemma
\ref{l:t51},  they are not homologous. Hence ${\alpha_\sigma}$ is not homotopic to
$\alpha_\tau=\overline {\alpha}_\sigma$.

     \section{Genus and cobordism for  based   matrices}\label{fi:g599912}

Throughout this section the symbol $R$ denotes a domain, i.e., a commutative ring with unit and with no 
zero-divisors. By a  based matrix 
 over $R$, we   mean a based matrix over the additive group of $R$.

\subsection{Genus of   based matrices}\label{fi:g5999}   We define a numerical invariant of a  based matrix 
$T=(G,s,b)$ over $R$ called its {\it genus} and denoted $\sigma (T)$.  For    subsets $X,Y\subset G$, set
$b(X,Y)=\sum_{g\in X, h\in Y} b(g,h)\in R$.  Clearly,  $b(X,Y)=-b(Y,X)$  and  $b(X,X)=b(\emptyset ,X)=0$ for all
$X,Y\subset G$. A  {\it (simple) filling} $\mathcal X$ of
$T$ is  a  finite family $\{X_i\}_i$ of disjoint  (possibly empty)  subsets of  
$G$ such that $\cup_i X_i=G$, $\# (X_i) \leq 2$ for all $i$, and one of $X_i$ is  the one-element set $\{s\}$.  The {\it
matrix} of  $\mathcal X=\{X_i\}_i$   is the  matrix $(b(X_i,X_j))_{i,j}$.  This is a skew-symmetric
square matrix (with zero diagonal) over
$R$. Its rank (the maximal size of a non-zero minor) is   an even non-negative integer; let  $\sigma(\mathcal X)$ denote
half of this rank.  Set 
$\sigma(T)=\min_\mathcal X
\sigma(\mathcal X) $ where $ \mathcal X$ runs over all   fillings of $T$.   Extending $b$ by linearity to the
$R$-module
$\Lambda=RG$ freely generated by $G$ and identifying a subset  $X\subset G$ with the vector $\sum_{g\in X}
g\in 
\Lambda$, we can interpret
$\sigma(T)$ as half the minimal rank of the restriction of $b$  to the submodules  of $\Lambda$ arising from    fillings
of
$T$.

Note that $\sigma (T)\geq 0$ and  $\sigma (T)=0$ if and only if $T$ has a   filling with zero matrix. 
In the latter case we say that  
$T$ is {\it hyperbolic}.

\begin{lemma}\label{l:gg1} The genus of a  based matrix is  a  homology  invariant.
		   \end{lemma}
                     \begin{proof} 
  By Remark \ref{fi:g52}.1, it suffices to prove that $\sigma(T)=\sigma(T')$ for any  based matrix
$T'=(G',s,b')$     obtained from a  based matrix 
$T=(G,s,b)$ by a move
$M_i$ with $i=1,3$.  The set $  G'-G$ consists of one  element if $i=1$ and of two elements if $i=3$.  Pick a   
  filling $\mathcal X=\{X_i\}_i$ of
$T$ such that  $\sigma(T)= 
\sigma(\mathcal X)$.  Consider the     filling ${\mathcal X}'=(G'-G)\cup  \{X_i\}_i$ of $T'$.  Its matrix   is obtained
from the one of
$\mathcal X$ by adjoining a row and a column.  For $i=1$, these row and column are zero so that   $  \sigma
(\mathcal X')
=
\sigma(\mathcal X)$.   For $i= 3$,  we have $b(G'-G,Y)=b(\{s\}, Y)$ for all  $Y\subset G$.  Since one of the sets $X_i$ equals
$\{s\}$,  we again obtain  $  \sigma
(\mathcal X')
=
\sigma(\mathcal X)$.   Hence
$ \sigma(T')\leq \sigma(\mathcal X') =\sigma(\mathcal X)=\sigma(T) $. 

To prove the opposite inequality, pick  a   
  filling $\mathcal X'=\{X_i\}_i$ of
$T'$ such that  $\sigma(T')= 
\sigma(\mathcal X')$. We shall construct  a   
  filling $\mathcal X $ of
$T$ such that  $\sigma(\mathcal X)\leq 
\sigma(\mathcal X')$.  This would imply   $\sigma(T)\leq 
\sigma(\mathcal X) \leq   
\sigma(\mathcal X')= \sigma(T')$.
Consider the case $i=1$.  One of the sets $X_i$ contains the 1-element set
$G'-G$.  We replace this
$X_i$ by $X_i -(G'-G)$ and keep all the other $X_i$. This gives a   filling  $\mathcal X$   of $T$ whose matrix
coincides with the matrix  of $\mathcal X'$.  Hence  $ 
\sigma(\mathcal X) =  
\sigma(\mathcal X') $.  Let now $i=3$.  If  one of the sets $X_i$ is equal to $G'-G=\{g_1,g_2\}$, then  removing this
$X_i$ from
$\mathcal X'$ we obtain a   filling  $\mathcal X$  of $T$.  As in the previous paragraph, 
$\sigma(\mathcal X)=\sigma(\mathcal X')$.     Suppose that the   elements $g_1,g_2$  of $G'-G$ belong to different
subsets, say $X_1, X_2$,  of the   filling
$\mathcal X'$.   Then the sets $X_i$ with $i\neq 1,2$ and $X= (X_1\cup X_2) - \{g_1,g_2\}$ form a    filling of
$T$.  Let $X_0$ be the  term of the   fillings $\mathcal X$ and  $\mathcal X'$ equal to $\{s\}$.  For
any $Y\subset G$,  
$$b(X,Y)=b'(X,Y)=b'(X_1,Y)+b'( X_2,Y)- b'(g_1,Y)- b'(g_2,Y)=b'(X_1,Y)+b'( X_2,Y)- b'(X_0,Y).$$ 
Applying this to $Y=X_i$ with $i\neq 1,2$,  we obtain that 
 the skew-symmetric  bilinear form  determined by the  matrix
of $\mathcal X$ is induced from  the  skew-symmetric   bilinear form  determined by the  matrix
of $\mathcal X'$ via the  linear map of the corresponding free $R$-modules sending the basis vectors $X$ and
$\{X_i\}_{i\neq 1,2}
$ respectively to 
$ X_1+X_2-X_0$ and $\{X_i\}_{i\neq 1,2}$.  Hence  
$\sigma(\mathcal X)
\leq   
\sigma(\mathcal X')$.
\end{proof}

\begin{corol} \label{aaal:gg25} For any   based matrix  $T$ over $R$, we have $\sigma(T_\bullet)=\sigma (T)$.
A   based matrix  over $R$  homologous to a hyperbolic based matrix is itself hyperbolic.
		   \end{corol}

\subsection{Genus for tuples of based matrices}\label{uussstrg52} The definition of the genus of a based matrix can be
extended to tuples of based matrices. Consider a tuple of $r\geq 1$ based matrices 
$T_1=(G_1,s_1,b_1),..., T_r=(G_r,s_r,b_r)$ over $R$.    Replacing $T_1,...,T_r$ by isomorphic based
matrices, we can assume that the sets
$G_1,...,G_r$ are disjoint.  Let $\Lambda=RG$ be the free $R$-module with basis 
$G=\cup_{t=1}^r G_t$. Let $\Lambda_s$ be the submodule of  $\Lambda$ generated by $s_1,...,s_r$.   We call a vector
$x\in
\Lambda$ {\it short} if   $x\in \Lambda_s$  or $x\in g+
\Lambda_s$ for some  $g\in 
G-\{s_1,...,s_r\} $  or $x \in g+h+
\Lambda_s$ for   distinct   $g,h \in G-\{s_1,...,s_r\}$.   
  A  {\it filling} of
$T_1,...,T_r$ is  a finite family $  \{\lambda_i\}_i$ of short vectors in $\Lambda$ such that $\sum_i
\lambda_i=\sum_{g\in G} g (\modu  \Lambda_s)$ and one of $\lambda_i$ is equal to $s_1+s_2+...+s_r$.
Note that each element of  $G-\{s_1,...,s_r\}$ appears in exactly one $\lambda_i $ with non-zero coefficient; this 
 coefficient is then
$+1$. The basis vectors
$ s_1,...,s_r $  may appear in several $\lambda_i$ with non-zero coefficients.  

The maps $\{b_t:G_t\times G_t\to R\}_t$ induce a skew-symmetric bilinear form $b=\oplus_t b_t:\Lambda\times
\Lambda\to
R$ such that 
$b(g,h)=b_t(g,h)$ for 
$g,h\in G_t$  and 
$b(G_t,G_{t'})=0$ for $t\neq t'$.   The {\it matrix}
of a filling
$\lambda=
\{\lambda_i\}_i$ of
$T_1,...,T_r$  is the  matrix
$(b(\lambda_i,\lambda_j))_{i,j}$.  This is a skew-symmetric square matrix over
$R$. Let  $\sigma(\lambda) \in \ZZ$ be half of its  rank.  Set  $\sigma(T_1,...,T_r)=\min_\lambda
\sigma(\lambda) $ where $\lambda$ runs over all fillings of $T_1,...,T_r$.  
Clearly $\sigma (T_1,...,T_r)\geq 0$ and  $\sigma (T_1,...,T_r)=0$ if and only if $(T_1,...,T_r)$ has a filling with
zero matrix.  In the latter case we call the sequence 
$T_1,...,T_r$  {\it hyperbolic}.

It is obvious that the genus
$\sigma(T_1,...,T_r)$  is preserved when   $T_1,...,T_r$ are permuted or 
replaced with isomorphic based matrices.   Also
$\sigma(-T_1,...,-T_r)=\sigma(T_1,...,T_r) $.  If
$T_r$ is a trivial based matrix, then
$\sigma(T_1,..., T_r)=\sigma(T_1,...,
T_{r-1})$ (because then the vector $s_r\in \Lambda$ lies in the annihilator of $b$).

For $r=1$, the notion of  a filling   is slightly wider than the notion of a simple filling  in  Section
\ref{fi:g5999}. However,  they give the same genus   and the same set of hyperbolic based matrices.

\begin{lemma} \label{ccxf} For any $1\leq t\leq r$ and any based matrices  $T_0, T_1,..., T_r$,
$$\sigma(T_1,..., T_r)\leq  \sigma(T_1,..., T_t, T_0)
+ \sigma(-T_0, T_{t+1},..., T_r).$$
\end{lemma}
\begin{proof} Consider for concreteness the case where $t=1$ and $r=2$, the general case is quite similar.  
We must prove that $\sigma(T_1, T_2)\leq  \sigma(T_1,  T_0)
+ \sigma(-T_0, T_{2})$.  Let
$T_i=(G_i,s_i,b_i)$ for $i=0,1,2$ and   $T'_0=(G'_0,s'_0,b'_0)$ be a  copy of $T_0$ where 
$G'_0=\{g' \,\vert \, g\in G_0\}$, $s'_0=(s_0)'$,
and $b'_0$ is defined by $b'_0(g',h')=b_0(g,h)$ for $g,h\in G_0$.   We can assume
that the sets
$G_1, G_0,G'_0, G_2$ are disjoint.  Let  $\Lambda_1 , \Lambda_0, \Lambda'_0, \Lambda_2 $ be  free $R$-modules
freely generated by
$G_1, G_0, G'_0,  G_2$, respectively, and let $\Lambda=  \Lambda_1\oplus  \Lambda_0\oplus \Lambda'_0\oplus 
\Lambda_2 $.
There is  a  unique skew-symmetric bilinear form $B=b_1\oplus b_0\oplus (-b'_0)\oplus b_2$ on $\Lambda$ such that the
sets
$G_1, G_0, G'_0, G_2\subset \Lambda$ are mutually orthogonal and the restrictions of $B$ to these subsets  are equal
 to $b_1, b_0, -b'_0,b_2$, respectively.

 Let $\Phi$   be the submodule of
$ \Lambda_0\oplus \Lambda'_0$ generated by    the vectors
$\{g+g'\}_{g\in G_0}$.  Set $L= \Lambda_1\oplus  \Phi \oplus \Lambda_2\subset \Lambda$.     Observe that the
projection $p:
L  \to \Lambda_1\oplus    \Lambda_2$  along $\Phi$ transforms $B$ into $b_1\oplus b_2$. 
Indeed, for any  $g_1,h_1\in G_1, g , h\in G_0, g_2,h_2\in G_2$, 
$$B(g_1+g  +g' +g_2, \,\,h_1+ h + h'+h_2)$$
$$=
b_1(g_1 ,h_1)+ b_0 (g ,  h)+ (-b'_0) (g', h')+ b_2(g_2,h_2)= 
b_1(g_1 ,h_1)+  b_2(g_2,h_2).$$

Pick    a filling  $ \{\lambda_i\}_i\subset \Lambda_1 \oplus \Lambda_0$ of
$(T_1, T_0) $ whose matrix has rank   $2\sigma(T_1,  T_0)$.  This means that   the restriction of $B$ to the
submodule $V_1\subset \Lambda_1\oplus \Lambda_0$ generated by $ \{\lambda_i\}_i$ has rank   $2\sigma(T_1,  T_0)$. 
Similarly, pick      a filling
$ 
\{\varphi_j\}_j\subset
\Lambda'_0\oplus \Lambda_2$    of
$(-T'_0, T_2) $ such that 
  the restriction of $B$ to the
submodule $V_2\subset \Lambda'_0\oplus \Lambda_2$ generated by $ \{\varphi_j\}_j$ has rank  $2\sigma (-T'_0, T_2)$.
  We claim that there is a  finite  set $\psi\subset (V_1+V_2)\cap L$ such that  $p(\psi) \subset
\Lambda_1\oplus    \Lambda_2$ is a filling of   
$(T_1, T_2) $. Denoting by
$V$ the submodule of   
$\Lambda_1\oplus    \Lambda_2$ generated by $p(\psi)$, we obtain  then   the desired inequality:
$$\sigma(T_1, T_2)\leq  \sigma (p(\psi))= (1/2)  \rk ((b_1\oplus b_2) \vert_V)= (1/2)  \rk (B \vert_{p^{-1}(V)}) \leq 
(1/2) 
\rk (B
\vert_{(V_1+V_2)\cap L} )
$$
$$\leq   (1/2) \rk (B
\vert_{V_1+V_2})=  (1/2)  \rk (B \vert_{V_1}) + (1/2) \rk (B \vert_{V_2})=\sigma(T_1,  T_0)
+ \sigma(-T_0, T_{2}).$$
Here the second inequality follows from the inclusion $p^{-1}(V)\subset  (V_1+V_2)\cap L + \Ker p\subset L$ and the fact
that
$\Ker p=\Phi$ lies in the annihilator of $B\vert_L$.

To construct $   \psi$, we modify     $ \{\lambda_i\}_i$  as
follows.  Let $\lambda_1$ be the vector  of this filling equal to $s_1+s_0$.  Adding appropriate 
multiples of
$\lambda_1$  to other  
$\lambda_i$ we can ensure that the basis vector $s_0 \in G_0 $ appears in all $\{\lambda_i\}_{i\neq
1}$  with coefficient 0.  This transforms  $ \{\lambda_i\}_i$ into a new filling of  $(T_1, T_0) $ which will be
from now on denoted    $\lambda= \{\lambda_i\}_i$. This transformation does not change the module
$V_1$ generated by $ \{\lambda_i\}_i$.   Similarly, we can assume that  a vector  $\varphi_1$ of the filling
$\varphi=\{\varphi_j\}_{j }$ is equal to 
$s'_0+s_2$ and   the basis vector $s'_0 \in G'_0 $ appears in all
$\{\varphi_j\}_{j\neq
1}$ with coefficient 0,

The
filling $\lambda$ gives rise to a  1-dimensional manifold     $\Gamma_\lambda$ with boundary $(G_1\cup G_0)-
\{s_1,s_0\}$.        Each
$\lambda_i $ having the form 
$g+h  (\modu R s_1)$  with $g,h\in  (G_1\cup G_0)-
\{s_1,s_0\}$ gives rise to a component of
$\Gamma_\lambda$ homeomorphic to $[0,1]$ and connecting   $g$ with $h$.  Each $\lambda_i $ having
the form 
$g  (\modu R s_1 )$ with $g \in  (G_1\cup G_0)-
\{s_1,s_0\}$ gives rise to a component of
$\Gamma_\lambda$ which is a copy of $[0, \infty)$  
where $0$ is identified with  $g$.  Other $\lambda_i$ and in particular $\lambda_1$ do not contribute to
$\Gamma_\lambda$.  The definition of a filling implies that
$\partial
\Gamma_\lambda= (G_1\cup G_0)-
\{s_1,s_0\}$.   Similarly,   the filling
$\varphi$ gives rise to a    1-dimensional manifold 
$\Gamma_\varphi$ with boundary 
$(G'_0\cup G_2)-\{s'_0, s_2\} $. We can assume that  $\Gamma_\lambda$ and 
$\Gamma_\varphi$ are disjoint.   Gluing  $\Gamma_\lambda$ to
$\Gamma_\varphi$ along the canonical identification
$G_0-\{s_0\} \to G'_0-\{s'_0\}, g\mapsto g'$, we obtain a 1-dimensional manifold, $\Gamma$, with $\partial
\Gamma=(G_1-\{s_1\})\cup  (G_2-\{s_2\})$.  Each component  $K$ of $\Gamma$ is glued from several components of 
$\Gamma_\lambda\amalg  \Gamma_\varphi$ associated with certain vectors  $\lambda_i\in V_1\subset  \Lambda_1\oplus
\Lambda_0
\subset \Lambda$ and/or $ \varphi_j
\in V_2\subset \Lambda'_0\oplus \Lambda_2 \subset \Lambda$. Let   $\psi_K\in
\Lambda$ be the sum of these vectors.  Observe that $\psi_K\in (V_1+V_2)\cap L$; the inclusion  $\psi_K\in
L$ follows from two facts: 
 (i) each point of $K\cap (G_0-\{s_0\})\approx K\cap (G'_0-\{s'_0\})$ is adjacent to   one component of  
$\Gamma_\lambda$ and 
 to one component of   $\Gamma_\varphi$ and (ii) $s_0  $ does not show up   in  
$\{\lambda_i\}_{i\neq 1}$  and   $s'_0   $ does not show up   in  
$\{\varphi_j\}_{j\neq
1}$. 
Set
$\psi_1=\lambda_1+\varphi_1=s_1+s_0+s'_0+ s_2
\in
\Lambda$. Clearly, $\psi_1 \in  (V_1+V_2)\cap L$.    
Set
$\psi=\{\psi_1\}\cup \{\psi_K\}_K$ where $K$ runs over the  components of $\Gamma$ with non-void boundary. 
 Let us check that $p(\psi) \subset
\Lambda_1\oplus    \Lambda_2$ is a filling of   
$(T_1, T_2) $.  Observe that for a compact component $K$ of $\Gamma$ with endpoints $g,h\in G_1\cup G_2$, we have 
$p(\psi_K)=g+h \, (\modu R s_1+R s_2)$.  For a non-compact component $K$ of $\Gamma$ with  one endpoint 
$g \in G_1\cup G_2$, we have 
$p(\psi_K)=g  (\modu R s_1+R s_2)$.  Thus all  vectors in the family $\psi $ are short and their sum is equal to
$\sum_{g\in G_1\cup G_2} g  (\modu R s_1+R s_2)$.  Also $p(\psi_1)=s_1+s_2$.  This means that $p(\psi)$ is a
filling of 
$(T_1, T_2)
$ so that $\psi$ satisfies all the required conditions.
\end{proof}

\subsection{Cobordism of based matrices}\label{u125452} Two based matrices $T_1,T_2$ over $R$  are {\it cobordant} if
$\sigma(T_1,-T_2)=0$.

\begin{theor} \label{859s}  (i) Cobordism is an equivalence relation on the set of isomorphism classes of based matrices. 

(ii) Homologous based matrices are cobordant.

(iii)  The genus   of a tuple of based matrices  is a cobordism invariant.

(iv)  A based matrix   is cobordant to a  trivial based matrix if and only if  it  is hyperbolic. 
\end{theor}
\begin{proof} (i)   For a  based matrix $T=(G,s,b)$, the based matrix $-T$ is isomorphic to the triple $ (G',s',b')$  where
$G'=\{g'\, \vert \, g\in G\}$ is a  disjoint copy of $G$ and  $b'(g',h')=-b(g,h)$ for any $g,h\in G$.  Consider  the filling
$\{g+g'\}_{g\in G}$ of  the pair $(T, -T)$.  The matrix of this filling is
0. Therefore 
$\sigma(T,-T)=0$ so that $T$ is cobordant to itself.  The symmetry of cobordism follows from the equalities
$\sigma(T_2,-T_1)=\sigma (-T_2,T_1)=\sigma(T_1,-T_2)$.  The transitivity  of cobordism follows from the inequalities
$$0\leq \sigma (T_1,-T_3)\leq  \sigma (T_1,-T_2)+\sigma (T_2,-T_3)$$
which is a  special case of Lemma \ref{ccxf}. 

(ii) Let a based  matrix $T' $ be  obtained from a based matrix $T =(G ,s ,b )$ by a move $M_i$
with
$i=1,2,3$. We can assume that  the underlying set of $T'$  is a  union of a disjoint copy $ \{h'\,\vert \,h\in G\}$ of
$G$ and one new element $g$  in the case $i=1,2$ or two new elements $g_1,g_2$ in the case 
$i=3$.   For
$i=1$ (resp.
$i=2$, $3$),  the vectors
$\{h+h'\}_{ h\in G}$  and the vector $g$ (resp. $g-s'$, $g_1+g_2-s'$)  form a filling of  the pair $(T,
-T')$. The matrix of this filling  is zero. Hence $\sigma(T,-T')=0$ so that $T$ is cobordant to
$T'$. 

(iii) We need to prove that    $\sigma(T_1,...,T_r)$ is preserved  when  $T_1,...,T_r$ are replaced with
cobordant based matrices.  By induction, it suffices to prove that 
$\sigma(T_1,...,T_{r-1} , T'_r)=\sigma(T_1,...,T_{r-1} , T_r)$ for any based matrix   $T'_r$  cobordant to $T_r$.   
Lemma \ref{ccxf} gives that 
$$\sigma(T_1,...,T_{r-1} , T_r)\leq \sigma(T_1,...,T_{r-1} , T'_r)+\sigma(-T'_r ,
T_r) =\sigma(T_1,...,T_{r-1} , T'_r).$$  Similarly, 
$\sigma(T_1,...,T_{r-1} , T'_r)\leq 
\sigma(T_1,...,T_{r-1} , T_r)$. Hence 
 $\sigma(T_1,...,T_{r-1} , T'_r)=\sigma(T_1,...,T_{r-1} ,
T_r)$. 

(iv)  If  a  based matrix $T$ is cobordant to a  trivial based matrix $T_0=(\{s_0\}, s_0, b=0)$, then $\sigma
(T)=\sigma(T_0)=0$ and therefore  
$T$ is hyperbolic. Conversely, if $T=(G,s,b)$ is hyperbolic, then it has a filling with zero matrix. Adding to this
filling the vector $s+s_0$, we obtain a filling of the pair $(T, T_0)$ with zero matrix. Hence $T$ is cobordant
to $-T_0=T_0$.
\end{proof}

\begin{corol}\label{aaal:gfglpmns5} For any   based matrices  $T_1,...,T_r$ over $R$, we have
 $\sigma((T_1)_\bullet,...,(T_r)_\bullet)=\sigma(T_1,..., T_r)$.
		   \end{corol}

 \subsection{Exercises}\label{sdm1}   1. Verify that the definitions of the genus of a  (single) based matrix  over $R$ given 
in Sections \ref{fi:g5999} and  \ref{uussstrg52} are equivalent.

2. Prove that the function  $(T_1,T_2)\mapsto \sigma(T_1,-T_2)$ defines a metric on the set of cobordism classes of based
matrices over $R$. 

3. Prove that $\sigma(T^-_1,..., T^-_r)=\sigma(T_1,...,T_r)$  for any based matrices
$T_1,...,T_r$ over $R$.

4. Prove that  for any  $ 1\leq t\leq r$ and any based matrices  $  T_1,..., T_r$, $T'_1,...,T'_q$ with $q\geq 1$,  
$$\sigma(T_1,..., T_r)\leq  \sigma(T_1,..., T_t, T'_1,...,T'_q)
+ \sigma(-T'_1,...,-T'_q, T_{t+1},..., T_r)+q-1.$$

5. Prove that  $  u(T_1)+...+ u(T_r)=0$  for 	any hyperbolic  tuple  of based matrices  $T_1,...,T_r$ over $\RR$.

  \section{Genus estimates and sliceness of strings}

\subsection{Genus estimates for strings}\label{fi:g5999568}    
Setting $R=\ZZ$, we can apply the definitions and results of Section  \ref{fi:g599912}  to the
based matrices of strings.
We  begin with an estimate relating
the  slice genus   of   strings to the genus of their based matrices.

 \begin{lemma}\label{th:vv53}  For any string $\alpha$, we have $\sigma (T_\bullet (\alpha))=\sigma (T(\alpha))\leq 2\,
sg(\alpha)$.
\end{lemma}
\begin{proof}  The equality $\sigma (T_\bullet (\alpha))=\sigma (T(\alpha))$ follows from Corollary
\ref{aaal:gg25}. We prove that $\sigma (T(\alpha))\leq 2\,
sg(\alpha)$.

Consider an  oriented 3-manifold
$M$, a compact oriented surface  $F$ of genus
$ k=sg(\alpha)$   bounded by a circle $S$,  and a proper map
$\omega:F \to M$ such that
 $\omega\vert_{S}:S \to
\partial M$ is a (generic) closed curve on $\partial M$ realizing $\alpha$.  Let $\inc:H_1(\partial M) \to H_1(M)$ be the
inclusion homomorphism and $\omega_\ast:H_1(F)\to H_1(M)$ be the homomorphism  induced by $\omega$. Set $L=
 \inc^{-1} (\omega_\ast(H_1(F)))\subset  H_1(\partial M)$.  Since the intersection form $B:H_1(\partial M) \times
H_1(\partial M) \to
\ZZ$ annihilates the  kernel of
$\inc$ and $\rk\, \omega_\ast(H_1(F)) \leq \rk  H_1(F)=2k$, we obtain  that the rank of the bilinear form
$B\vert_L:L\times L \to \ZZ$ is smaller than or equal to $4k$.

Consider   the based matrix 
$T=T(\alpha) =(G,s,b)$ of $\alpha$.   As in the proof of Lemma \ref{l:nblddddd}, the map $\omega$ gives rise to an 
involution
$\nu$ on the set $\Join\! \! (\omega(S))=\arr(\alpha)= G-\{s\}$.  This defines a simple filling $ \mathcal X$ of   
$T $ consisting of  
$\{s\}$ and the orbits of $\nu$.  The proof of Lemma
\ref{l:nblddddd} shows that    $\sum_{x\in X} [\omega_x] \in L$ for any  orbit $X $  of $\nu$.
The homology class $[\omega(S)] \in H_1(\partial M)$ also lies in $L$ because $\inc ([\omega(S)])=0$. 
The matrix    of $ \mathcal X$    is obtained by evaluating $B$ on the   vectors $[\omega(S)]$ and $\{\sum_{x\in X}
[\omega_x]\}_X$ where $X$ runs over the  orbits of $\nu$. Therefore the rank of this matrix  is smaller than or equal to 
$\rk (B\vert_L)\leq 4k$. Hence 
$\sigma(T)\leq 2k=2\, sg(\alpha)$.
\end{proof}

The following  theorem  provides an algebraic obstruction to the sliceness of a string.

\begin{theor}\label{th:vwww53}  For a  slice string $\alpha$, the based matrices  $ T(\alpha) $ and $T_\bullet (\alpha)$ are
hyperbolic.
\end{theor}

This theorem is a direct consequence of  the previous lemma and the definitions.  We complement this theorem with the following result whose
proof is postponed to Section \ref{129ppp8}.

\begin{theor}\label{bnknk3}  Based
matrices  of cobordant strings are cobordant.
\end{theor}

\subsection{Genus estimates for sequences of strings}\label{fgb8}  We   generalize Lemma \ref{th:vv53} to sequences of
strings.  

  \begin{lemma}\label{th91373}  For any strings $\alpha_1,...,\alpha_r$,   $$ \sigma
(T_\bullet (\alpha_1),...,T_\bullet (\alpha_r ))=\sigma (T(\alpha_1),...,T(\alpha_r ))\leq
2\, sg(\alpha_1,...,\alpha_r).$$
\end{lemma}
\begin{proof}    The equality
$ \sigma
(T_\bullet (\alpha_1),...,T_\bullet (\alpha_r ))$ $=\sigma (T(\alpha_1),...,T(\alpha_r ))$ follows from Corollary \ref{aaal:gfglpmns5}.  
The rest of the proof is similar to the proof of Lemma \ref{th:vv53}. Consider
an    oriented 3-manifold
$M$, a compact (oriented) surface $F$   of genus
$  k=sg(\alpha_1,...,\alpha_r)$   bounded by $r$ circles $S_1,...,S_r$,   and a proper map
$\omega:F\to M$ such that
  the  maps $\omega\vert_{S_t}:S_t \to
\partial M $ with $t=1,...,r$ are disjoint (generic) closed curves on $\partial M$ realizing $\alpha_1,..., \alpha_r$,
respectively.   Let $\inc:H_1(\partial M) \to H_1(M)$ be the
inclusion homomorphism and $\omega_\ast:H_1(F)\to H_1(M)$ be the homomorphism  induced by $\omega$. 
The  group   $H_1(F)$  is generated by the homology classes of
$S_1,...,S_r\subset F$ and a subgroup   $H\subset H_1(F)$ isomorphic to $\ZZ^{2k}$. Set
$L=
\inc^{-1} (\omega_\ast(H))\subset  H_1(\partial M)$.  Since the intersection form $B:H_1(\partial M) \times
H_1(\partial M) \to
\ZZ$ annihilates the  kernel of
$\inc$, we obtain  that  
$$\rk\,  (B\vert_L:L\times L \to \ZZ) \leq 2 \rk\, \omega_\ast(H) \leq 2\rk  H=4k.$$

For $t=1,...,r$, consider   the based matrix $T_t=(G_t,s_t,b_t)$ of $\alpha_t$.  Set $G=\cup_t G_t$.    As in the proof
of Lemma
\ref{l:nblddddd}, the map
$\omega$ gives rise to an  involution
$\nu$ on the set $\Join\! \! (\omega(\partial F)) = G-\{s_1,...,s_r\}$.    The proof of Lemma
\ref{l:nblddddd} shows that    for any  orbit $X $  of $\nu$, we have  $\sum_{x\in X} [\omega_x] \in   \inc^{-1}
(\omega_\ast(H_1(F))) $.  Adding to $\sum_{x\in X} [\omega_x]$ an appropriate    linear
combination
$\sum_t n_{X,t}  [\omega(S_t)]$ of the homology classes $[\omega(S_1)], ..., [\omega(S_r)]
\in H_1(\partial M)$  with $n_{X,t}\in \ZZ$ we obtain an element of $ L$.  
Consider the vector  $ \sum_{x\in X} x +\sum_t n_{X,t}  s_t  $ in the lattice $\ZZ G $ freely generated by  
$G$.  These vectors corresponding to all  orbits $X$  of $\nu$  
 together with  the vector $s_1+...+s_r\in \ZZ G $ form a 
filling   of   the tuple 
$T_1,...,T_r $. The matrix    of this filling     is obtained by evaluating $B$ on the  homology classes
 $\{\sum_{x\in X} [\omega_x] +\sum_t n_{X,t}  [\omega(S_t)] \in H_1(\partial M)\}_X$  and  $[\omega(S_1)]+
...+ [\omega(S_r)]\in H_1(\partial M)$.  Since all these homology classes belong to
$L$, 
  the rank of this matrix 
is smaller than or equal to $4k$. Thus 
$\sigma (T(\alpha_1),...,T(\alpha_r ))\leq 2k=
2\, sg(\alpha_1,...,\alpha_r)$.
\end{proof}

\begin{theor}\label{tqpaw53}  If a sequence of strings    is slice, then  the sequence of their based
matrices  and the sequence of their primitive based
matrices are hyperbolic.
\end{theor}

This theorem is a direct consequence of  the previous lemma and the definitions. 

\subsection{Proof of Theorem \ref{bnknk3}}\label{129ppp8} If strings $\alpha$, $\beta$ are cobordant, then $sg(\alpha, {\overline
\beta}^-)=0$. By Lemma
\ref{th91373}, $\sigma (T(\alpha), T({\overline \beta}^-))=0$.  As we know,  $T({\overline \beta}^-)=-T(\beta)$. Thus,
 $\sigma (T(\alpha),  -T(\beta))=0$ so that $T(\alpha)$ is cobordant to ~$T(\beta)$.

\subsection{Secondary obstructions to   sliceness}\label{1hhhgp8} We introduce    invariants of   strings which   may  give further
obstructions to   sliceness (cf. Question 2 in Section \ref{questions}). Consider a string $\alpha$ with core circle  $S$ and canonical realization 
$\omega:S\to \Sigma_\alpha  $  as in Section
\ref{fi:g31}. Let
$\Sigma$ be the closed oriented surface obtained by gluing 2-disks to all components of $\partial \Sigma_\alpha$.  Pick an integer $p\geq 2$
and set  $R=\ZZ/p\ZZ$ and 
$H= H_1(\Sigma; R)$.  The $R$-module 
$H$  is generated by the set  $s\cup \{[e]\}_{e\in
\arr(\alpha)}$ where   the homology classes of loops on $\Sigma$ are taken with coefficients in $R$ and $s=[\omega(S)]\in H$, cf.   Section
\ref{fi:g32}.  Consider the  intersection form
$B_R:H \otimes H  \to R$.   For
$h\in H$, consider the  string
$\alpha_h$ formed  by  
$S$ and the arrows
$e\in
\arr (\alpha)$ such that $B_R([e],h)=0$.  The invariants of  $\alpha_h$ can be viewed as  invariants of $\alpha$ parametrized by $p$ and $h$. In
particular, we can consider the 1-variable polynomial $ u(\alpha_h)$, the based matrix  $ T(\alpha_h)$, etc.

In the next lemma, a {\it Lagrangian}  
is a  group $L\subset H$   equal to its annihilator  $\Ann (L)=\{g\in H\,\vert \, B_R(L,g)=0\}$. If $p$ is   prime, then    each
Lagrangian $L\subset H$ is a direct summand of $H$ and $H/L\approx L$.

 \begin{lemma}\label{th955555553}  If $\alpha$ is slice, then 
there is a Lagrangian  $L\subset H$ such that $s\in L$ and the string $\alpha_h$ is
slice for all $h\in L$. Moreover,  there is an involution
  on the set $ \arr(\alpha) $ such that  for any  its orbit $X$,    $\sum_{e\in X} [e] \in L$.
\end{lemma}
\begin{proof}   If $\alpha$ is slice, then there are a compact oriented 3-manifold $M'$ and  a realization 
$\omega' :S\to \partial M'$     of  $\alpha$    contractible in $M'$. 
The pair  
$(\partial M', \omega')$ can be obtained from   $ \omega:S\to \Sigma$ of $\alpha$    by   1-surgeries  on  $\Sigma-\omega(S)$. Attaching  the
corresponding solid 1-handles to
$\Sigma \times 0 \subset \Sigma\times [0,1]$ we obtain an oriented 3-manifold  $N$ 
such that  $\partial N= (-\partial M') \cup \Sigma$  and the curves $  \omega', \omega$ are
homotopic in $N$.    Gluing $N$ to $M'$ along $\partial M'$ we obtain a compact oriented 3-manifold $M$ such that $\partial M=\Sigma$ and
$\omega$ is contractible in $M$. 
Consider  the boundary  homomorphism $\partial : H_2(M,\partial M;R) \to H_1(\partial M;R)  =H$ and  the inclusion
homomorphism
$i:H=H_1(\partial M;R) \to H_1(M; R)$. Set $L=\Ima (\partial)=\Ker (i)$. 
  It is well known that
  $L  $ is a Lagrangian.  For completeness, we outline a proof. An element $g\in H$ belongs to $ \Ann (L)$ iff  $B_R (  \partial x,g )=0$ for every 
$x\in H_2(M,\partial M;R)$. By the Poincar\'e duality,  there is a unique $\tilde x \in H^1(M;R)=\Hom (H_1(M;R), R)$ such that $x=\tilde x \cap
[M]$.  Then  
$B_R ( \partial x, g)=x\cdot i(g)=
\tilde x (i(g))$  where $\cdot$ is the   intersection pairing  
$H_2(M,\partial M;R)\times H_1(M ;R)\to R$.   Therefore $g\in \Ann (L)$ iff  $i(g)$ is annihilated by all homomorphisms $H_1(M;R) \to R$. 
This holds iff
$i(g)=0$, that is iff
$g\in L$. 

Pick  
$h\in L$ and  pick   any $x $ in  $\partial^{-1} (h)\subset  H_2(M,\partial M;R)$. The   cohomology class $\tilde
x\in H^1(M;R)$ defines a
$p$-fold covering  
$\tilde M \to M$. Since $\omega:S\to \Sigma=\partial M$ is contractible in $M$, it lifts   to a loop $\tilde \omega: S\to \partial \tilde M$ 
contractible in $\tilde M$. By the   equality
$\tilde x (i(g))= B_R(h,g)$ for $ g\in H$,  the underlying string of 
$\tilde \omega$ is $\alpha_h$.  Therefore  $\alpha_h$ is slice. Constructing an involution  
 on  $ \arr(\alpha) $ as in the proof of Lemma \ref{l:nblddddd}  (where $F$ is a 2-disk) we obtain the last claim of the lemma.
\end{proof}

 This lemma implies that for  all $h\in L$,    the based matrix $T(\alpha_h)$ is hyperbolic and  $u (\alpha_h)=0$.

   \section{Lie cobracket for strings}\label{cob}

We introduce a Lie cobracket in  the  free module generated by homotopy classes of strings.
This induces a Lie bracket in the module of homotopy invariants of strings and other related algebraic structures.  

Throughout the section, we fix a 
commutative ring 
 with unit $R$.

		         \subsection{Lie   coalgebras}\label{su:61} 
We recall  here the  notion of a Lie coalgebra  dual to the one of a Lie algebra. To 
this end, we first  reformulate the notion of a Lie algebra.   For an $R$-module $L$,   denote by $\Perm_L$ the permutation 
$x\otimes 
y\mapsto y\otimes x$ in $L^{\otimes 2}=L\otimes  L$ and by $\tau_L$ the 
permutation $x\otimes y\otimes z\mapsto z\otimes x\otimes y$ in 
$L^{\otimes 
3}=L\otimes L \otimes   L$. Here and below  
$\otimes=\otimes_R$.  
 A Lie algebra over $R$ is an $R$-module $L$ endowed with an $R$-homomorphism (the 
Lie bracket)  $\theta:L^{\otimes 2}\to L$ such that   $\theta \circ 
\Perm_L=-\theta$ (antisymmetry) and    $$\theta \circ (\id_L\otimes 
\theta) 
\circ (\id_{L^{\otimes 3}}+\tau_L+\tau_L^2)=0 \in \Hom_R(L^{\otimes 3},L) $$
(the Jacobi identity).  Dually, a Lie coalgebra over $R$ is an $R$-module 
$A$ 
endowed with an $R$-homomorphism (the Lie cobracket)  $\nu:A\to 
A^{\otimes 2}  $ 
such that   $  \Perm_A \circ \nu=-\nu$   and  
			 \begin{equation}\label{jac}(\id_{A^{\otimes 
3}}+\tau_A+\tau_A^2)  \circ (\id_A\otimes \nu) \circ \nu=0 \in \Hom_R(A,  A^{\otimes 3}).
\end{equation}

For a Lie coalgebra $(A,\nu:A\to A^{\otimes 2})  $ over     $R$ and
an integer 
$n\geq 
1$, set  
   $$\nu^{(n)}=(\id_A^{\otimes (n-1)} \otimes \nu) \circ \cdots \circ 
  (\id_A  \otimes \nu) \circ \nu:A \to A^{\otimes (n+1)}.$$
In particular, $\nu^{(1)}=\nu$.
Following \cite{tu}, Section 11, we call a  Lie coalgebra $(A,\nu)$ over $R$  {\it spiral}, 
if $A$ is 
free as an  $R$-module and the filtration $\Ker \nu^{(1)}\subset \Ker 
\nu^{(2)}\subset 
\cdots $ exhausts $A$, i.e., $A=\cup_{n\geq 1} \Ker \nu^{(n)}$. 

A Lie 
coalgebra $(A,\nu)$ gives rise to the {\it dual Lie algebra}
  $A^*=\Hom_R(A,R)$ where the Lie bracket  $ A^*\otimes A^* \to 
A^*$ is  the homomorphism 
dual to $\nu$.  For $u,v\in A^*$, the value of $[u,v]\in A^*$ on     $x\in A$ is computed by
$$[u,v] (x)=\sum_i u(x^{(1)}_i) \, v (x^{(2)}_i)\in R$$
for any (finite) expansion $\nu(x)=\sum_i x^{(1)}_i \otimes x^{(2)}_i\in A\otimes A$.

A {\it homomorphism} of Lie coalgebras $(A,\nu)\to (A',\nu')$ is an $R$-linear 
homomorphism $\psi:A\to A'$ such that
$(\psi \otimes \psi) \nu (a)= \nu' \psi(a)$ for all $a\in A$. 
It is clear that the dual homomorphism $\psi^*:(A')^*\to A^*$ is a 
homomorphism 
of Lie algebras.

  \subsection{Lie   coalgebra of strings}\label{su71}
Let $\str$ be the 
set of 
homotopy classes of virtual strings and let $\str_0 \subset \str$ be its subset formed by the  homotopically non-trivial
classes.  Let
$\A_0=\A_0(R)$ be the free 
$R$-module freely generated by $\str_0$. We shall 
provide $\A_0$ with the 
structure 
  of a Lie coalgebra. 
  
 We begin with notation. For a  string $\alpha$,
let $\langle 
\alpha 
\rangle $ denote its class in $\str_0$ if $\alpha$ is homotopically non-trivial and set $\langle 
\alpha 
\rangle=0\in \A_0$ if $\alpha$ is homotopically trivial.   
 For an arrow $e=(a,b)$ of a string $\alpha$, denote
by $\alpha^1_e$ 
the string 
obtained from $\alpha$ by removing all arrows except
those with both 
endpoints 
in the interior of  the arc $ab$. (In particular, $e$ is removed.) Similarly, denote by $\alpha^2_e$ the string
obtained from 
$\alpha$ by  
removing all arrows except those with both endpoints
in the interior of  
$ba$. 
 Set
\begin{equation}\label{cobr}\nu (\langle \alpha
\rangle)= \sum_{e\in 
\arr(\alpha)}  \langle \alpha^1_e\rangle\otimes 
 \langle \alpha^2_e\rangle-\langle
\alpha^2_e\rangle\otimes 
 \langle \alpha^1_e\rangle  \in \A_0\otimes \A_0.\end{equation}

  \begin{lemma}\label{l:72} The $R$-linear homomorphism
$\nu: \A_0\to \A_0\otimes 
\A_0$ given 
on the generators of $\A_0$ by Formula \ref{cobr} is a
well-defined Lie 
cobracket. The Lie coalgebra  $(\A_0 ,\nu)$     is spiral.
		   \end{lemma}
                     \begin{proof} To show that $\nu$
is well-defined 
we must 
verify that   $\nu (\langle \alpha
\rangle)$ does not 
change under 
the homotopy moves (a$)_s$, (b$)_s$, (c$)_s$ on
$\alpha$. The  arrow 
added by 
(a$)_s$ contributes $0$ to the cobracket by the definition of    $\langle...\rangle$. The 
contribution of all the other arrows is preserved.
Similarly, the two 
arrows 
added by (b$)_s$ contribute  opposite terms to the cobracket  which is therefore preserved. 
Under (c$)_s$, all 
arrows 
contribute  the same before and after the move.	 
	
 The equality $\Perm_{\A_0} \circ \nu=-\nu$ is
obvious. We now 
verify Formula \ref{jac}. Let $\alpha$ be a string
with core circle 
$S$. We can expand   $(\id \otimes \nu)
( \nu (\langle 
\alpha 
\rangle))$ as a sum of   expressions $z(e,f)$
associated with     
ordered 
pairs of unlinked arrows $e,f\in \arr(\alpha)$.    
Note that the 
endpoints of 
$e,f$ split $S$ into four arcs meeting only at their
endpoints. The endpoints of $e$ (resp. $f$) bound one of these arcs, say $x$ (resp. $y$).   
 The 
other two arcs form $ S- (x\cup y)$ and lie \lq\lq
between" $e$ and 
$f$. Denote 
by $  \beta $ 	     (resp. $ \gamma $, $ \delta $) the
string obtained 
from 
$\alpha$ by removing all arrows except those with both
endpoints in the 
interior 
of $x$ (resp. of $y$, of   $ S- (x\cup y)$).  Set
$\varepsilon=+1$ if 
$e$ and 
$f$ are co-oriented, i.e., if their tails bound a
component of  $ S- 
(x\cup y)$. 
It is easy to see that
		     $$z(e,f)=  \varepsilon  (\langle
\beta\rangle \otimes 
\langle \delta\rangle \otimes  \langle \gamma\rangle 
- \langle 
\beta\rangle \otimes  \langle \gamma\rangle \otimes
 \langle 
\delta\rangle ).$$
A direct computation using this formula gives
$$(\id_{{\A_0}^{\otimes 3}}+\tau_{\A_0}+\tau_{\A_0}^2)
(z(e,f)+z(f,e))=0.$$		 
Thus $   \id_{{\A_0}^{\otimes 3}}+\tau_{\A_0}+\tau_{\A_0}^2
$ annihilates  
  $(\id \otimes \nu) ( \nu (\langle \alpha \rangle))$.
Hence $\nu$ is a 
Lie 
cobracket. The spirality of $(\A_0,\nu) $  follows from
the obvious fact 
that
 $\nu^{(n)}(\langle \alpha\rangle)=0$ for any string
$\alpha$ of rank $\leq 
n$. 
(Actually a stronger assertion holds: $\nu^{(n)}(\langle
\alpha\rangle)=0$ 
for any 
string $\alpha$ of rank $\leq 4n+2$.) \end{proof}

        Let
$\A=\A (R)$ be the free 
$R$-module freely generated by $\str$.  Since $\str=\str_0 \cup \{O\}$  where $O\in \str$ is the homotopy class of a trivial
string,  $\A=\A_0 \oplus R O$.  The Lie cobracket $\nu$ in  $\A_0$ extends to $\A$ by $\nu (O)=0$.

The   Lie cobrackets  in $\A_0$  and $\A$ induce 
  Lie brackets 
in  $\A_0^* 
=\Hom_R (\A_0 , R)$ and  $\A^* 
=\Hom_R (\A , R)$.  
Examples below show  that these  Lie cobrackets  and   Lie brackets
 are non-zero.  Clearly, $\A^*=  \A_0^*  \oplus R$ where the Lie bracket in $R$ is zero. 
The elements of
$\A^*
$ bijectively correspond to maps $\str \to R$, i.e., to 
$R$-valued 
homotopy invariants of strings. Thus, such invariants
form a Lie 
algebra.

	\subsection{Examples.}\label{su99}   (1) If   $\rank \alpha \leq 6$, then   $\nu (\langle
\alpha
\rangle)=0$. This follows from the fact that any string     of rank $\leq 2$ is homotopically trivial.  

(2) For any $p,q\geq 1$, we have $\nu (\langle\alpha_{p,q}\rangle)=0$.

(3)  Consider the string 
$ \alpha_\sigma$ of rank 7 where $\sigma$ is the permutation $ (123) (4) (576)$ of the set $\{1,2,\ldots, 7\}$. 
It follows from the definitions that
 $\nu(\langle \alpha_\sigma \rangle)= \langle \alpha_{1,2} \rangle\otimes \langle \alpha_{2,1} \rangle- \langle \alpha_{2,1}
\rangle
\otimes \langle \alpha_{1,2} \rangle$.
 As we know,  $\alpha_{1,2}$ and $ \alpha_{2,1}$   are homotopically non-trivial  strings      representing  distinct
generators of $\A$. Hence $\nu(\langle \alpha_\sigma\rangle)\neq 0$.   This example can be used to show    that the
product of strings is not commutative even up to homotopy: 
 there are  strings $\gamma,
\delta$ such that a product of $\gamma, \delta$ is not homotopic to a product of $\delta, \gamma$.    Drawing a picture
of  $\alpha_\sigma$, one  observes that $\alpha_\sigma$ is a product of $\delta=\alpha_{2,1}$ with a string,
$\gamma$, of rank 4 obtained  from $\alpha_{1,2}$ by  adding a \lq\lq small" arrow.  
Since $\gamma$ has a small arrow, it is easy to form a product of $\gamma$ with $\delta$ also having a
small arrow.  The resulting string, $\beta$,  is homotopic to a string of rank 6. Hence $\nu (\langle \beta\rangle )=0$.
Therefore $\alpha_\sigma$ is not homotopic to $\beta$.

(4) In generalization of the previous example pick any integers $p,q,p',q'\geq 1$ such that $p+ q\geq 3, p'+q'\geq 3$.
 Consider the string 
$\alpha=\alpha_\sigma$ of rank $m=p+q+p'+q'+1$  where $\sigma$ is the permutation   of the set
$\{1,2,\ldots, m\}$ defined by
$$\sigma(i)=\left\{\begin{array}{ll}
i+q ,~ {\rm {if}} 
\,\,\, 
1\leq i\leq p \\
\noalign{\smallskip}
i-p
,~ 
{\rm
{if}} \,\,\, p<i\leq p+q
\\
\noalign{\smallskip}
i,~ 
{\rm
{if}} \,\,\, i=p+q+1
\\
\noalign{\smallskip}
i+q'
,~ 
{\rm
{if}} \,\,\, p+q+1<i\leq p+q+1+p'  
\\
\noalign{\smallskip}
i-p'
,~ 
{\rm
{if}} \,\,\, p+q+1+p'  <i\leq m.
\end{array} \right.$$
  It follows from the definitions that
$$\nu(\langle \alpha\rangle)=  \langle \alpha_{p',q'} \rangle \otimes  \langle \alpha_{p,q} \rangle -  \langle \alpha_{p,q}
\rangle \otimes  \langle \alpha_{p',q'} \rangle .$$
Clearly,  $\nu(\langle \alpha\rangle)\neq 0$ unless $p= p'$ and  $q= q'$.

(5) Consider the numerical
invariants $u_1, u_2,\ldots \in  \A^* $ constructed in Section \ref{sn:g21}.  For  $p,p'\geq 1$, we  compute the
value of
$[u_p,u_{p'}]\in 
\A^*
$ on the  string $\alpha=\alpha(p,p',q,q')  $ defined in the previous example. Assume for  concreteness that   the numbers 
$p,p',q,q' $ are pairwise distinct. Then
$$[u_p,u_{p'}] (  \alpha )=u_p(\alpha_{p',q'} ) \,u_{p'} (\alpha_{p,q})- u_p(\alpha_{p,q} )\,  u_{p'} (\alpha_{p',q'}) =0-
(-q) (-q')=-qq'.$$
Hence $[u_p,u_{p'}]\neq 0$ for $p\neq p'$.

  \subsection{Filtration of $\A_0$.}\label{su735} Assigning to a string its homotopy rang and homotopy genus (see  
Section \ref{sn:g133}) we obtain two maps $hr, hg:\str_0 \to \ZZ$.  For   $r, g\geq 0$, set $$\str_{r,g}=\{\alpha \in
\str_0 \,
\vert \, hg(\alpha)\leq r,\,\,\,\, hg(\alpha)\leq g\}.$$
This set is finite since there is only a finite number of strings of rank $\leq r$.  The set  $\str_{r,g}$ generates a 
submodule  of $\A_0$ denoted $\A_{r,g}$.  This submodule is a free $R$-module of   rank $\#(\str_{r,g})$. 
Clearly, 
\begin{equation}\label{fil}\nu ( \A_{r,g})\subset \bigoplus_{p,q\geq 0, p+q<r} \A_{p,g} \otimes \A_{q,g} \subset  
\A_{r,g} \otimes \A_{r,g}.\end{equation}
Thus, each $\A_{r,g}$ a  Lie coalgebra.  The   inclusions $\A_{r,g}\hookrightarrow \A_{r',g'}$ for   $r\leq r', g\leq
g'$ make the family  $\{\A_{r,g}\}_{r,g}$ into a direct spectrum of Lie coalgebras. The equality $\A_0=\cup_{r,g}
\A_{r,g}$ shows that $\A_0=\injlim \{\A_{r,g}\}$.

The Lie cobracket in $\A_{r,g}$ induces
 a Lie bracket 
in   $\A_{r,g}^* 
=\Hom_R (\A_{r,g} , R)$.  Formula \ref{fil} implies that this Lie algebra is nilpotent. Restricting  maps  $\str_0 \to R$ to
$\str_{r,g}$ we obtain a Lie algebra homomorphism 
	$\A_0^*\to  \A_{r,g}^* $.  It is clear that $\A_0^*=\projlim \{\A^*_{r,g}\}$.

	\subsection{Relations with Lie coalgebras of
curves.}\label{su75}
	Let $\Sigma$ be a connected  surface and 
	$\hat \pi$ be the set of homotopy classes of closed
curves on $\Sigma$. 
(It can be identified with the set of conjugacy
classes in 
	$\pi=\pi_1(\Sigma)$.) 	 There is a  map $\psi : \hat \pi \to \str$ sending
each homotopy 
class 
of curves into the homotopy class of the underlying
strings.  Clearly, 
 $\psi  (\hat \pi )= \cup_r\,  \str_{r,g} $
where $g= g(\Sigma) $ is the genus of $\Sigma$.
Observe that the mapping class group of $\Sigma$ acts on $\hat \pi$ in the obvious way and  $\psi$ factors through the
projection of $\hat \pi$ to the set of orbits of this action.

Let $Z =Z(R)$ be the free
$R$-module with 
basis $\hat \pi$.  The map $\psi: \hat \pi \to \str$ induces an $R$-linear homomorphism $ Z\to \A$ whose image is equal to 
 $  \cup_r\,  \A_{r,g} $.  Composing this homomorphism with the   projection  $\A=\A_0\oplus R O \to \A_0$  we obtain 
an $R$-linear homomorphism  $\psi_0:Z\to \A_0$. 

The author defined in \cite{tu},
Section 8 a 
structure of a spiral
Lie coalgebra in $Z$. (In fact $Z$ is a Lie bialgebra,
but we do not 
need it.)     A direct comparison of the
definitions shows that the  map  $\psi_0 :Z\to \A_0$  
is a homomorphism 
of Lie 
coalgebras.

 \subsection{Associated algebraic structures.}\label{su:6111}   In this section we suppose that $R\supset \QQ$.  A spiral Lie
coalgebra
$(A,\nu)$ over $R$ naturally gives rise to a group  $\Exp A^*$   and  a Hopf algebra   $S(A)$ over $R$, see \cite{tu},
Section 11. For completeness, we recall here these constructions.

Observe first that the  dual Lie algebra $A^*=\Hom_R(A,R)$   has the following completeness property. 
Consider 
the lower central series $A^*=A^{*(1)}\supset A^{*(2)}\supset \cdots$ of 
$A^*$ 
where $A^{*(n+1)}= [A^{*(n)}, A^*]$ for $n\geq 1$. Let $a_1,a_2,\ldots 
\in A^*$ 
be an infinite sequence such that for any $n\geq 1$ all terms of the 
sequence 
starting from a certain place belong to $A^{*(n)}$. Clearly, if $x \in 
\Ker 
\nu^{(n)}$ and $a\in A^{*(n+1)}$, then $a(x)=0$.  Since $A=\cup_{n } 
\Ker 
\nu^{(n)}$, 
the sum $a(x)=a_1(x)+a_2(x)+\cdots$ contains only a finite number of 
non-zero 
terms for every $x\in A$.  Therefore $a(x)$ is a well-defined element of $R$. The formula $x\mapsto a(x):A\to R$   defines  
an element  of
$   A^* $  denoted $a_1+ a_2+ \cdots$ and called the (infinite) 
sum of 
$a_1,a_2,\ldots$.  A similar argument shows that 
 $\cap_n  A^{*(n)}=0$ and 
the natural Lie algebra homomorphism   $A^*\to \projlim_n (A^*/A^{*(n)})$ 
is an 
isomorphism.

  For   $a,b\in A^*$,    
consider the  sum
$$\mu (a,b)= a+b + \frac{1}{2} [a,b]+ \frac{1}{12} ([a,[a,b]] +[b, 
[b,a]])+ 
\cdots \in A^*$$
where the right-hand side is the Campbell-Hausdorff series for $\log 
(e^a e^b)$, 
see \cite{se}. The resulting mapping $\mu:A^*\times A^*\to A^*$ is a 
group 
multiplication in $A^*$. Here $a^{-1}=-a$ and $0$ is the group unit. 
The group  
$(A^*,\mu) $ is denoted $\Exp A^*$. Heuristically, this is 
the \lq\lq Lie group" with Lie algebra $A^*$. The equality  $A^*= \projlim_n (A^*/A^{*(n)})$ 
implies that   the  group $\Exp A^*$ 
is 
pro-nilpotent.

Consider the symmetric (commutative and 
associative) algebra of $A$:
$$S=S(A)=\oplus_{n\geq 0} S^n (A).$$   Here $S^0(A)=R$, $S^1(A)=A$, and $S^n(A)$  
is the 
$n$-th symmetric tensor power of $A$ for $n\geq 2$. The unit $1\in R= 
S^0(A)$ is 
the unit of $S$. The group multiplication $\mu:A^*\times A^*\to A^*$ 
  induces a comultiplication   $S\to S\otimes S$  as follows. Since 
$A$ is a 
free $R$-module, the natural map $A\to (A^*)^*$ extends to an embedding 
of $S$ 
into the algebra of $R$-valued functions on $A^*$. We can identify
 $S$ with the image of this embedding. Similarly, we can identify 
$S\otimes S$ 
with  an algebra of $R$-valued functions on $A^*\times A^*$.  It is 
easy to 
observe that for any $x\in S$, we have $x\circ \mu\in S\otimes S$. 
Indeed, it 
suffices to prove this for   $x\in A$. Then  
$x\in \Ker 
\nu^{(n)}$ for some $n$ so that $x$ annihilates all but finite number of terms of the  
Campbell-Hausdorff series. Our claim follows then from the duality 
between the 
Lie bracket in $A^*$ and the Lie cobracket $\nu$. For example, if $n=3$ 
and 
$\nu^{(2)}(x)=\sum_i \alpha_i\otimes \beta_i \otimes \gamma_i \in 
A^{\otimes 3}$, 
then
 $$x\circ \mu=x\otimes 1+1\otimes x+ \frac{1}{2} \nu(x)+\frac{1}{12}
 \sum_i (\alpha_i  \beta_i \otimes \gamma_i+ \gamma_i\otimes \alpha_i  
\beta_i )
.$$
 The formula $\Delta (x)=x\circ \mu$ defines a coassociative 
comultiplication in 
$S$. It has a  counit $S\to R$ defined as  the projection to $S^0(A)=R$. The antipode $S\to S$ is the 
algebra 
homomorphism sending any $x\in A$ into $-x\in A$.  A routine check 
shows that 
  $S$   is a  (commutative)
Hopf 
algebra.
 Heuristically, it should be viewed as the Hopf algebra of  $R$-valued functions on
the 
group $\Exp A^*$ or as the Hopf   dual of the universal enveloping   
algebra of 
  $A^*$. 
 
The construction of     $\Exp A^*$ and    $S(A)$ can be  generalized as follows.   Pick   $h\in R$ and observe that
 the  mapping $ h \nu :A\to A\otimes A$ is  a Lie cobracket in
$A$.  It induces the Lie bracket $[,]_h=h [,]$ in $A^*$ where $[,]$ is the Lie bracket induced by $\nu$.  The corresponding
multiplication
$\mu_h$  in
$A^*$ is given by
$$\mu_h (a,b)= a+b + \frac{h}{2} [a,b]+ \frac{h^2}{12} ([a,[a,b]] +[b, 
[b,a]])+ 
\cdots$$
This multiplication   makes $A^*$ into a group denoted  $\Exp_h A^*$. As above, $\mu_h$  induces a comultiplication in
the symmetric algebra $S=S(A)$. This makes $S$  into a Hopf algebra  over $R$ denoted $ S_h(A)$.  For
$h=1$,  
we obtain the same objects as in the previous paragraphs.  Note for the record that for any $h\in R$, the formula  
$a\mapsto ha: A^*\to A^*$ defines a group homomorphism $\Exp_h A^*\to \Exp  A^*$.  If $h\in R$ is a non-zero-divisor,
this homomorphism is injective.

 We can apply  the constructions  of this subsection to any $h\in R$ and to the spiral Lie coalgebras       $\A_0, \A,
\A_{r,g}  $, $ Z $ considered above.   The equality $\A=\A_0 \oplus R $ implies that 
 $\Exp_h  \A^*=\Exp_h   \A_0^* \times \underline R $ where  $\underline R$ is the additive
group of 
$R$.    The  group $\Exp_h\A_{r,g}^* $ and the Hopf algebra $ S_h ( \A_{r,g} )$  
are quotients of  
$\Exp_h
\A_0^* $ and   $S_h (\A_0) $, respectively. 
 The homomorphism     $\psi_0:Z\to \A_0$ extends by multiplicativity to  a  Hopf
algebra homomorphism   $S_h(Z)\to 
S_h (\A_0)$.  Dualizing $\psi_0$,  we obtain     a mapping
 $\A_0^*\to Z^*$ which is  a Lie algebra
 homomorphism  and at the same time   
  a group 
homomorphism 
$\Exp_h\A_0^*\to \Exp_hZ^*$.

\section{Virtual strings versus  virtual  knots}\label{vks}

 Virtual knots were introduced by L. Kauffman \cite{Ka} as a  
generalization of  classical knots.  We relate them to  virtual   strings by showing that each virtual knot gives rise to a
polynomial on virtual strings with coefficients in the   ring $\QQ [z]$.  As a technical tool, we
introduce a skein algebra of virtual knots and compute it in terms of  strings.

    \subsection{Virtual knots.}\label{vks:1}    We   define virtual knots  in terms of arrow diagrams following
\cite{gpv}.  An {\it arrow diagram} is a virtual string whose arrows are endowed
with   signs $\pm$.  By the core circle and the endpoints of  an arrow diagram, we mean the core circle and the endpoints of 
the   underlying   virtual string. 
The sign of an  arrow $e$ of an arrow diagram  is denoted
$\sign (e)$.   Homeomorphisms of arrow diagrams are defined as the homeomorphisms of the underlying strings preserving
the signs of all arrows.  The homeomorphism classes of       arrow diagrams will be  also called  arrow diagrams.  

We describe  
three moves (a$)_{ad}$, (b$)_{ad}$, (c$)_{ad}$ on arrow diagrams where $ad$ stands for
\lq\lq arrow diagram". Let $\alpha$ be an arrow diagram with core circle $S$. Pick two  distinct points 
$a,b \in S$  such that the (positively oriented) arc   $ab\subset S$ is disjoint from the set of   
endpoints of 
$\alpha$. The move (a$)_{ad}$   
adds to 
$\alpha$  the arrow $(a,b)$ with sign $+$ or $-$.   This move has two forms determined by the sign $\pm$.   The move
(b$)_{ad}$  acts on 
$\alpha$ as follows. Pick two   arcs   on $S$ disjoint from each 
other 
and   from the   endpoints of $\alpha$. Let $a,a'$ be the 
endpoints of 
the first arc (in an arbitrary order) and $b,b'$ be the endpoints of 
the second 
arc. The move adds 
to $\alpha$ two arrows $(a,b)$ and $(b',a')$ with opposite signs.  This move has eight forms depending on the
choice of the sign of $(a,b)$,  two possible choices for  $a$, and two possible choices for  $b$.   (This list of
eight forms of   (b$)_{ad}$ contains two equivalent pairs so that in fact the move  (b$)_{ad}$ has only six forms.)
 The
move (c$)_{ad}$     applies  to 
$\alpha$   when $\alpha$ has three arrows  with signs 
$((a^+,b),+), ( (b^+,c),+), ( (c^+,a),-)$ where  $a, a^+, b,b^+,c, c^+\in S$ such that
the   arcs   $aa^+$,   $bb^+$,     $cc^+$ are disjoint from each 
other and 
from the other   endpoints of $\alpha$. The   move (c$)_{ad}$  replaces 
these three  arrows    with the arrows  $ ((a,b^+),+), ( (b,c^+),+), ((c,a^+),-)$.  

By definition, a {\it virtual  knot}  is an equivalence class of arrow diagrams with respect to the equivalence relation
generated by  the  moves (a$)_{ad}$, (b$)_{ad}$,
(c$)_{ad}$ and  homeomorphisms. Note that our set of moves is somewhat different from the one in 
\cite{gpv} but  generates the same equivalence relation (cf. below). 

In the sequel the virtual knot represented by an arrow diagram $D$
will be denoted $[D]$.   A {\it trivial arrow diagram} having no arrows represents the {\it trivial virtual knot}. 

Forgetting the signs of arrows, we  can associate with  any arrow diagram $D$ its underlying virtual string $\underline
D$.   This induces   a \lq\lq forgetting" map  $K\mapsto \underline K$ from   the set of virtual knots 
   into the set of virtual strings. This map is surjective bur not  injective.  The theory of virtual knots is considerably reacher
than the theory of virtual strings. For instance,   the   fundamental  group of a virtual knot  \cite{Ka}  allows to
distinguish virtual knots with the same underlying strings.

 Note finally that the definition of  an $r$-th covering of a string
  in Section \ref{covvvr1} extends to virtual knots: one keeps only  arrows $e$ of an arrow diagram such that $n(e)\in r\ZZ$ and of course one
keeps their signs. 

   \subsection{From knots to virtual knots.}\label{vks:2}  Arrow diagrams are closely related  to the   standard
  knot diagrams on surfaces. An (oriented)
  knot diagram on an (oriented) surface $\Sigma$ is a (generic oriented)  closed curve on $\Sigma$ such that at
each its double point one of the branches of the curve passing through this point is distinguished. The distinguished branch is
called an {\it overcrossing} while the second branch passing through the same point is  called an {\it undercrossing}.  A knot diagram on
$\Sigma=\Sigma\times
\{0\}$ determines an (oriented)  knot in the cylinder 
$\Sigma\times
\RR$ by   pushing the overcrossings into  $\Sigma\times (0,\infty) $.  

Any knot diagram $d$ gives rise to an arrow diagram  $D(d)$  as follows.  First of all,  the closed curve underlying $d$ gives
rise to a virtual string, see Section \ref{sn:g12}. We provide each arrow  of this string with the sign of the corresponding
double point of $d$. This sign is $+$ (resp. $-$)  if the  pair  (a positive tangent vector to the
overcrossing 
 branch, a positive tangent vector to the undercrossing 
 branch) is positive  (resp. negative)  with respect to the orientation of $\Sigma$.  Our definition of the
arrow diagram associated with $d$ differs   from the one in \cite{gpv}: their arrow diagram  is obtained  from   ours by
reversing    all arrows  with  sign $-$.

There is a canonical mapping from
the set of isotopy classes of (oriented) knots in $\Sigma\times
\RR$ into the set of virtual knots.   It  assigns to a   knot    $K\subset \Sigma\times
\RR$ the virtual knot $[D(d)]$ where $d$  is a   knot diagram  on
$\Sigma$  presenting a knot   in  $\Sigma\times \RR$  isotopic to $K$.    The virtual knot  $[D(d)]$ does not depend  on 
  the choice  of $d$.  This follows from the fact that two  knot diagrams on
$\Sigma$  presenting isotopic knots  in  $\Sigma\times \RR$   can be obtained from each other  by ambient
isotopy in
$\Sigma$ and the     Reidemeister moves. Recall the standard list of the   Reidemeister
moves:  (1)  a  move adding a twist on the right (resp. left) of a branch; (2)  a move  pushing a
branch over another  branch and creating two crossings;  (3) a move pushing a branch over a crossing.  This list is
redundant. In particular,  the left move of type (1)    can be presented as a composition of  type (2) moves  and the inverse
to a right move of type (1).  One  move  of type (3)  together with moves of
type  (2)  is sufficient  to generate all moves of type (3) corresponding to various orientations on the branches  (see, for
instance, 
\cite{turr}, pp. 543--544).   As the generating move of type (3) we take the move  (c${)}^{-}$ described in Section
\ref{sn:g13}.   It remains   to observe that  the moves     (a$)_{ad}$, (b$)_{ad}$,
(c$)_{ad}$ on arrow diagrams are  exactly the moves induced by the right Reidemeister moves of type (1),  the  
 Reidemeister moves of type (2), and the move (c${)}^{-}$.
		     
 \subsection{Skein algebra of virtual knots.}\label{vks:3}    Let 
$R=\QQ[z]$ be the ring of polynomials in one variable $z$ with  rational coefficients.   Consider the polynomial  algebra $
R[\mathcal  K]$ generated by 
  the set of virtual knots $
\mathcal  K$.  This is a commutative associative algebra with  unit whose elements are polynomials in elements of $\mathcal  K$ with coefficients in
$R$.  
 We now introduce certain elements of $R[\mathcal  K]$ called   skein relations.

   Pick an arrow diagram $D$ with core circle $S$ and pick an arrow $e=(a,b)$ of $D$ with sign $+$ (here $a,b\in S$).  Let
$D^-_e$ be the same arrow diagram with the sign of $e$ changed to $-$.  Let $D'_e$ be the arrow diagram obtained from
$D$ by removing all arrows with at least one endpoint on the arc $ba\subset S$. Let $D''_e$ be the arrow diagram obtained
from
$D$ by removing all arrows with at least one endpoint on the arc $ab\subset S$.  The {\it skein relation} corresponding to
$(D,e)$ is
 $[D]-[D^-_e]-z [D'_e] [D''_e]\in R[\mathcal  K]$.

The ideal of  the algebra $ R[\mathcal  K]$ generated by the trivial virtual  knot  and the skein relations
(determined by all the  pairs
$(D,e)$ as above) is called the {\it skein ideal}.   The quotient of  
$ R[\mathcal  K]$  by this ideal   is called the {\it skein algebra of virtual knots} and denoted  $\mathcal {E}$. The  next
theorem computes
$\mathcal {E}$ in terms of strings. Recall the set 
 $\str_0$   of  non-trivial 
homotopy classes of virtual strings, cf. Section \ref{su71}.

\begin{theor}\label{th:eee}
                    There is a canonical $R$-algebra isomorphism $\nabla: \mathcal {E} \to R[\str_0]$ where $R[\str_0]$ is the
polynomial algebra  generated by 
 $
\str_0$.
                     \end{theor}

This  theorem allows us to associate with any virtual knot $K$ a polynomial $\nabla(K)\in R[\str_0]$.   It will be clear from
the definitions that
$$\nabla(K)= \langle \underline K \rangle+\sum_{n\geq 2} z^{n-1}  \,\nabla_n (K)$$
where $ \nabla_n (K) $ is a homogeneous element of  $\QQ [\str_0]$ of degree $n$ which is non-zero only for a finite set of
$n$.  Combining
$\nabla$ with homotopy invariants of strings we obtain  invariants of virtual knots.  For example, composing $\nabla$ with the
algebra homomorphism
$R[\str_0] \to R[t]$ sending the homotopy class of a string $\alpha$ into the polynomial $u(\alpha)(t) $, we obtain an
algebra homomorphism  $\mathcal {E} \to R[t]=\QQ[z,t]$. This gives a   2-variable polynomial invariant of virtual knots.   Further polynomial 
invariants of virtual knots can be similarly obtained  from the higher $u$-polynomials defined in Section \ref{covvvr1}.   

The constructions above can be applied to the virtual knot   derived from a geometric knot   $K\subset \Sigma\times
\RR$   in
Section
\ref{vks:2}. The resulting polynomial $\nabla(K)\in R[\str_0]$  is  invariant  under  the action on knots of orientation preserving 
 homeomorphisms  
$\Sigma\times
\RR\to \Sigma\times
\RR$
induced by orientation preserving homeomorphisms   $\Sigma\to \Sigma$.     The polynomial $\nabla(K)$ is interesting only in the case when the
genus of $\Sigma$  is at least
$2$.  This is due to the fact that the strings realized by curves on a surface  of genus 0 or 1 are homotopically trivial.

Theorem \ref{th:eee} will be proven in the next section.  Here we   give   an explicit expression for the value of  $\nabla$
on the generator  
$[D]\in \mathcal {E}$ represented by an arrow diagram
$D$. We need a few definitions. 
  The endpoints of  the
arrows of  $D$ split the core circle of $D$ into  (oriented) arcs called the {\it edges} of $D$.  Denote the set of  
edges of $D$ by
$\edg(D)$.   Each endpoint  $a$ of an arrow of $D$ is adjacent to two edges  $a_-,a_+\in \edg(D)$, respectively
incoming and   outgoing with respect to $a$. 
For an integer $n\geq 1$, an {\it $n$-labeling} of $D$ is a   map  $f:\edg(D)\to \{1,2,...,n\}$ satisfying the
following condition: for any  arrow $e=(a,b)$ of $D$,  either 

 (i) $f(a_+)=f(a_-), 
f(b_+)=f(b_-)$ or 

 (ii) 
$f(a_+)= f(b_-) \neq  f(a_-)= f(b_+)$ and  $\sign( f(a_-)-f(a_+))=\sign (e)$.  

The arrows $e$ as in (ii) are called {\it $f$-cutting arrows}.  The number of   $f$-cutting
arrows 
  of 
$D$ is denoted $\vert f\vert$ and the   number of $f$-cutting
arrows 
  of 
$D$ with $\sign =-1$  is denoted $\vert f\vert_-$.
Note that the value of $f$ on two adjacent 
edges
$a_-,a_+\in
\edg(D)$ may differ only when 
$a$ is  an endpoint of an  $f$-cutting arrow.  Therefore   $\vert f\vert \geq \# f(\edg
(D))-1$.  For $i=1,...,n$,  let ${\underline D}_{f,i}$   be the string  obtained from $D$ by removing all arrows except  the
arrows
$ (a,b)$ with $f(a_+)=f(a_-)= 
f(b_+)=f(b_-)=i$ (and forgetting the signs of the arrows).    

Let $\lbl_n (D)$ be the set of $n$-labelings $f$
of
$D$ such that   $f(\edg
(D))=\{1,...,n\}$,   $\vert f\vert=n-1$, and the  $f$-cutting arrows of $D$ are  pairwise unlinked
(in the sense of Section 
\ref{sn:g21}).  Then
\begin{equation}\label{etae}\nabla([D])=\sum_{n=1} \sum_{f\in  \lbl_n (D)} \frac {(-1)^{\vert f\vert_-} z^{n-1}}{n!}
\prod_{i=1}^n
\langle {\underline D}_{f,i}  
\rangle  \in R[\str_0].\end{equation}
The expression on the right-hand side is finite since $\lbl_n (D)=\emptyset$ for $n>\# \edg(D)$.  The set
 $\lbl_1 (D)$ consists of only one element $f=1$ so that  the   free term of
$\nabla([D])$  is   $\langle {\underline
D} 
\rangle$.

\section{Proof of Theorem \ref{th:eee}}\label{vksss}

The proof of Theorem \ref{th:eee}  largely follows the proof of Theorems 9.2 and 13.2 in \cite{tu}.  We therefore expose
only the main lines of the proof.  The key point  behind Theorem \ref{th:eee} is the existence of a natural comultiplication in
$\mathcal E$ and we define it first.  Then we   construct $\nabla$ and   prove that it is an isomorphism.

 \subsection{Comultiplication in $\mathcal E$.}\label{vks:4} We need   to
study more extensively the  labelings of   arrow diagrams defined at the end of the previous section.
Let $D$ be an arrow diagram with core circle $S$.  Each  $n$-labeling $f$ of   
$D$  gives rise to  $n$ 
monomials  $D_{f,1}, ..., D_{f,n} \in \mathcal
{E}$  as follows.    Identifying  $a=b$ for every $f$-cutting arrow  $(a,b)$ of $D$, we
transform  
$S$ into a 4-valent graph,  $\Gamma^f$, with $\vert f\vert$  vertices.   The  projection $S\to
\Gamma^f$  maps the non-$f$-cutting arrows of
$D$ into \lq\lq arrows" on  $\Gamma^f$, i.e.,  into ordered pairs  of (distinct)  generic points of $\Gamma^f$. The
labeling
$f$   induces a labeling of the edges of
$\Gamma^f$ by the numbers $1,2,...,n$.  It follows from the definition of  a labeling that for each $i=1,...,n$, the
union of edges of
$\Gamma^f$ labeled with
$i$ is  a disjoint union of $r_i=r_i(f)\geq 0$   circles   $S^i_1,..., S^i_{r_i}$. 
The orientation of   $S$ induces an
orientation of the edges of $\Gamma_f$ and of these circles.  We  transform each  
circle  $S^i_j$ with $j=1,..., r_i$  into an arrow diagram by adding to it all the arrows of $\Gamma^f$ with both endpoints
on
$S^i_j$.  The signs of these arrows are by definition the signs of the corresponding non-$f$-cutting arrows of $D$.   
Set $$D_{f,i}= \prod_{j=1}^{r_i} \, [S^i_j] \in \mathcal E.$$

For    any $n\geq  2$, denote  
$\Lbl_n  (D)$ the set of
$n$-labelings
$f$ of
$D$ such that  the $f$-cutting arrows of $D$ are  pairwise unlinked.  The latter condition can be reformulated in terms of
the numbers $r_1(f),...,r_n(f)$ introduced above: $f\in \Lbl_n  (D)$ if and only if $r_1(f)+ \cdots + r_n (f)= \vert f\vert
+1$.   For
$f\in \Lbl_n (D)$, set  
$$\Delta(D,f)= (-1)^{\vert f\vert_-} z^{\vert f\vert}\,  D_{f,1}  \otimes   D_{f,2}  \otimes \cdots  \otimes D_{f,n} \in
\mathcal {E}^{\otimes n}$$
where $\mathcal {E}^{\otimes
n}$ is the tensor product over $R$ of $n$ copies of $\mathcal {E}$.

  By a comultiplication in $\mathcal E$, we mean a coassociative algebra homomorphism $\Delta:\mathcal
{E}\to
\mathcal {E}{\otimes }\mathcal {E}$. (The coassociativity means that $(\id \otimes \Delta) \Delta= (\Delta \otimes \id
)\Delta$.)   We claim that the formula 
$$\Delta ([D])=   \sum_{f\in  \Lbl_2 (D)} \Delta(D,f)\in \mathcal {E}{\otimes
}\mathcal {E}$$
   extends  by multiplicativity to a well-defined  comultiplication    in $\mathcal E$. 
This can be deduced from     \cite{tu},  Theorem  9.2  or  proven directly repeating the same arguments.  We explain
how to deduce our claim from    \cite{tu}.  Comparing the  definition of 
$\Delta ([D])$ with the comultiplication in the algebra of skein classes of knots in (a surface)$\times \RR$ given in  \cite{tu},
we observe that they correspond to each other provided $D$  underlies a knot diagram on the surface.  (The variables
$h=h_1, \overline h=h_{-1}$ used in  \cite{tu} should be replaced with $0$ and $z$, respectively.  After the substitution 
$h=0$,  we can consider only labelings  satifying - in the notation of  \cite{tu} - the condition $\vert \vert f \vert \vert =-
\vert f\vert$ which  translates here as the assumption that   the  $f$-cutting arrows of $D$ are  pairwise unlinked.)   
The  results of 
 \cite{tu} imply that  if
  a  move 
(a$)_{ad}$, (b$)_{ad}$, (c$)_{ad}$  on
$D$ underlies a Reidemeister move on a knot diagram, then  $\Delta ([D])$ is preserved under  this move. Since any
arrow diagram $D$ underlies a knot diagram on a surface  and any move  (a$)_{ad}$, (b$)_{ad}$, (c$)_{ad}$  on
$D$ can be induced  by a  Reidemeister move, we conclude that  $\Delta ([D])$ is 
invariant under the moves  (a$)_{ad}$, (b$)_{ad}$, (c$)_{ad}$  on
$D$.  Therefore the formula $[D] \mapsto \Delta([D])$  yields a well-defined mapping $\mathcal K \to     \mathcal
{E}{\otimes }\mathcal {E}$. This mapping  uniquely extends    to an  algebra
homomorphism 
$ R[\mathcal
{K}] \to
\mathcal {E}{\otimes }\mathcal {E}$. 
The  results
of 
 \cite{tu} imply that  for  an arrow diagram   $D$  underlying a knot diagram  on a surface and  any arrow $e$ of $D$ with
$\sign (e)=+$, the   skein relation  
 $  [D]- [D^-_e]-z  [D'_e]  \,  [D''_e]  $ lies in the kernel of the latter homomorphism.   The condition that  
$D$  underlies a knot diagram    is verified for all $D$.   Therefore the  homomorphism $ R[\mathcal
{K}] \to
\mathcal {E}{\otimes }\mathcal {E}$ annihilates the skein ideal and induces  an  algebra
homomorphism 
$ \Delta: \mathcal
{E}\to
\mathcal {E}{\otimes }\mathcal {E}$.  
The coassociativity of $\Delta$  follows from the easy 
formulas
$$(\id \otimes \Delta) \Delta ([D])=  \sum_{f\in  \Lbl_3 (D)} \Delta(D,f)= (\Delta \otimes \id
)\Delta ([D])$$
 (cf.   \cite{tu}, p. 665).  More generally, 
   for any $n\geq 2$, the value  on $[D] \in \mathcal E$ of the iterated homomorphism 
$$\Delta^{(n)}= (\id^{\otimes (n-1)} \otimes \Delta)\circ  (\id^{\otimes (n-2)} \otimes \Delta) \circ ...\circ 
(\id  \otimes \Delta) \Delta: \mathcal {E}\to \mathcal {E}^{\otimes
 (n+1)}$$
  is computed by
$$\Delta^{(n)} ([D])= \sum_{f\in  \Lbl_{n+1} (D)} \Delta(D,f).$$
Note for the record
that   each arrow diagram $D$ admits   constant 2-labelings $f_1, f_2$   taking values $ 1,2$ on all edges,  respectively. 
The corresponding  summands of $\Delta ([D])$ are 
$\Delta(D,f_1)=[D] \otimes 1$ and $\Delta(D,f_2)= 1  \otimes [D]$.

  \subsection{Homomorphism $\nabla: \mathcal {E} \to R[\str_0]$.}\label{vks:6} There are two obvious $R$-linear
homomorphisms  $\varepsilon:  \mathcal {E}\to  R$ and  $ q:  \mathcal {E} \to R[\str_0]$.   The homomorphism
$\varepsilon$ sends
$1\in  \mathcal {E}$ into $1\in R$ and sends all   virtual knots and their  non-void products into $0$. 
 The homomorphism
$q$ sends
  $1$ and  all  products of $\geq 2$ virtual knots into
$0$  and sends  a virtual knot $K$ into $\langle \underline K \rangle$.    Tensorizing $q$ with itself,  we obtain for all 
$n\geq 1$ a homomorphism $q^{\otimes n}: \mathcal {E}^{\otimes n} \to R[\str_0]^{\otimes n}$.  Let
$s_n: R[\str_0]^{\otimes n}
\to  R[\str_0]$ be the $R$-linear homomorphism sending $a_1\otimes \cdots \otimes a_n$ into    $(n!)^{-1} a_1\cdots
a_n$. Set 
$$\nabla=\varepsilon + q+ \sum_{n\geq 2}  \, s_n \,q^{\otimes n} \Delta^{(n-1)}:  \mathcal {E} \to R[\str_0] $$
where  $\Delta^{(1)}=\Delta$.  It is clear that $\nabla$ is  $R$-linear. The same  argument as in   \cite{tu},  Lemma 13.4
shows that $\nabla$ is an algebra homomorphism.  Computing $\nabla$ on
the skein class of an arrow diagram $D$, we obtain 
$$\nabla([D])=\sum_{n\geq 1}   \, s_n \,q^{\otimes n}   \sum_{f\in  \Lbl_{n} (D)} \Delta(D,f)= \sum_{n\geq 1}       \sum_{f\in  \Lbl_{n} (D)} 
\frac {(-1)^{\vert f\vert_-} z^{\vert f\vert}}{n!} \,  \prod_{i=1}^n  q(D_{f,i}).$$
Note that  $q(D_{f,i})=0$ unless $r_i(f)=1$ in which case  $q(D_{f,i})=\langle \underline D_{f,i}\rangle$. 
For a labeling $f\in \Lbl_n(D)$ the equalities $r_1(f)=...=r_n(f)=1$ are equivalent to the inclusion $f\in \lbl_n(D)$.  
This  yields  Formula \ref{etae}.

Observe that   $\nabla([D])$ is a sum of   $\langle \underline
D\rangle$ and a polynomial  in strings  of rank $<\rank D$.  An  induction on the rank of strings shows that the image of
$\nabla$ contains all strings. Therefore $\nabla$ is surjective. 

 The proof of the injectivity of $\nabla$  is based on the
following   lemma. 

\begin{lemma}\label{l:jae} There is a $\QQ$-valued function  $\eta$ on the set of isomorphism classes of (finite) oriented
trees such that the following holds:

(i) if $T$ is  a tree with one vertex and no edges, then $\eta(T)=1$;

(ii) if an oriented tree $T'$ (resp. $U$) is obtained from an oriented tree $T$ by reversing the orientation of an edge
$e$ (resp. by contracting $e$ into a point), then $\eta(T)+\eta(T')+\eta(U)=0$;

(iii)   if an oriented tree $T'$ (resp. $T''$) is obtained from an oriented tree $T$ by replacing two distinct edges with
common origin $ab, ac$ by $ab, bc$ (resp. by $ac, cb$) and if $U$ is obtained from $T$ by identifying $b$ with $c$ and
$ab$ with $ac$, then $\eta(T)=\eta(T')+ \eta(T'')+\eta(U)$.
\end{lemma}

In this lemma  by an edge $ab$ we mean an {\it oriented} edge directed from $a$ to $b$.

Lemma \ref{l:jae} was first established in \cite{tu}, Theorem 14.1  where it is also shown that   $\eta$  is unique 
(we shall not need this).  The construction   in \cite{tu} is indirect and does not provide an explicit  formula for $\eta$.  Such
a formula was pointed out by Fran\c cois Jaeger  \cite{ja}.  The following proof of Lemma \ref{l:jae}  is a simplified version
of the proof given by   Jaeger 
\cite{ja}.

                     \begin{proof} By a {\it forest} we shall mean  a disjoint union of a finite family of  finite oriented trees.  The set
of vertices of a forest $F$ is denoted $V(F)$. For a forest
$F$ and an integer $n\geq 1$,  denote by  $C_n(F)$ the set of surjective mappings $f:V(F) \to \{1,...,n\}$ such that for
every edge $ab$ of $F$ we have $f(a)<f(b)$.  This set is empty for  $n> \#(V(F))$.  Set
$$\eta(F)= \sum_{n\geq 1} \frac {(-1)^{n+1}} {n}\, \# (C_n(F)) \in \QQ.$$
We claim that $\eta$ satisfies all the conditions of the lemma.  Condition (i) is obvious. Condition (iii) is a  direct  corollary of
the definitions. Indeed for all $n$,  the set $C_n(T)$ splits as a disjoint union of the sets $C_n(T'), C_n(T''), C_n(U)$. Hence
$\# (C_n(T))=\# (C_n(T'))+\# (C_n(T''))+\# (C_n(U))$ and   $\eta(T)=\eta(T')+ \eta(T'')+\eta(U)$.  It remains to verify 
(ii).   Let $F$ be obtained from $T$ by removing the interior of the edge $e$.  For all $n$,   
the set $C_n(F)$ splits as a disjoint union of the sets $C_n(T), C_n(T'), C_n(U)$. Hence 
$\# (C_n(F))=\# (C_n(T))+\# (C_n(T'))+\# (C_n(U))$ and  $\eta(F)=\eta(T)+ \eta(T')+\eta(U)$. Thus we need only to
prove that  $\eta(F)=0$ for every forest $F$ with two components $T_1,T_2$.

For  non-negative integers $n,k_1,k_2 $, denote by $C_n(k_1,k_2)$ the set of pairs $(l_1,l_2)$ where for
$i=1,2$, $l_i$ is an order-preserving
injection from $\{1,..., k_i\}$ into $\{1,...,n\}$ and   
$l_1(\{1,..., k_1\})\cup l_2(\{1,..., k_2\})= \{1,..., n\}$. Having $g_1\in C_{k_1} (T_1)$, $g_2\in C_{k_2} (T_2)$
and having  $(l_1,l_2)\in C_n (k_1,k_2)$ we   define a mapping  $f=f(g_1,g_2,l_1,l_2):V(F) \to \{1,...,n\}$ by
  $f(v)=l_1 g_1(v)$ for $v\in V(T_1)$ and  $f(v)=l_2 g_2(v)$ for $v\in V(T_2)$.
Clearly,  $f\in C_n(F)$.  It is obvious that any      $f\in C_n(F)$ can be uniquely presented in the form 
$f=f(g_1,g_2,l_1,l_2)$ where $g_i\in C_{k_i} (T_i)$ with $k_i=\#( f(V(T_i)))\geq 1 $ for $i=1,2$.  Therefore
$$\eta(F)=  \sum_{n\geq 1}   \frac
{(-1)^{n+1}} {n}\, \left ( \sum_{k_1,k_2 \geq 1} \,\sum_{g_1\in C_{k_1} (T_1),g_2\in C_{k_2} (T_2)}  \# (C_n
(k_1,k_2)) \right )$$
$$ =
\sum_{k_1,k_2 \geq 1} \# (C_{k_1} (T_1)) \, \# (C_{k_2} (T_2)) 
\sum_{n\geq 1}   \frac {(-1)^{n+1}} {n}\, \# (C_n (k_1,k_2)) .$$
Thus  it is enough  to prove that for all $k_1 \geq 1,k_2\geq 1$, the numbers $c_n (k_1,k_2)=\# (C_n (k_1,k_2))$ verify
\begin{equation}\label{fra} \sum_{n\geq 1}   \frac {(-1)^{n+1}} {n}\, c_n (k_1,k_2)=0.\end{equation}
Clearly, $c_n (k_1,k_2)$ is the number of pairs $(S_1,S_2)$ where $S_1,S_2$ are subsets of $\{1,...,n\}$ such that 
$S_1\cup S_2=\{1,...,n\}$, $\#(S_1)=k_1$, $\#(S_2)=k_2$. In particular, $c_n (k_1,k_2)=0$  if $k_1+k_2<n$ or  
$k_1>n$ or $k_2>n$.   For any
$n\geq 1$ and commuting variables
$x, y $, 
$$ (x+y+xy)^n=\sum_{  k_1,k_2\geq 0 } c_n (k_1,k_2)\, x^{k_1} y^{k_2}.$$
Therefore  
$$\log (1+x+y+xy)= \sum_{n\geq 1} \frac{ (-1)^{n+1} }{n} (x+y+xy)^n
$$
$$= \sum_{n\geq 1} \,\,\sum_{  k_1,k_2\geq 0 }  \frac {(-1)^{n+1}}{n} c_n (k_1,k_2) \, x^{k_1} y^{k_2}
= \sum_{  k_1,k_2\geq 0 } \,\,    \left (\sum_{n\geq 1} \frac {(-1)^{n+1}}{n} c_n (k_1,k_2) \right )\, x^{k_1} y^{k_2}.$$
Since
$$ \log (1+x+y+xy) =\log ((1+x)(1+y))=\log (1+x) +\log (1+y),$$
the terms with $k_1 \geq 1,k_2\geq 1$ in the above series must vanish.  This gives  Formula \ref{fra}. \end{proof} 

\subsection{The injectivity of $\nabla: \mathcal {E} \to R[\str_0]$.}\label{vks:624} We begin by associating with any
virtual string $\alpha$ an element $\zeta (\alpha)\in  \mathcal {E}$.  Let $S$ be the core circle of $\alpha$. 
A {\it  surgery} along an   arrow  $(a,b)\in \arr(\alpha)$ consists in picking two    (positively oriented) arcs    
$aa^+, bb^+\subset S$ and  then quotienting the  complement  of their interiors  $S-((aa^+)^\circ \cup (bb^+)^\circ)$ by
$a=b^+, b=a^+$.  It is understood that the arcs    
$aa^+, bb^+ $ are small enough not to    contain  endpoints of $\alpha$ besides $a,b$, respectively.  Such a surgery
transforms
$S$ into  two disjoint oriented circles.  We make each of them into a string by adding  all the arrows of
$\alpha$ with both endpoints  on the arc
$a^+b$ (resp. on $ba^+$).  (The arrows of $\alpha$ with one endpoint on $ab$ and the other one on $ba$
disappear under  surgery.)

Let us call a
set
$F\subset
\arr(\alpha)$ {\it special} if the arrows of $\alpha$ belonging to $F$ are pairwise unlinked.  
Applying surgery inductively to all arrows of
$\alpha$ belonging to a special set 
$F$, we transform
$\alpha$ into    $n=\#(F)+1$ strings.  Providing all the arrows of these   strings
with sign $+$, we obtain $n$ arrow diagrams    $D^F_1,..., D^F_n$.  Note that  
 they
  have together  at most     $\#(\arr(\alpha))-\#(F) $ arrows. We now  define an oriented graph
$\Gamma_F$. The vertices of $\Gamma_F$ are the symbols   $v_1,..., v_n$.  Two verices $v_i, v_j$
are related by an oriented edge   leading from $v_i$ to $v_j$ if there is an   arrow  $(a,b)\in F$ such that
the arcs $aa^+, bb^+\subset S$ involved in the surgery along this arrow lie  on the core circles of  $D^F_i, D^F_j$,
respectively. It is easy to see that $\Gamma_F$ is a tree. Set
$$\zeta(\alpha)=\sum_{F\subset \arr(\alpha)} \eta(\Gamma_F)\, z^{\#(F)} \prod_{i=1}^{\#(F)+1} [ D^F_i]\in \mathcal
{E}$$ where $F$ runs over all special subsets of $\arr(\alpha)$.  The summand corresponding to $F=\emptyset$ is the
string
$\alpha$ itself with sign  $+$ on all arrows.  
 
The key property of $\zeta(\alpha)  \in \mathcal
{E}$ is   its invariance under the   basic homotopy moves on $\alpha$. This follows from  \cite{tu}, Lemma 15.1.1
in the case where the moves are realized geometrically by   homotopy of a curve   realizing $\alpha$  on a surface.   Since
the homotopy moves can be always realized geometrically, the result follows. 
The mapping $\alpha\mapsto \zeta(\alpha)$  extends by multiplicativity to an algebra homomorphism $
R[\str_0]\to \mathcal
{E}$ denoted also  $\zeta$.

We can now prove the injectivity of $\nabla$. For $r\geq 0$, denote by $B_r$ the $R$-submodule of $\mathcal E$ 
additively generated  by  monomials $[D_1] [D_2] \cdots [D_n]$ such that the total number of arrows in the arrow
diagrams $D_1, D_2, ...,D_n$ is less than or equal to $r$.  Clearly, $0=B_0\subset B_1\subset ...$ and $\cup_r B_r=
\mathcal E$.  Pick $b= [D_1]
[D_2]
\cdots [D_n]\in B_r$.  Using the skein relation in $\mathcal E$ it is easy to see that    $b (\modu B_{r-1})\in
B_r/B_{r-1}$ does not depend on the signs of the arrows of $D_1,...,D_n$.  This observation, Formula \ref{etae}  and the
definition   of $\zeta $ imply  that 
$(\zeta
\nabla )(b) -b \in B_{r-1}$. Therefore $(\zeta \nabla-\id)^r(b) =0$. The inclusion $b\in \Ker \nabla$ would imply  
$b=0$. Thus  $B_r \cap \Ker \nabla  =0$.   Since $\cup_r B_r=
\mathcal E$, we obtain   $\Ker \nabla=0$. 

\subsection{More on $\mathcal E$.}\label{vks:624ddd}    The comultiplication   $\Delta$ defined  in Section \ref{vks:4}
makes
$\mathcal E$ into a Hopf algebra over $R$. Its counit is the    homomorphism  $\varepsilon:  \mathcal {E}\to R$ 
 used in the definition of  $\nabla$.   For an arrow diagram $D$, denote by
$\tilde D$ the same diagram with opposite signs  on all
arrows.  The transformation $[D] \mapsto -[\tilde D]$ preserves the skein relation  and therefore induces an algebra
automorphism of
$\mathcal E$.  This automorphism is an antipode for $\mathcal
E$. This follows from the corresponding theorem for the skein algebras of curves on surfaces conjectured in \cite{tu} and
proven  in \cite{crr}  and independently in 
\cite{pr}.     In the construction of the Hopf algebra
$\mathcal E$ instead of the ground ring
$R=\QQ[z]$ we can use  $\ZZ[z]$.  It is only to construct the homomorphisms $\nabla$ and  $\zeta$   that
we need    $\QQ$.

  Consider the Hopf algebra $ S_z ( \A_0 )$ derived as in Section \ref{su:6111} from the spiral Lie coalgebra $\A_0$, the
ring  
$R=\QQ[z]$ and the element $h=z\in R$. Note that  
$ S_h( \A_0 )= R[\str_0]$ as algebras.

\begin{theor}\label{th:fge} The 
           homomorphism  $\nabla: \mathcal {E} \to R[\str_0]=S_z ( \A_0 )$  is an isomorphism of Hopf algebras.
                     \end{theor}

The proof of this theorem follows  the lines of \cite{tu}, Section 12 and Lemma 13.5; we omit the details.

	\section{Open strings} 

\subsection{Definitions.}\label{openssd}   Replacing the circle  in the definition of a virtual  string  by an oriented
one-dimensional manifold $X$ we obtain a   {\it virtual string with core manifold} $X$.   The definition   of  homotopy
extends   to strings with core manifold  $X$ word for word.  Of special interest are strings with core manifold homeomorphic
to 
 $ [0,1]$;  we call them  {\it open strings}. In this context it is natural to call     virtual strings with core manifold homeomorphic
to 
 $S^1$
 {\it closed} strings.   

Open strings  underlie  (generic)
paths on   surfaces connecting  distinct  points on the boundary.  Gluing the endpoints of 
the core interval, we can transform any open string $\mu$ into a closed   string $\mu^{cl}$,   its {\it closure}.   Similarly
to paths, open strings  can be  multiplied   via the gluing of their core intervals along one   endpoint.  Repeating word for word the definitions of 
Section \ref{covvvr1} we obtain for all $r\geq 1$ the  notion of an  $r$-th covering  of an open string (this is again an open string).

\subsection{Polynomials of open strings.}\label{openss}   For an   open string $\mu$ with core manifold $[0,1]$,  we can
define two 
 polynomials
$u^+(\mu)$ and $u^-(\mu)$.    Observe that the set $\arr(\mu)$ of arrows of $\mu$ is a
disjoint union   $\arr^+(\mu)\cup \arr^-(\mu)$ where $\arr^{+} (\mu)$ (resp.  $\arr^{-} (\mu)$)  is the set of arrows
$ (a,b)\in \arr(\mu)$ with $a,b\in [0,1]$ such that $a<b$ (resp. $b<a$).   For $e \in \arr(\mu)$,   set $n(e)=n(e^{cl})\in
\ZZ$
 where  $e^{cl}$ is the corresponding arrow of   $\mu^{cl}$. 
 For   $k\geq 1$,  set  $$u_k^{\pm} (\mu)=\#\{e\in 
\arr^{\pm}(\mu)\,\vert \,n(e) =k\}
- \#\{e\in \arr^{\mp}(\mu)\,\vert \,n(e) =- k\} \in \ZZ.$$
This number and the polynomials  $u^{\pm} (\mu)=\sum_{k\geq 1} u^{\pm}_k(\mu)\, t^k$ are homotopy invariants of
$\mu$.   Clearly, $u(\mu^{cl})=u^{+} (\mu)+u^{-} (\mu)$ and $u^{\pm} (\mu \nu)=u^{\pm} (\mu)+u^{\pm} (\nu)$ for
any  open strings $\mu, \nu$.  Using  
$u^{\pm}$,  it is easy to give examples of non-homotopic open strings with homotopic closures. Using the coverings as in Section \ref{covvvr1},
we can define for open strings \lq\lq higher versions" of $u^{\pm}$ parametrized by finite sequences of positive integers.

 \subsection{Cobordism of  open  strings}\label{coboop} An open string $\mu$ is {\it slice} if its closure $\mu^{cl}$ is a
slice (closed) string.  Theorem \ref{th:103654} implies  that the sliceness is a homotopy property of open strings.

An open string $\mu$  is {\it ribbon} if its core interval admits an orientation reversing involution  transforming   arrows 
of $\mu$ into arrows of $\mu$ with opposite orientation. The closure of a ribbon open string is a ribbon
(closed) string. Therefore ribbon open strings are slice. 

We associate with   every open string $\mu$ an open string $\mu'$
obtained from $\mu$ by reversing orientation on the core interval and on all arrows. Clearly
$(\mu')'=\mu$ and $ (\mu\nu)'=\nu'\mu'$ for  any open strings $\mu,\nu$.  

We say that   open strings $\mu, \nu$ are   {\it cobordant} and write $\mu \sim_c \nu$ if
$\mu\nu'$ is slice.    An open  string    is cobordant to a trivial open string  (with no
arrows) if and only if  it   is slice.  

	\begin{lemma}\label{coboooo} (i) Cobordism is an equivalence relation on the set of  open strings.

(ii) Homotopic open strings are cobordant. 

(iii) If two open strings are cobordant, then their closures are cobordant. 

(iii) If two open strings are cobordant, then their $r$-th coverings are cobordant for all $r\geq 1$.

			   \end{lemma} 
                     \begin{proof}  For any open string  $\mu$, the  product 
$\mu\mu' $ is ribbon and therefore slice.  Thus $\mu \sim_c \mu$. If  $\mu \sim_c \nu$, then 
$(\mu\nu')^{cl}$ is slice. The closed string  $ (\nu\mu')^{cl}$ is obtained from $(\mu\nu')^{cl}=((\nu\mu')')^{cl}$ by the
involution
$\alpha\mapsto \overline \alpha^-$. Hence $ (\nu\mu')^{cl}$ is slice and $\nu \sim_c \mu$.

To proceed we need the following property:  if $\mu, \nu,\delta $ are open strings whose product $\mu
\nu \delta $ is slice and if $\nu$ is slice, then so is $\mu \delta$.  Indeed, observe  that    $(\mu
\nu \delta )^{cl } $ is homeomorphic to   $(\nu \delta \mu)^{cl }$.   Since $(\nu \delta
\mu)^{cl}$ is a product of $\nu^{cl}$ and  $(  \delta
\mu)^{cl}$, the cancellation property mentioned in Section \ref{fbm}  implies   that  $(  \delta
\mu)^{cl}$ is slice.  Since  $(  \delta
\mu)^{cl}$ is homeomorphic to $( \mu \delta
)^{cl}$, the latter  string is slice. Hence 
$\mu \delta$ is   slice.

We can now prove the transitivity of cobordism.  
If  $\mu
\sim_c
\nu$,
$\nu
\sim_c
\delta$, then  $\mu \nu'$ and  $\nu \delta'$ are slice.  Since   products of slice closed strings are slice,  
the products of slice open strings are slice. Thus,   $\mu \nu' \nu \delta'$ is slice. Since $\nu'\nu$ is slice,  so is
$\mu\delta'$. Therefore $\mu \sim_c \delta$. 

If  $\mu$ is homotopic to $\nu$, then $(\mu \nu')^{cl}$ is homotopic to $(\nu \nu')^{cl}$. Therefore 
$(\mu \nu')^{cl}$ is slice so that $\mu\sim_c \nu$.  This proves (ii).  We leave (iii) and (iv) as an exercise for the reader.
\end{proof}

Multiplication of open strings induces a  multiplication in the set of cobordism
classes of open strings that makes this set into  a group denoted $\mathcal O$. Using the formulas $u^{\pm} (\mu')= - u^{\pm}
(\mu)$ and Theorem \ref{th:0388} we obtain that
$u(\mu^{cl})=u^+(\mu)+u^-(\mu)\in \ZZ[t]$ is an additive cobordism invariant of  open strings.

 \subsection{Graded based matrices}\label{fi:g5512}     A {\it graded (skew-symmetric) based matrix} over an abelian group
$H$
  is a  based   (skew-symmetric) matrix $(G, s, b)$
over $H$ endowed with a spltting  of 
$G-\{s\}$ as a union of  two  disjoint subsets $G^+$ and
$G^-$.  The {\it underlying based matrix}  of a graded based matrix
$(G,s,b)$ is   obtained   by forgetting the splitting
$G-\{s\} =G^+\cup G^-$.   The {\it negation} $-T$
of  $T=(G,s,b) $ is the triple $(G,s,-b)$ with the same splitting   $G-\{s\}=G^+\cup G^-$.   

We define annihilating elements, core elements, and complementary elements of  a graded based matrix $T=(G,s,b)$  
as in  Section
\ref{fi:g51} with the following additional requirements: annihilating elements must  lie in $G^+$, core elements  must lie in
$G^-$ and for any pair of complementary elements, one of them lies in $G^+$ and the second one in $G^-$.  All other
dfefinitions  and results of  Section \ref{fi:g51} extend to this setting with obvious changes.  In particular, we have a notion
of homology for graded based matrices over 
$H$.

  If
$H\subset
\RR$, then the two  1-variable polynomials
  $$ u^\pm(T) \,(t)=\sum_{g\in G^\pm, b(g,s)  > 0}  t^{ b(g,s) } -\sum_{g\in G^\mp, b(g,s)  < 0}  t^{ - b(g,s) } $$
are homology invariants of
a  graded based matrix
$T=(G,s,b)$.

For an  open string
$\mu$, the based matrix of  its closure $\mu^{cl}$  is   graded via  the splitting
$\arr(\mu^{cl}) = \arr(\mu) = \arr^+(\mu) \cup \arr^- (\mu)$.  This defines a graded based matrix   $T(\mu)$ over
$\ZZ$. The homology class of  $T(\mu)$  is an invariant of
the homotopy class of $\mu$.  Clearly,  $u^\pm (\mu)= u^\pm (T(\mu))$.  Note also that $T(\mu')=-T(\mu)$.

 \subsection{Addition of graded based matrices}\label{fi33330152}  
We define   addition for graded based matrices  which mimics the   product of    strings.  Let $T_i=(G_i, s_i, b_i)$ be
a graded based matrix over  an abelian group $H$ where $i=1,2$.  We define the  sum $T_1\oplus T_2=(G,s,b)$ as follows.  
Set  $G^{\pm}= G_1^{\pm} \amalg G_2^{\pm}$ and 
$G=\{s\}
\amalg G^+\amalg G^-$.  For    $g\in G-\{s\}=G^{+} \amalg G^{-}$,
set
$\varepsilon_g =1$ if $g\in G^+$ and $\varepsilon_g =0$ if $g\in G^-$.  The skew-symmetric mapping $b:G\times G \to
H$ is   defined  as follows:    for   
$g \in G_i-\{s_i\} $ with
$i \in \{1,2\}$,    set  $b(s,g)=b_i(s_i,g), b(g,s)=b_i(g,s_i) $;   for  any 
$g \in G_i-\{s_i\}, h \in G_j-\{s_j\}$ with
$i,j\in \{1,2\}$,   set 
$$ b(g,h) =\left \{ \begin {array} {ll} b_i(g,h),~ {\rm 
{if}}\,\,\, 
i=j, \\ 
\varepsilon_g b_j(s_j,h)- \varepsilon_h b_i(s_i,g) ,~ {\rm {if}} \,\,\, i\neq j. 
\end {array}  
\right.
$$
The direct sum of
graded based matrices is   commutative and associative (up to isomorphism).

Let $R$ be a domain.  A graded based matrix  over $R$  is {\it hyperbolic} if its underlying based matrix is hyperbolic. 

\begin{lemma}\label{lujuj1} For any graded based matrix  over   $R$, its direct sum  with its negation  is hyperbolic.
The direct sum of two hyperbolic graded based matrices over $R$  is hyperbolic. 
		   \end{lemma}
                     \begin{proof} Let $T_1=(G_1,s_1,b_1)$ be a graded based matrix over $R$. Let  $T_2=(G_2,s_2,b_2)$ be
a  copy of
$T_1$ where 
$G_2=\{g' \,\vert \, g\in G_1\}$, $s_2=(s_1)'$,    
and $b_2$ is defined by $b_2(g',h')=b(g,h)$ for $g,h\in G_1$.  We verify that  the direct sum $T=T_1 \oplus
(-T_2)=(G,s,b)$  is hyperbolic.  Consider the  subsets 
  $\{s\}$ and  
$\{g,  g'\}_{g\in G_1-\{s_1\}}$ of $G$.  These subsets form a simple filling of    $G$.  The matrix of this filling  is zero.
Indeed, for   $g \in G_1-\{s_1\}$,
$$b(s, \{g,g'\})=b(s, g )+b(s,  g')= b_1(s_1,g)+ (-b_2) (s_2,g')= b_1(s_1,g)- b_2((s_1)',g')=0.$$
For   $g,h\in G_1-\{s_1\}$, 
$$b(\{g,g'\},\{h,h'\})=b(g,h)+b(g,h')+b(g',h) +b(g',h')$$
$$=b_1(g,h)+ \varepsilon_g (-b_2)(s_2,h')- \varepsilon_{h'} b_1(s_1,g)
 + \varepsilon_{g'}  b_1 (s_1,h)- \varepsilon_{h} (-b_2) (s_2,g') +(-b_2)(g',h')$$
$$=b_1(g,h) - \varepsilon_g b_1(s_1,h)- \varepsilon_{h} b_1(s_1,g)
 + \varepsilon_{g} b_1 (s_1,h)+ \varepsilon_{h} b_1(s_1,g) -b_1(g,h)=0.$$
The second claim of the lemma is an exercise on the definitions; we leave it to the reader. \end{proof}

Quotienting the monoid of graded based matrices  over   $R$ by hyperbolic matrices, we obtain an abelian group 
$\mathcal G (R)$.  We call it the {\it group of cobordisms} of  graded based matrices  over  $R$.  

Assigning to an open string
its graded based matrix we obtain an additive homomorphism $\mathcal O\to \mathcal G(\ZZ)$. The group  $\mathcal
G(\ZZ)$ is non-trivial. This is clear from the  existence of   non-trivial additive homomorphisms  $u^\pm: \mathcal
G(\ZZ)\to \ZZ[t]$. 

The notion of a graded based matrix over $R$ and the addition of such matrices may seem   artificial from the algebraic viewpoint.
Possibly, a more satisfactory (although equivalent) language would   describe  a graded based matrix over $R$ as    a   free
$R$-module  $V$ of finite rank endowed with a vector  in the dual module $V^*=\Hom_R (V,R)$,  with a distinguished basis    partitioned into two
disjoint subsets, and with a
$R$-valued skew-symmetric bilinear form $V\times V \to R$. 
To pass from the  definition above to this one, we associate with $(G, s, b)$ the free $R$-module $V$ with basis $G-\{s\}$, 
the  element of $V^* $ sending any   $g\in 
G-\{s\}$ to $b(s,g)$, the partition $G-\{s\}=G^+\cup G^-$, and the 
skew-symmetric bilinear form $V\times V \to R$ induced by $b$.

 \subsection{The algebra of   open strings.}\label{a21aee}  Let  $R$ be a 
commutative ring 
  with unit and $\otimes=\otimes_R$.
  A (left) {\it module}  over a Lie algebra
$(L, [ \,,]:L^{\otimes 2} \to L)$ over
$R$ is  an
$R$-module
$M$ endowed with an
$R$-linear homomorphism $\rho:L\otimes M \to M$ such that \begin{equation} \label{drf} \rho ([\,,]\otimes \id_M)= \rho
(\id_L
\otimes
\rho )(\id_{L\otimes L \otimes M}- 
 \Perm_L\otimes \id_M) :L\otimes L \otimes M \to M\end{equation}
where $\Perm_L$ is  the permutation 
$x\otimes 
y\mapsto y\otimes x$ in $L^{\otimes 2} $.  (Formula \ref{drf} is equivalent to the usual identity
$[x,y]m=x(ym)-y(xm)$ for
$x,y\in L, m\in M$.)  Dually,  a  {\it comodule}  over
a Lie coalgebra 
$(A,\nu:A\to A^{\otimes 2})$ over
$R$ is  an
$R$-module
$M$ endowed with an
$R$-linear homomorphism $\rho:M\to A\otimes M  $ such that 
\begin{equation} \label{drref} (\nu \otimes \id_M) \rho =  (\id_{A\otimes A \otimes M}-
\Perm_A\otimes \id_M) (\id_A
\otimes
\rho)  \rho:M\to A\otimes A \otimes M.\end{equation} 
Such
$M$  is automatically a module over the dual Lie algebra
$A^*=\Hom_R(M,R)$: an element  $a\in A^*$ acts on $M$  by the endomorphism $\varphi_a:M\to M$ sending  $m\in M$ to 
$ - (a
\otimes
\id_{M})
\rho(m)\in R\otimes M=M$. 

 A comodule
$ (M,\rho)$ over a Lie coalgebra $A$  is {\it spriral} if $M=\cup_{n\geq 1} \Ker \rho^{(n)}$
where 
$$\rho^{(n)}=(\id_A^{\otimes (n-1)} \otimes \rho) \circ \cdots \circ 
  (\id_A  \otimes \rho) \circ \rho:M \to A^{\otimes n}\otimes M.$$  If  $R\supset \QQ$  and  
$A, M$ are 
spiral, then  
 the  action of $A^*$ on $M$ integrates into a group action of the group 
$\Exp  A^*$ on  $M$: an element  $a\in \Exp  A^*=A^*$ acts on $m\in M$ by   $$a m=e^{\varphi_a} (m)=m+\sum_{k\geq
1}   (\varphi_a)^k(m)/k!. $$   Note that for   $m\in \Ker \rho^{(n)}$ the sum on the
right-hand side has at most $n$ non-zero terms.

  Let
$\M=\M (R)$ be the free 
$R$-module freely generated by the 
set of 
homotopy classes of open virtual strings.  We    provide $\M$ with the 
structure of a comodule over the  
  Lie coalgebra  of closed strings $\A_0$.  Let $\langle \mu \rangle$ be the generator of $\M$ represented by an open 
string
$\mu$. For an arrow $e\in \arr(\mu)$,
  a surgery along $e$ defined as in Section \ref{vks:624} transforms $\mu$ into a disjoint union of a  closed  string
$\alpha_e$ and an open string $\beta_e$.  Set  
$$\rho (\langle \mu \rangle) =\sum_{e\in \arr_+(\mu)}   \langle \alpha_e \rangle \otimes \langle \beta_e
\rangle  - \sum_{e\in \arr_-(\mu)}   \langle \alpha_e \rangle \otimes \langle \beta_e
\rangle  \in \A_0 \otimes \M.$$
A direct computation shows that this  gives a well-defined  $R$-linear homomorphism $\rho:\M\to \A_0\otimes \M 
$ satisfying Formula \ref{drref}.  Thus $\M$ is a   comodule over  $\A_0$.   
Combining $\rho$ with the inclusion $\A_0 \subset \A$ we obtain that $\M$  is a
comodule over   $\A$ as well.   It is easy to see that $\M$ is spiral.   If  $R\supset \QQ$, then the construction
 above gives a group action of
$  \Exp  \A^*$ on  $\M$.  

\subsection{Exercises.}\label{a21ass}  1. Multiplication of open strings makes $\M$ into an associative algebra
with unit. Check  that  the group $\Exp   \A^*$ acts on  $\M$ by algebra automorphisms. 

2.   Let $cl:\M\to \A$  be the   $R$-linear homomorphism induced by closing open strings.  Check that for any open
string
$\mu$,
 we have $\nu(\langle cl(\mu)  \rangle)=(\id_{\A\otimes \A}  -\Perm_{\A} )  (\id_\A \otimes cl) \rho
(\langle \mu \rangle)  $.

	\section{Questions} \label{questions}

    1.    Which  primitive based matrices  $T_\bullet$ can  be realized as $T_\bullet (\alpha)$ for a string $\alpha$ ?     A
necessary  condition 
pointed out in Section \ref{sn:g22}  says that $(u (T_\bullet ))'(1)=0$. 
  Note that for the based matrix 
$T(\alpha)=(G,s,b)$, we have $\vert b(e,f)\vert \leq \# 
(G)-2$ for 
all $e,f \in G$.  
     This   however yields no conditions on the primitive based matrices of  strings, since such a matrix 
$T_\bullet=(G_\bullet,s_\bullet,b_\bullet)$ may  arise from a 
string of a rank  $\gg \# (G_\bullet)$.

2. Can one detect  non-slice strings with hyperbolic based matrices   using the secondary obstructions of Section \ref{1hhhgp8} ?

3. Is it true that slice   strings are stably ribbon, i.e., that for any slice string $\alpha$ there is a ribbon
string $\beta$ such that a product of  $\alpha$ and $\beta$ is homotopic to a ribbon  string ?  Is it true for open strings ? 
A positive answer to the second question would imply a positive answer to the first question.

4.   Classify all strings of small rank (say, $\leq
6$) up to homotopy and/or up to cobordism.

5. Is it true that every  string is homotopic to  a string  of type $\alpha_\sigma$ for  some permutation $\sigma$ ?
If not,  is  it      true up to cobordism ?

6. Is  multiplication of open strings    commutative  up to   homotopy ?  If not, is it commutative     up to cobordism ? 

7.   Compute the group  $\mathcal O$ of cobordism classes of  open strings. 

8. Compute the group $\mathcal G(\ZZ)$ of cobordism classes of  graded based matrices over $\ZZ$.

9.    Is there an  invariant of   virtual  knots combining the skein invariant
$\nabla$   with the Kontsevich universal finite type invariant of knots ?  This might lead to  mixed arrow-chord
diagrams. 

10.  Study   invariants of virtual strings that change 
in a 
controlled way (say by constants) under the moves (a$)_s$, (b$)_s$,  (c$)_s$, 
cf. the 
theory of Arnold's invariants of plane curves.  

11. Generalize the invariants of virtual knots introduced in this paper to virtual links.

                     \end{document}